\documentclass{article}

\usepackage[final]{neurips_2025}

\usepackage[utf8]{inputenc} 
\usepackage[T1]{fontenc}    
\usepackage{hyperref}       
\usepackage{url}            
\usepackage{booktabs}       
\usepackage{amsfonts}       
\usepackage{nicefrac}       
\usepackage{microtype}      
\usepackage{xcolor}         

\usepackage{amsmath}
\usepackage{bm}
\usepackage{algorithm}
\usepackage{algorithmic}
\usepackage{subcaption}
\usepackage{graphicx}
\usepackage{multirow}
\usepackage{makecell}
\usepackage{enumitem}

\title{A Gradient Guided Diffusion Framework for Chance Constrained Programming}

\author{%
	Boyang Zhang \\
	School of Advanced Interdisciplinary Sciences\\
	University of Chinese Academy of Sciences\\
	Beijing 100049, China \\
	\texttt{zhangboyang23@mails.ucas.ac.cn} \\
	\And
	Zhiguo Wang\textsuperscript{$*$} \\
	Department of Mathematics\\
	Sichuan University \\
	Chengdu 610065, China \\
	\texttt{wangzhiguo@scu.edu.cn} \\
	\AND
	Ya-Feng Liu\thanks{Corresponding authors.} \\
	Ministry of Education Key Laboratory of Mathematics and Information Networks\\
	School of Mathematical Sciences \\
	Beijing University of Posts and Telecommunications \\
	Beijing 102206, China \\
	\texttt{yafengliu@bupt.edu.cn} \\
}

\begin{document}

	\maketitle

	\begin{abstract}
		
		Chance constrained programming (CCP) is a powerful framework for addressing optimization problems under uncertainty. In this paper, we introduce a novel \textbf{G}radient-\textbf{G}uided \textbf{D}iffusion-based \textbf{Opt}imization framework, termed GGDOpt, which tackles CCP through three key innovations. First, GGDOpt accommodates a broad class of CCP problems without requiring the knowledge of the exact distribution of uncertainty—relying solely on a set of samples. Second, to address the nonconvexity of the chance constraints, it reformulates the CCP as a sampling problem over the product of two distributions: an unknown data distribution supported on a nonconvex set and a Boltzmann distribution defined by the objective function, which fully leverages both first- and second-order gradient information. Third, GGDOpt has theoretical convergence guarantees and provides practical error bounds under mild assumptions. By progressively injecting noise during the forward diffusion process to convexify the nonconvex feasible region, GGDOpt enables guided reverse sampling to generate asymptotically optimal solutions. Experimental results on synthetic datasets and a waveform design task in wireless communications demonstrate that GGDOpt outperforms existing methods in both solution quality and stability with nearly 80\% overhead reduction.
		
		Our code is available at https://github.com/boyangzhang2000/GGDOpt.

	\end{abstract}
	
	\section{Introduction}
	
	\subsection{Problem formulation}
	
	Chance constrained programming (CCP) is an efficient modeling paradigm for optimization problems with uncertain constraints, which finds wide applications in diverse fields, such as finance (\cite{bonami2009exact}), robot control (\cite{calafiore2006scenario}), and wireless communications (\cite{wang2014outage}). In this paper, we consider a CCP with the following form:
	\begin{equation}
		\label{eq1}
		\begin{aligned}
			\min_{\bm{x}}\quad&f(\bm{x})\\
			\mathrm{s.t.}\quad&\bm{x}\in\mathcal{X}_{\rho},
		\end{aligned}
	\end{equation}
	where $f:\mathbb{R}^n\rightarrow\mathbb{R}$ is a differentiable objective function and 
	$\mathcal{X}_{\rho}$ is the chance (or probabilistic) constraint set defined by
	\begin{equation}
		\label{eq:cons}
		\mathcal{X}_{\rho} = \Big\{\bm{x}\in \mathbb{R}^n \mid \text{Prob}_{\bm{h}}\{\bm{g}(\bm{x},\bm{h})\geq\bm{0}\}\geq1-\rho\Big\}. 
	\end{equation}
	In the above, $\bm{h}$ is a random vector with probability distribution $P$ supported on a set $\Xi\subset \mathbb{R}^d$, $\rho\in (0,1)$, $\bm{g}=(g_1,g_2,\ldots,g_m):\mathbb{R}^n\times\Xi\rightarrow\mathbb{R}^m$, and $\text{Prob}(A)$ denotes the probability of an event $A$. Problem (\ref{eq1}) is generally challenging to solve for the following two reasons. First, evaluating the probability term $\text{Prob}_{\bm{h}}\{\bm{g}(\bm{x},\bm{h})\geq0\}$ typically involves a high-dimensional integration, which is computationally intractable. Second, even when $\bm{g}$ is linear, the feasible set $\mathcal{X}_{\rho}$ remains nonconvex, further complicating the optimization. 
	
	\subsection{Related works}
	
	Apart from very special cases where $\mathcal{X}_{\rho}$ can be transformed into a convex formulation under strong assumptions (\cite{kataoka1963stochastic}, \cite{lagoa2005probabilistically}, \cite{henrion2007structural}, \cite{prekopa2013stochastic}), there are two popular approaches to tackling general problem (\ref{eq1}), which are Convex Approximation (CA) method and Sample Average Approximation (SAA) method. The CA method seeks to construct a tractable inner approximation of $\mathcal{X}_{\rho}$, but it typically requires the information of the \emph{exact} distribution $P$, often assuming that $P$ belongs to specific families such as Gaussian or log-concave distributions (\cite{ben2000robust}, \cite{bertsimas2004price}, \cite{lagoa2005probabilistically}, \cite{nemirovski2007convex}). In contrast, the SAA method approximates $P$ using an empirical distribution based on sampled data, reformulating the CCP as a binary integer program (\cite{ahmed2008solving}, \cite{pagnoncelli2009sample}, \cite{adam2016nonlinear}). However, this reformulation remains computationally intractable. These restrictive assumptions on the underlying distribution $P$, along with the high computational cost, significantly limit the practical applicability of CCP.
	
	One important question to ask is: \textbf{can we design a general framework to efficiently solve CCP when the underlying distribution $P$ is unknown?} The answer to the above question is particularly crucial in our interested case where samples can be efficiently drawn from $\mathcal{X}_{\rho}$, albeit the explicit formulation of $\mathcal{X}_{\rho}$ is unavailable. This motivates us to seek high-quality solutions to the CCP problem (\ref{eq1}) from a new perspective via sampling-based methods (\cite{wibisono2018sampling}, \cite{ma2019sampling}, \cite{lee2021structured}, \cite{chen2022improved}, \cite{seyoum2025beyond}). The core idea of applying sampling-based methods to solve CCP problems lies in reformulating the original nonconvex CCP with intractable constraints as a sampling problem from an unknown distribution. This reformulation leverages probabilistic techniques to handle the challenging constraints through stochastic sampling rather than deterministic evaluation.
	
	Notably, generative models are designed to approximate unknown data distributions based on observed samples, enabling the generation of new data points from the learned approximation. In particular, diffusion models have emerged as a powerful family of generative models, offering high-quality sample generation, stable training dynamics, and scalability to high-dimensional problems (\cite{ho2020denoising}). The sampling process based on score estimation enables diffusion models to generalize to conditional distributions, thereby generating samples that satisfy requirements through conditional information guidance (\cite{ho2022classifier}). As a powerful generative artificial intelligence (AI) technology, diffusion model has been successfully deployed across various domains, such as, image generation ( \cite{yue2023resshift}, \cite{huang2025diffusion}), inverse problems (\cite{chung2022improving}, \cite{chung2022diffusion}, \cite{songsolving}), and optimization (\cite{krishnamoorthy2023diffusion}, \cite{li2024diffusion}, \cite{wu2024diff}, \cite{kong2024diffusion},  \cite{liang2025diffusion}). Recently, \cite{guogradient} introduced a novel form of gradient guidance to adapt pre-trained diffusion models for user-specified tasks.
	
	Despite their success in various domains, diffusion models have rarely been explored in the context of CCP. The possible reason behind might be that tackling CCP problems via diffusion models generally requires efficient sampling from a composite distribution, the product of an unknown data distribution (associated with the constraint) and a known Boltzmann distribution (induced by the objective function), but the training data is only available from the unknown component. This makes the application of diffusion models to CCP both novel and nontrivial.

	\subsection{Our contributions}
	
\begin{figure}[t]
		\centering
		\begin{subfigure}[b]{0.98\textwidth}
			\centering
			\includegraphics[width=\linewidth]{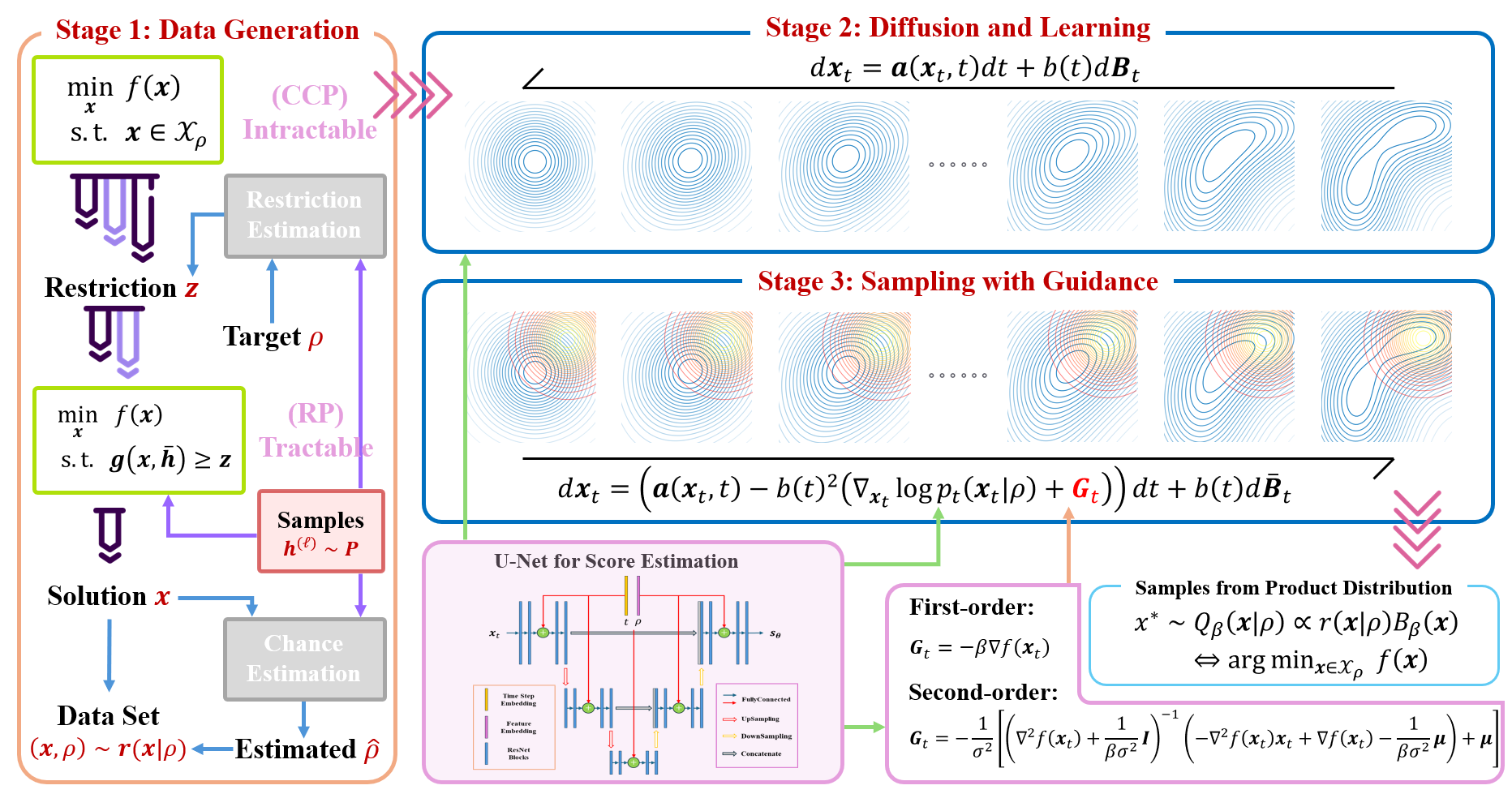}
		\end{subfigure}
		\caption{A framework of GGDOpt. (1) Generate a training set of points satisfying the chance constraint by solving a deterministic restricted problems. (2) Train a diffusion model with classifier-free guidance to learn the score of the conditional distribution. (3) Perform the reverse diffusion process with additional gradient guidance to sample from the product of the data distribution and the Boltzmann distribution.}
		\label{fig:ggdopt}
	\end{figure}

	In this paper, we propose GGDOpt (see Figure \ref{fig:ggdopt}), a novel \textbf{G}radient-\textbf{G}uided \textbf{D}iffusion-based \textbf{Opt}imization framework for solving problem (\ref{eq1}), with the following originality:
	\begin{itemize}[leftmargin=10pt]
		\item \textbf{Applicable to broader problem domains.} Built on the basis of diffusion model with classifier-free guidance and optimization via sampling, GGDOpt accommodates a broad class of CCP problems without requiring the knowledge of the exact distribution of uncertainty—relying solely on a set of samples.
		\item \textbf{Problem reformulation with a novel paradigm.} GGDOpt reformulates the CCP problem as a sampling task over the product of two distributions: an unknown data distribution implicitly defined by the constraint and a Boltzmann distribution induced by the objective function with a full utilization of first- and second-order information of the underlying CCP.
		\item \textbf{Feasibility-aware data generation and efficient guided sampling.} To generate high-quality training data that satisfy the chance constraint, GGDOpt solves a deterministic restricted problems by standard optimization techniques. The solutions are used to guide the training of the conditional diffusion model, effectively capturing the geometry of the feasible region. To sample from the product distribution, we develop a gradient-guided reverse process derived in closed form based on the structure of the product distribution. Compared with \cite{guogradient}, our guidance terms do not require backpropagation through the neural network.
		\item \textbf{Theoretical convergence and practical evaluation.} Regarding the sampling process as a reverse time stochastic differential equation (SDE), GGDOpt is shown to generate asymptotically optimal solutions as the time step and inverse temperature go to infinity. A practical error bound is also provided with two components: the limited time length error and limited inverse temperature error.
	\end{itemize}
	
	\subsection{Organization}
	
	The remainder of the paper is organized as follows. In Section \ref{sec:prob}, a reformulation of CCP problem (\ref{eq1}) is provided via sampling, and a gradient guidance-based score estimation schedule is provided with both first- and second-order information. A novel GGDOpt framework for solving problem (\ref{eq1}) is given in Section \ref{sec:GGDOpt}. Theoretical convergence and experimental results are presented in Section \ref{sec:conv} and Section \ref{sec:expe}, respectively. The conclusion is drawn in Section \ref{sec:conc}.
	
	\section{Problem reformulation via sampling}
	\label{sec:prob}
	Let $r(\bm{x}|\rho) = \mathbb{I}_{\mathcal{X}_{\rho}}(\bm{x})$ denote the indicator function of the chance constraint $\mathcal{X}_{\rho}$. Let $B_{\beta}(\bm{x})\propto e^{-\beta f(\bm{x})}$ represent the Boltzmann distribution associated with the objective function $f(\bm{x})$, where $\beta>0$. The resulting sampling task is to draw samples from the following target distribution:
	\begin{equation}
		\text{sample } \bm{x} \sim Q_{\beta}(\bm{x}|\rho)\propto r(\bm{x}|\rho)B_{\beta}(\bm{x}).
	\end{equation}
	
	Intuitively, the distribution $Q_{\beta}(\bm{x}|\rho)$ assigns higher probability density to regions where the objective function $f(\bm{x})$ takes smaller values. Under certain regularity conditions (\cite{kong2024diffusion}), as $\beta \rightarrow \infty$, the sampling distribution $Q_{\beta}(\bm{x}|\rho)$  asymptotically concentrates around the global minimizer of the CCP in (\ref{eq1}). Therefore, the CCP (\ref{eq1}) admits the following equivalent reformulation:
	\begin{equation}
		\displaystyle
		\bm{x}^*=\mathop{\arg\min}_{\bm{x}}\ \left\{f(\bm{x})+\mathbb{I}_{\mathcal{X}_{\rho}}(\bm{x})\right\}\iff \text{sample } \bm{x}^*\sim Q_{\beta}(\bm{x}|\rho),\ \beta\to\infty.
	\end{equation}	
	A natural way would be to directly employ Langevin dynamics for sampling from distribution $Q_{\beta}(\bm{x}|\rho)$. However, the unknown nature of component $r(\bm{x}|\rho)$ prevents the derivation of an exact expression of the score function. Fortunately, we can obtain a set of feasible samples $\{\bm{x}^{(i)},\rho^{(i)}\}_{i=1}^N$, which are drawn from the unknown distribution $r(\bm{x}|\rho)$. More details on this will be presented in Subsection \ref{subsec:data}. This motivates us to leverage diffusion models to directly learn the product distribution $Q_{\beta}(\bm{x}|\rho)\propto r(\bm{x}|\rho)B_{\beta}(\bm{x})$, where $r(\bm{x}|\rho)$ is unknown but $B_{\beta}(\bm{x})$ is explicitly known.

	\subsection{Diffusion models}
	
	Given observed samples $\bm{x}_0$ from a distribution of interest, the goal of a diffusion model is to learn to model its true data distribution $p_0(\bm{x}_0)$. 
	Once learned, we can generate new samples from our approximate model at will. The diffusion model builds a diffusion process by defining a forward SDE starting from $p_0(\bm{x}_0)$ as follows:
	\begin{equation}\label{eqn_forward}
		d\bm{x}_t=\bm{a}(\bm{x}_t,t)dt+b(t)d\bm{B}_t,
	\end{equation}
	where $t\in[0, T]$, $\bm{B}_t$ is the standard Wiener process (a.k.a., Brownian motion), $\bm{a}(\cdot,t):\mathbb{R}^d\to\mathbb{R}^d$ is a vector valued function called the drift coefficient, and $b(\cdot):\mathbb{R}\to\mathbb{R}$ is a scalar function known as the diffusion coefficient. 
	
	By starting from samples of $\bm{x}_T\sim p_T(\bm{x}_T)$ and reversing the process, we can obtain samples $\bm{x}_0\sim p_0(\bm{x}_0)$ . The reverse of a diffusion process is also a diffusion process, running backwards in time and given by the following reverse-time SDE:
	\begin{equation}\label{eqn_reverse}
		d\bm{x}_t=\left(\bm{a}(\bm{x}_t,t)-b(t)^2\nabla_{\bm{x}_{t}}\log p_t(\bm{x}_t)\right)dt+b(t)d\bm{\bar{B}}_t,
	\end{equation}
	where $\bm{\bar{B}}_t$ is a standard Wiener process when the time flows backwards from $T$ to 0. The only unknown term $\nabla_{\bm{x}_t}\log p_t(\bm{x}_t)$ is the score function of the marginal
	density $p_t(\bm{x}_t)$.

	To estimate $\nabla_{\bm{x}_t}\log p_t(\bm{x}_t)$, we can train a time-dependent score-based model $\bm{s}_{\bm{\theta}}(\bm{x}_t,t)$ with
	\begin{equation}
		\begin{aligned}			&\bm{\theta}^*=\mathop{\arg\min}_{\bm{\theta}}\  \mathbb{E}_{t\sim\mathcal{U}[0,T]}\left\{\lambda_t\mathbb{E}_{\bm{x}_0}\mathbb{E}_{\bm{x}_t|\bm{x}_0}\left[\left\|\bm{s}_{\bm{\theta}}(\bm{x}_t,t)-\nabla_{\bm{x}_t}\log p_{0t}(\bm{x}_t|\bm{x}_0)\right\|_2^2\right]\right\},
		\end{aligned}
	\end{equation}
	where $p_{0t}(\bm{x}_t|\bm{x}_0)$ is the transition kernel and can be obtained by the forward process (\ref{eqn_forward}). When $\bm{a}(\cdot,t)$ is affine, the transition kernel is always a Gaussian distribution, where the mean and variance are often known in closed forms (\cite{sarkka2019applied}). With sufficient data and model capacity, score matching ensures that the optimal solution $\bm{s}_{\bm{\theta}^*}(\bm{x}_t,t)$ approximates $\nabla_{\bm{x}_t}\log p_t(\bm{x}_t)$ for almost all $\bm{x}_t$ and $t$.

	\subsection{Gradient guidance}
	\label{sec:guidance}
	
	A direct application of diffusion models to CCP (\ref{eq1}) is infeasible, as this requires sampling from the product distribution $Q_{\beta}(\bm{x}|\rho)\propto r(\bm{x}|\rho)B_{\beta}(\bm{x})$, whereas only samples from $r(\bm{x}|\rho)$ are accessible. Therefore, obtaining a precise characterization of the score function of $Q_{\beta}(\bm{x}|\rho)$ and its diffused version is crucial.
	
	For a given data set $\mathcal{D}=\{(\bm{x}^{(i)},\rho^{(i)})\}_{i=1}^{N}$, we use its empirical $p_{0}(\bm{x}_0|\rho)$ to approximate the unknown distribution $r(\bm{x}_0|\rho)$ and denote $\tilde{p}_0(\bm{x}_0|\rho)\propto p_{0}(\bm{x}_0|\rho)B_{\beta}(\bm{x}_0)$. The diffused distribution is then given by the forward process (\ref{eqn_forward}), i.e., 
	\begin{equation}
		\begin{aligned}
			&p_t(\bm{x}_t|\rho) = \int_{\bm{x}_0} p_{0t}(\bm{x}_t|\bm{x}_0)p_{0}(\bm{x}_0|\rho)d\bm{x}_0, \\
			&\tilde{p}_t(\bm{x}_t|\rho) = \int_{\bm{x}_0} p_{0t}(\bm{x}_t|\bm{x}_0)\tilde{p}_0(\bm{x}_0|\rho)d\bm{x}_0 \propto \int_{\bm{x}_0} p_{0t}(\bm{x}_t|\bm{x}_0)p_{0}(\bm{x}_0|\rho)B_{\beta}(\bm{x}_0)d\bm{x}_0.
		\end{aligned}
	\end{equation}
	
	In order to sample with the reverse process (\ref{eqn_reverse}), we need to characterize the score function of the diffused product distribution $\nabla_{\bm{x}_t}\log \tilde{p}_t(\bm{x}_t|\rho)$, which is given by the following theorem.
	
	\textbf{Theorem 1.} For any given $\beta>0$, there exists $\hat{\bm{x}}_0(\bm{x}_t)$ such that the score function of the diffused product distribution can be formulated as
	\begin{equation}
		\nabla_{\bm{x}_t}\log \tilde{p}_t(\bm{x}_t|\rho) = \nabla_{\bm{x}_t}\log p_t(\bm{x}_t|\rho) \underbrace{- \beta  \nabla_{\bm{x}_t}f\big(\hat{\bm{x}}_0(\bm{x}_t)\big)}_{\text{gradient guidance $\bm{G}_t$}},
	\end{equation}
	where $\nabla_{\bm{x}_t}\log p_t(\bm{x}_t|\rho)$ is the score function of the diffused data distribution and $\hat{\bm{x}}_0(\bm{x}_t)$ satisfies
	\begin{equation}
		f(\hat{\bm{x}}_0(\bm{x}_t)) = -\frac{1}{\beta}\log\Big(\int_{\bm{x}_0} p_{t0}(\bm{x}_0|\bm{x}_t,\rho)B_{\beta}(\bm{x}_0)d\bm{x}_0\Big).
	\end{equation}
	
	Theorem 1 demonstrates that sampling from the product distribution can be accomplished by introducing a gradient guidance term during the sampling process of the original data distribution, which has a strong connection between the posteriori $p_{t0}(\bm{x}_0|\bm{x}_t,\rho)$ and the Boltzmann distribution $B_{\beta}(\bm{x}_0)$. 
	
	Next, we present a special case where the gradient guidance terms admit explicit expressions.
	
	\textbf{Corollary 1.} Assume that $p_{t0}(\bm{x}_0|\bm{x}_t,\rho) = \mathcal{N}(\bm{x}_0|\bm{\mu}_{0|t},\sigma_{0|t}^2\bm{I})$, then we have the following results.
	\begin{itemize}[leftmargin=10pt]
		\item \textbf{First-order guidance:} For $f\in\mathcal{C}^{1}(\mathbb{R}^n,\mathbb{R})$, we get
		\begin{equation}
			\label{eq:first}
			\bm{G}_t = - \beta \nabla_{\bm{x}_t} f(\bm{x}_t).
		\end{equation}
		\item \textbf{Second-order guidance:} For $f\in\mathcal{C}^{2}(\mathbb{R}^n,\mathbb{R})$, we get
		\begin{equation}
			\label{eq:second}
			\bm{G}_t = -\frac{1}{\sigma_{0|t}^2}\left[\bm{H}^{-1}\left((-\nabla_{\bm{x}_t}^2 f(\bm{x}_t) \bm{x}_t+\nabla_{\bm{x}_t} f(\bm{x}_t)) -\frac{1}{\beta \sigma_{0|t}^2}\bm{\mu}_{0|t}\right)+\bm{\mu}_{0|t}\right],
		\end{equation}
		where $\bm{H} = \nabla_{\bm{x}_t}^2 f(\bm{x}_t)+\frac{1}{\beta\sigma_{0|t}^2}\bm{I}$.
	\end{itemize}

	It is worthwhile noting that, for $p_{0}(\bm{x}_0|\rho)=\mathcal{N}(\bm{x}_0|\bm{\mu}_0,\sigma_0^2\bm{I})$ and the Gaussian transition kernel, the assumption in Corollary 1 holds and the parameters $(\bm{\mu}_{0|t}, \sigma_{0|t})$ can be expressed explicitly. In practice, we can use Tweedie's formula (\cite{efron2011tweedie}) to obtain an estimator of $\bm{\mu}_{0|t}$, and treat the variance as a hyper parameter; see Subsection \ref{subsec:samp} for details on this. 	Although the second-order guidance requires computing the inverse of a general Hessian matrix, which may be computationally expensive, it brings faster convergence and better variance reduction.

	\section{GGDOpt for CCP}
	\label{sec:GGDOpt}
	
	In this section, we give our GGDOpt framework for CCP (\ref{eq1}). The whole process can be divided into three stages: data generation, diffusion and learning, and sampling with guidance. More specifically, in the data generation stage, a collection of points satisfying the chance constraint is generated to characterize the nonconvex feasible set. The diffusion and learning stage progressively inject noise to convexify the nonconvex feasible region and learn the score function of the conditional distribution in order to perform sampling. After learning, the sampling with guidance stage iteratively runs the reverse process with an extra gradient guidance to sample from the product distribution, which will asymptotically converge to an optimal solution to problem (\ref{eq1}). Next, we present the details of the three stages in GGDOpt one by one.

	\subsection{Stage 1: data generation}
	\label{subsec:data}
	
	First we give an efficient approach to generate high-quality data that satisfy the chance constraint while maintaining lower objective values. Suppose that we have a set of samples $\{\bm{h}^{(\ell)}\}_{\ell=1}^L$, denote the empirical mean $\bar{\bm{h}} = \frac{1}{L}\sum_{\ell=1}^{L}\bm{h}^{(\ell)}$. Notice that in most of cases, it's much easier to solve the following deterministic restricted problem (RP) with a fixed $\bar{\bm{h}}$:
	\begin{equation}
		\label{eq:res}
		\begin{aligned}
			\min_{\bm{x}}\quad&f(\bm{x})\\
			\mathrm{s.t.}\quad& \bm{g}(\bm{x},\bar{\bm{h}})\geq\bm{z}_i,
		\end{aligned}
	\end{equation}
	where $\bm{z}_i\geq \bm{0}$ is a given restriction, $i=1,\ldots,N$. Let $\bm{x}(\bm{z}_i)$ denote the solution to problem (\ref{eq:res}) for a given $\bm{z}_i$. As the smallest element $z_{min}$ in $\bm{z}_i$ increases, the probability of the nonlinear constraint $g(\bm{x}(\bm{z}_i),\bm{h})\geq 0$ also increases.
	Then, solving problem (\ref{eq:res}) allows us to generate high-quality data that satisfies the chance constraint for arbitrary $\rho\in(0,1)$ while enjoys low objective values.
	
	Since the distribution of the random variable $\bm{h}$ is unknown, referring SAA method, we approximate the chance constraint using the empirical distribution over samples $\{\bm{h}^{(\ell)}\}_{\ell=1}^L$. Then, after getting $\bm{x}(\bm{z}_i)$, we have
	\begin{align}
		\text{Prob}_{\bm{h}}\{\bm{g}(\bm{x}(\bm{z}_i),\bm{h})\geq\bm{0}\}\approx
		\underbrace{\frac{1}{L}\sum_{l=1}^{L}\ell_{0/1}(\bm{g}(\bm{x}(\bm{z}_i),\bm{h}^{(\ell)}))}_{1-\rho^{(i)}},
	\end{align} 
	where $\ell_{0/1}(\bm{g})=1$ if $\bm{g}\geq\bm{0}$ and $\ell_{0/1}(\bm{g})=0$ otherwise. By calculating the empirical $\rho^{(i)}$, an asymptotic approximation of the underlying probability is obtained, requiring no assumption on the underlying distribution $P$. In the appendix, we give a tight lower bound for the probability constraint $\text{Prob}_{\bm{h}}\{\bm{g}(\bm{x}(\bm{z}_i),\bm{h})\geq\bm{0}\}$ if the variance and the mean of the random variable $\bm{h}$ are known, which is helpful to obtain a better approximation $\rho^{(i)}$.
	
	Let $\bm{x}^{(i)}:=\bm{x}(\bm{z}_i)$ and repeating the above process, i.e., solving problem (\ref{eq:res}) and estimating $\rho^{(i)}$, and gradually increasing $\bm{z}_i$, we can generate a collection of data points $\mathcal{D} = \{\bm{x}^{(i)},\rho^{(i)}\}_{i=1}^N$, which are then used to train our GGDOpt in the next stages.

	\subsection{Stage 2: diffusion and learning}
	
	From Theorem 1, we observe that the score function of the diffused product distribution has two terms, the conditional score $\nabla_{\bm{x}_t}\log p_t(\bm{x}_t|\rho)$ and the gradient guidance term $\bm{G}_t$ for which explicit forms of first- and second-order guidances have been derived in Corollary 1. Then the challenge reduces to learning the conditional score $\nabla_{\bm{x}_t}\log p_t(\bm{x}_t|\rho)$.
	
	In practice, naively conditioning a standard diffusion model by appending the conditioned variable at each step of the sampling process does not work well, as the model often ignores the conditioned information. Related works on conditional score estimation have been studied in (\cite{dhariwal2021diffusion}, \cite{dhariwal2021diffusion}, \cite{ho2022classifier}). Here we propose to use the classifier-free guidance (\cite{ho2022classifier}) to give an approximation of $\nabla_{\bm{x}}\log p_t(\bm{x}|\rho)$. 
	
	Instead of training a separate classifier model, classifier-free guidance choose to train an unconditional score estimator to approximate $\nabla_{\bm{x}_t}\log p_t(\bm{x}_t)$ together with the conditional score estimator to approximate $\nabla_{\bm{x}_t}\log p_t(\bm{x}_t|\rho)$. Specificity,
	we train a single model $\bm{s}_{\bm{\theta}}(\bm{x}_t,t,\rho)$, and the conditioning information $\rho$ is randomly discarded as empty set $\emptyset$ with probability $p_{uncond}$ to train unconditionally. Then the conditional score $\nabla_{\bm{x}_t}\log p_t(\bm{x}_t|\rho)$ is estimated by
	\begin{equation}
		\nabla_{\bm{x}_t}\log p_t(\bm{x}_t|\rho) \approx (1+w) \bm{s}_{\bm{\theta}}(\bm{x}_t,t,\rho) -w\bm{s}_{\bm{\theta}}(\bm{x}_t,t,\emptyset),
	\end{equation}
	for a given weight parameter $w$. Specifically, for the given data set $\mathcal{D}$ and network $\bm{s}_{\bm{\theta}}(\bm{x}_t,t,\rho)$ parameterized by $\bm{\theta}$, the training objective is defined as
	\begin{equation}
		\label{eq:loss}
		\text{Loss}(\bm{\theta}) = \mathbb{E}_{t\sim\mathcal{U}[0,T]}\left\{\mathbb{E}_{\bm{x}_0,\rho}\mathbb{E}_{\bm{x}_t|\bm{x}_0}\left[\left\|\bm{s}_{\bm{\theta}}(\bm{x}_t,t,\rho)-\nabla_{\bm{x}_t}\log p_{0t}(\bm{x}_t|\bm{x}_0)\right\|_2^2\right]\right\},
	\end{equation}
	and trained with Adam (\cite{kingma2014adam}). The training process of GGDOpt is given in Algorithm \ref{alg1}.

	\subsection{Stage 3: sampling with guidance}
	\label{subsec:samp}
	
	Given the forward process (\ref{eqn_forward}), the corresponding reverse process is given by the following reverse-time SDE with trained $\bm{s}_{\bm{\theta}}(\bm{x}_t,t,\rho)$ and gradient guidance $\bm{G}_t$:
	\begin{equation}
		\label{eq:rev}
		d\bm{x}_t=\left[\bm{a}(\bm{x}_t,t)-b(t)^2\big(\tilde{\bm{s}}_{\bm{\theta}}(\bm{x}_t,t,\rho)+\bm{G}_t \big)\right]dt+b(t)d\bm{\bar{B}}_t,
	\end{equation}
	where
	\begin{equation}
		\label{eq:condition}
		\tilde{\bm{s}}_{\bm{\theta}}(\bm{x}_t,t,\rho)= (1+w) \bm{s}_{\bm{\theta}}(\bm{x}_t,t,\rho) -w\bm{s}_{\bm{\theta}}(\bm{x}_t,t,\emptyset).
	\end{equation}

	For the first-order gradient guidance $\bm{G}_t$ in (\ref{eq:first}), we directly use the gradient of the objective scaled by a hyper parameter $\beta$. For the second-order gradient guidance (\ref{eq:second}), we need to give the posterior mean and variance $(\bm{\mu}_{0|t},\sigma_{0|t}^2)$. Here we use Tweedie's formula (\cite{efron2011tweedie}) to get an estimator of the posterior mean as follows: 
	\begin{equation}
		\bm{\mu}_{0|t} = \mathbb{E}\left[\bm{x}_0 | \bm{x}_t,\rho\right] = \frac{1}{\sqrt{\bar{\alpha}_t}} (\bm{x}_t + (1-\bar{\alpha}_t)\tilde{\bm{s}}_{\bm{\theta}}(\bm{x}_t,t,\rho)),
	\end{equation}
	with priori $p_{0t}(\bm{x}_t|\bm{x}_0) = \mathcal{N}(\bm{x}_t |\sqrt{\bar{\alpha}_t} \bm{x}_0, (1-\bar{\alpha}_t)\bm{I})$ for a specific noising schedule $\bar{\alpha}_t$. 
	
	While Tweedie’s formula theoretically provides both the posterior mean and covariance, $\boldsymbol{\Sigma}_{0|t}=(1-\bar{\alpha}_t)(\boldsymbol{I}+(1-\bar{\alpha}_t)\nabla^2\log p(\boldsymbol{x}_t))$, computing the covariance requires evaluating the Hessian of $\log p(\boldsymbol{x})$. In our framework, the score function $\boldsymbol{s}_{\boldsymbol{\theta}}$ is parameterized by a neural network, and computing its second derivatives involves backpropagation through the network’s Jacobian, which is computationally expensive, especially in high dimensions. To strike a balance between performance and efficiency, we choose to treat the covariance as a tunable hyper parameter $\sigma^2$. In the appendix, we give a detailed comparison between the fully Tweedie-based method and our approach to show that using a fixed variance can be a practical and robust alternative.
	
	Then the second-order guidance can be calculated by
	\begin{equation}
		\bm{G}_t = -\frac{1}{\sigma^2}\left[(\nabla^2 f(\bm{x}_t)+\frac{1}{\beta\sigma^2}\bm{I})^{-1}\left((-\nabla^2 f(\bm{x}_t) \bm{x}_t+\nabla f(\bm{x}_t)) -\frac{1}{\beta \sigma^2}\bm{\mu}_{0|t}\right)+\bm{\mu}_{0|t}\right],
	\end{equation}
	and the sampling process of GGDOpt is given in Algorithm \ref{alg2}.

	\begin{figure}[t]
		\centering
		\begin{minipage}[t]{0.49\textwidth}
			\centering
			\begin{algorithm}[H]
				\caption{Training of GGDOpt}
				\label{alg1}
				\begin{algorithmic}[1]
					\REQUIRE
					$\{(\bm{x}^{(i)},\rho^{(i)})\}_{i=1}^{N}\sim p_0(\bm{x}|\rho)$.
					\ENSURE
					$\bm{s}_{\bm{\theta}^*}(\bm{x},t,\rho)$.
					\STATE \textbf{repeat}
					\STATE \ \ Load $(\bm{x}_0,\rho_0)\sim p_0(\bm{x}|\rho)$.
					\STATE \ \ Set $\rho\leftarrow \emptyset$ with probability $p_{uncond}$.
					\STATE \ \ Sample $t\sim\mathcal{U}[0,T]$.
					\STATE \ \ Generate $\bm{x}_t\sim p_{0t}(\bm{x}_t|\bm{x}_0)$.
					\STATE \ \ Take gradient descent step on (\ref{eq:loss}).
					\STATE \textbf{until} converged.
				\end{algorithmic}
			\end{algorithm}
		\end{minipage}
		\begin{minipage}[t]{0.49\textwidth}
			\begin{algorithm}[H]
				\caption{Sampling of GGDOpt}
				\label{alg2}
				\begin{algorithmic}[1]
					\REQUIRE
					$\bm{s}_{\bm{\theta}^*}(\bm{x},t,\rho)$, objective $f$.
					\ENSURE
					$\bm{x}_0^*$.
					\STATE $\bm{x}_T\sim p_T$.
					\FOR{$t=T,...,1$}
					\STATE Calculate $\tilde{\bm{s}}_{\bm{\theta}}(\bm{x}_t,t,\rho)$ with (\ref{eq:condition}).
					\STATE Calculate $\bm{G}_t$ with (\ref{eq:first}) or (\ref{eq:second}) .
					\STATE Take guided sampling step with (\ref{eq:rev}).
					\ENDFOR
					\RETURN $\bm{x}_0^* =\bm{x}_{0}$.
					\vspace{1pt}
				\end{algorithmic}
			\end{algorithm}
		\end{minipage}
	\end{figure}

	\section{Convergence analysis}
	\label{sec:conv}
	
	In this section, we give the convergence analysis of the proposed GGDOpt framework in both theoretical and practical aspects. We show that: theoretically, the samples generated by the sampling process will concentrate around the points with the lowest function values within the support of the data distribution; and practically, the gap between the expected function values of generated samples and the optimal value will be bounded by two components.
	
	\subsection{Theoretical convergence}
	
	As provided by (\cite{pidstrigach2022score}), under mild assumptions, the sampling distribution of the standard diffusion model will have the exact same support as the data distribution. But what if we introduce an extra gradient guidance term? For a given $\rho$, denote $\mathcal{D}_{\rho}=\{\bm{x}^{(i)}\mid (\bm{x}^{(i)},\rho^{(i)})\in\mathcal{D}, \rho^{(i)}\leq \rho\}$ as the approximated feasible set of $\mathcal{X}_{\rho}$. The following theorem says that in our settings, as $T\rightarrow \infty$ and $\beta\rightarrow\infty$, the samples of GGDOpt will concentrate around the points with the lowest function values within the support of the data distribution $\mathcal{D}_{\rho}$ for any given $\rho$.

	\textbf{Theorem 2.} For any given $\rho\in(0,1)$, suppose that there exists a constant $\delta$ such that the error in the score estimation can be bounded as: 
	\begin{equation}
		\label{eq:scoreerror}
		\|\tilde{\bm{s}}_{\bm{\theta}}(\bm{x}_t,t,\rho) + \bm{G}_t - \nabla_{\bm{x}_t}\log \tilde{p}_t(\bm{x}_t|\rho)\| \leq \delta, \quad \forall~\bm{x}_t.
	\end{equation}
	For samples $\tilde{\bm{x}}_{sample}\sim p_{sample}(\bm{x}_0|\rho)$ generated by the reverse process
	\begin{equation}
		d\bm{x}_t=\left[\bm{a}(\bm{x}_t,t)-b(t)^2\big(\tilde{\bm{s}}_{\bm{\theta}}(\bm{x}_t,t,\rho)+\bm{G}_t \big)\right]dt+b(t)d\bm{\bar{B}}_t,
	\end{equation}
	with prior $p_{prior} = \mathcal{N}(\bm{0},\bm{I})$, affine drift coefficients $\bm{a}(\cdot,t)$, and
	\begin{equation}
		\tilde{\bm{s}}_{\bm{\theta}}(\bm{x}_t,t,\rho)= (1+w) \bm{s}_{\bm{\theta}}(\bm{x}_t,t,\rho) -w\bm{s}_{\bm{\theta}}(\bm{x}_t,t,\emptyset),
	\end{equation}
	as $T\rightarrow \infty$, $p_{sample}(\bm{x}_0|\rho)$ will have the same support as $\tilde{p}_{0}(\bm{x}_0|\rho)$. Further, as $\beta\rightarrow\infty$, $\tilde{\bm{x}}_{sample}$ will concentrate around $\bm{x}^* = \mathop{\arg\min}_{\bm{x}\in\mathcal{D}_{\rho}} f(\bm{x})$.
	
	The assumption in the score estimation error (\ref{eq:scoreerror}) quantifies the approximation accuracy of the trained score network relative to the true score function. It depends on the training quality of the neural network and the expressiveness of the model class. This type of assumption is common in the theoretical analysis of diffusion models (see, e.g., \cite{pidstrigach2022score}, \cite{de2021diffusion}) and is used to establish convergence results in generative modeling and sampling.
	
	\subsection{Practical error bound}
	
	In practice, the forward process cannot reach the stationary distribution and the training is not perfect. This results in the failure of the sample distribution to strictly concentrate on the data points. This will lead to two components of errors: the limited time length error $I_1$ and limited inverse temperature error $I_2$, which are given as follows:
	\begin{equation}
		\begin{aligned}
			|\mathbb{E}[f(\tilde{\bm{x}}_{sample})]-f(\bm{x}^{*})|&\leq|\underbrace{\mathbb{E}[f(\tilde{\bm{x}}_{sample})]-\mathbb{E}[f(\bm{x}^\pi)]|}_{I_1} + \underbrace{|\mathbb{E}[f(\bm{x}^\pi)]-f(\bm{x}^*)|}_{I_2}.
		\end{aligned}
	\end{equation}
	In the above, $\tilde{\bm{x}}_{sample}$ is sampled from the reverse process (\ref{eq:rev}), $\bm{x}^\pi$ follows the strong solution $p^{\pi}$ to the Fokker-Planck equation of (\ref{eq:rev}), and $\bm{x}^* = \mathop{\arg\min}_{\bm{x}\in\mathcal{D}_{\rho}} f(\bm{x})$. Next, we will give practical error bounds of both the two components with finite $T$ and $\beta$.
	
	\textbf{Assumption 1.} We assume the following conditions hold:
	\begin{itemize}[leftmargin=10pt]
		\item The forward process is given by $d\bm{x}= b(t)d\bm{B}_t$;
		\item The reverse process starts in $p_{prior} = \mathcal{N}(\bm{m}_T, \bm{\Sigma}_T)$ where $\bm{m}_T = \mathbb{E}[\tilde{p}_0(\bm{x}_0|\rho)]$ and $\bm{\Sigma}_T = \text{Cov}(\tilde{p}_0(\bm{x}_0|\rho))+T\cdot \bm{I}$;
		\item The objective function $f(\bm{x})$ satisfies $\|\nabla_{\bm{x}}f(\bm{x}) \|_2\leq C_1\|\bm{x} \|_2+C_2$.
	\end{itemize}
	
	The first two conditions in Assumption 1 correspond to the VE SDE in (\cite{songscore}) and are primarily used to characterize the discrepancy between the end distribution and the prior distribution. The third assumption is common in the convergence analysis of stochastic optimization and sampling algorithms (see, e.g., \cite{raginsky2017non}). In practice, Assumption 1 holds for a broad class of functions, including smooth bounded functions and quadratic objectives, which frequently arise in real-world optimization problems.
	
	\textbf{Theorem 3.} Under Assumption 1, denote $\sigma^{(k)}, k=1,\ldots,n$, the eigenvalues of $\bm{\Sigma}_T$. For any given $\rho\in(0,1)$, denote $N_{\rho}=|\mathcal{D}_{\rho}|$ and $\bm{x}^* = \mathop{\arg\min}_{\bm{x}\in\mathcal{D}_{\rho}} f(\bm{x})$. Then for any given $T>0$ and $\beta>0$, the optimization error can be bounded by 
	\begin{equation}
		\begin{aligned}
			|\mathbb{E}[f(\tilde{\bm{x}}_{sample})]-f(\bm{x}^{*})|\leq \underbrace{C_I\big(\sqrt{C_T}+\big(C_T/2\big)^{1/4}\big)}_{I_1} + \underbrace{(N_\rho-1)\max_{\bm{x}\in\mathcal{D}_{\rho}}|f(\bm{x})-f(\bm{x}^*)|e^{-\beta \delta_{\rho}}}_{I_2},
		\end{aligned}
	\end{equation}
	where $C_T= \frac{1}{2}\log\left(\prod_{k=1}^{n}(\sigma^{(k)}/T)\right)$ and $C_I, \delta_{\rho}$ are constants.
	
	Theorem 3 provides a non-asymptotic convergence result of GGDOpt with limited time length and inverse temperature. As $T\rightarrow\infty$ and $\beta\rightarrow\infty$, the optimization error goes to zero and GGDOpt is shown to generate asymptotically optimal solutions.
	
	\section{Experimental results}
	\label{sec:expe}
	
	In this section, we perform numerical experiments on both synthetic datasets and a wireless communications waveform design problem. To generate the data, we employ CVX (\cite{grant2008cvx}) to solve the restricted problem (\ref{eq:res}). In the diffusion and learning stage, we set $T = 1000$ with a linear noise schedule $\eta(t)$ ranging from 0.0001 to 0.02, and let $\bm{a}(\bm{x},t) = -\frac{1}{2}\eta(t)\bm{x}$ and $b(t) = \sqrt{\eta(t)}$. In the sampling with guidance stage, we evaluate both first- and second-order gradient guidances via implementing a DDIM-based technique (\cite{song2020denoising}) with a descaled time step $T'=100$ for accelerated sampling. We employ two variants of the U-Net model (\cite{ronneberger2015u}) as our score estimator: U-Net-1D for the linear chance constrained problem and both for robust waveform design. Additional experimental details are provided in the supplementary materials.
	
	\subsection{Linear chance constrained problem}
	
	Consider the following linear chance constrained problem:
	\begin{equation}
		\label{eq:linear}
		\begin{aligned}
			\min_{\bm{x}\in\mathbb{R}^n}\quad&\frac{1}{2}\bm{x}^{\top}\bm{x}+\bm{b}^{\top}\bm{x}\\
			\mathrm{s.t.}\quad&\text{Prob}_{\bm{c}\sim p_{\bm{c}}}\{\bm{c}^{\top}\bm{x}+d\geq0\}\geq 1-\rho,
		\end{aligned}
	\end{equation}
	where $p_{\bm{c}} = \mathcal{N}(\bm{c};\bar{\bm{c}},\bm{I})$ and $(\bm{b},\bar{\bm{c}},d,\rho)$ are hyper parameters selected from a test set. The above problem can be reformulated as a second-order conic (SOC) program, for which CVX (\cite{grant2008cvx}) is used for solution. To generate training data, we solve the restricted version of problem (\ref{eq:linear}) for $N=1000$ values of $z$ linearly spaced in the interval $\left[0, 0.5\right]$. Then we execute the reverse process with first- and second-order gradient guidance to generate samples.

	We compare our proposed GGDOpt against different types of SAA methods for solving the problem, using the corresponding CVX solutions as performance benchmarks. Each algorithm was executed 100 times (except CVX). The results with $n=8$ are presented in Table \ref{tab:linear}. 
	
	\begin{table}[htbp]
		\caption{Comparison results on the linear chance constrained problem (\ref{eq:linear})}
		\label{tab:linear}
		\renewcommand{\arraystretch}{1.5}
		\resizebox{\textwidth}{!}{
			\begin{tabular}{lccccc}
				\hline
				Method                  & Repeat & FvalMean & FvalStd & FvalMedian  & Runtime \\ \hline
				SOC\_CVX (\cite{grant2008cvx})     & 1      & \textbf{-0.6586}  & 0    & -0.6586     & 0.3214  
				\\ \hline
				SAA\_CVaR (\cite{nemirovski2007convex})      & 100    & -0.5893  & 0.0248  & -0.5869  & 0.3063 \\
				SAA\_MIP (\cite{pagnoncelli2009sample})      & 100    & -0.6281  & 0.0157  & -0.6318  & 15.4502 \\
				SAA\_PDCA (\cite{wang2023proximal})          & 100    & -0.6389  & 0.0314  & -0.6408  & 0.6276 \\
				SAA\_SNSCO (\cite{zhou20240})       & 100    & 0.8051   & 3.4014  & -0.6371   & 0.2793 				
				\\ \hline
				GGDOpt\_WithoutGuidance & 100    & 0.3481   & 0.5486  & 0.2798        & 0.0465  \\
				GGDOpt\_First-order      & 100    & -0.6483  & \textbf{0.0051}  & -0.6488       & \textbf{0.0486}  \\
				GGDOpt\_Second-order     & 100    & \textbf{-0.6491}  & \textbf{0.0056}  & \textbf{-0.6503}       & \textbf{0.0507}  \\ \hline
			\end{tabular}
		}
	\end{table}
	
	The results in Table \ref{tab:linear} demonstrate that, compared to the SOC\_CVX method, which requires explicit knowledge of the underlying distribution, GGDOpt can approximately find the global minimizer with only samples from distribution $p_{\bm{c}}$ while simultaneously achieving significant overhead reduction. Compared to SAA methods, GGDOpt achieves superior performance in terms of lower function values and enhanced numerical stability under the effect of gradient guidance. 
	
	As expected, the runtime increases with the problem dimension. However, both the first- and second-order versions of GGDOpt remain consistently faster than the baseline SAA\_PDCA method across all dimensions. Moreover, the increase in runtime is moderate, indicating that our approach scales favorably even in high-dimensional settings. 
	
	Furthermore, as the runtime increases with the problem dimension, both the first- and second-order versions of GGDOpt reduce the computational time by approximately 80\% compared with , offering substantial efficiency improvements. More detailed experimental results on larger problem scale and computational costs are listed in the appendix.

	\subsection{Robust waveform design}
	
	Consider the following robust waveform design problem (\cite{wang2014outage})
	\begin{equation}
		\label{eq:robust}
		\begin{aligned}
			\operatorname*{min}_{\boldsymbol{S}_{1},\ldots,\boldsymbol{S}_{K}\in\mathbb{R}^{N_{t}\times N_{t}}}\ &\sum_{i=1}^K\mathrm{Tr}(\boldsymbol{S}_i)\\
			\mathrm{s.t.}\quad \quad \ \ \ &\mathrm{Prob}_{\boldsymbol{h}_{i}\sim\mathcal{N}(\bar{\boldsymbol{h}}_{i},\boldsymbol{C}_{i})}\{\mathrm{R}_{i}\geq r_{i}\}\geq 1-\rho_{i},i=1,2,\ldots,K,\\
			&\boldsymbol{S}_{1},\ldots,\boldsymbol{S}_{K}\succeq\boldsymbol{0},i=1,2,\ldots,K,
		\end{aligned}
	\end{equation}
	where $N_t$ is the number of antennas at the base station and $K$ is the total number of users. For each user $i$, $\bm{S}_i \succeq \bm{0}, \bm{h}_i, R_i$ and $r_i\geq0$ denote the signal covariance matrix (to be designed), the random channel vector, the achievable rate, and the desired rate target, respectively.	
	
	Firstly, we use U-Net-2D as the score estimator. Notice that during the data generation, all the solutions to the restricted problem (\ref{eq:res}) exhibit a rank-one structure (\cite{huang2007complex}, \cite{chang2008approximation}, \cite{huang2020quadratic}). Remarkably, the generated samples maintain this rank-one property (with dominant eigenvalue accounting for >99\% of the total eigenvalue) after training, suggesting that the solutions to the robust waveform design problem (\ref{eq:robust}) inherently reside on a rank-one manifold with extremely high probability (\cite{wang2014outage}), which GGDOpt successfully captures. This implies that rank-one decomposition can be effectively applied after generation, enabling the use of U-Net-1D as a score estimator to reduce computational costs in both training and sampling process.
	
	Table \ref{tab:waveform} summarizes the comparison results of GGDOpt and two state-of-the-art methods for solving problem (\ref{eq:robust}) with $N_t = 16$ and $K = 3$, where the worst probabilities that the chance constraints satisfy for $K$ users are underlined. Notably, both baseline methods rely on explicit knowledge of the underlying distribution, whereas GGDOpt operates solely based on samples.  The results show that GGDOpt outperforms existing convex approximation methods, achieving superior feasible solutions outside the convex restriction of the feasible set, while significantly reducing computational overhead. Complete experimental details are provided in the appendix.

	\begin{table}[htbp]
		\caption{Optimization methods comparison for robust waveform design}
		\label{tab:waveform}
		\renewcommand{\arraystretch}{1.5}
		\resizebox{\textwidth}{!}{
			\begin{tabular}{clllll}
				\hline
				Method                                                                                          & Metric      & $\rho=0.05$         & $\rho=0.10 $        &$\rho=0.15$         & $\rho=0.20 $      \\ \hline
				\multirow{3}{*}{\begin{tabular}[c]{@{}c@{}}Sphere   Bounding\\ \cite{ben2000robust} \end{tabular}}           & Probability & \underline{0.99}; 0.99; 0.99 & \underline{0.99}; 0.99; 0.99 & \underline{0.99}; 0.99; 0.99 & \underline{0.99}; 0.99; 0.99 \\
				& FuncValue   & 0.1374           & 0.1366           & 0.1361           & 0.1357           \\
				& Runtime     & 1.4688           & 1.4375           & 1.4113           & 1.3875           \\ \hline
				\multirow{3}{*}{\begin{tabular}[c]{@{}c@{}}Bernstein-type   Inequality\\ \cite{wang2014outage} \end{tabular}} & Probability & 0.96; \underline{0.95}; 0.96 & 0.93; \underline{0.93}; 0.93 & 0.91; \underline{0.91}; 0.92 & 0.90; \underline{0.90}; 0.91 \\
				& FuncValue   & 0.1260           & 0.1253           & 0.1248           & 0.1244           \\
				& Runtime     & 1.2938           & 1.2813           & 1.2593           & 1.2652           \\ \hline
				\multirow{3}{*}{\begin{tabular}[c]{@{}c@{}}GGDOpt\\ First-order guidance \end{tabular}}                                                             & Probability & 0.99; \underline{0.95}; 0.99 & 0.92; 0.98; \underline{0.91} & 0.93; \underline{0.86}; 0.94 & 0.87; \underline{0.81}; 0.91 \\
				& FuncValue   & 0.1279           & 0.1265           & 0.1254           & 0.1247           \\
				& Runtime     & 0.0691           & 0.0628           & 0.0603           & 0.0635           \\ \hline
				\multirow{3}{*}{\begin{tabular}[c]{@{}c@{}}GGDOpt\\ Second-order guidance \end{tabular}}                                                            & Probability & 0.97; \underline{0.95}; 0.96 & \underline{0.90}; 0.94; 0.90 & 0.88; \underline{0.85}; 0.86 & 0.88; \underline{0.80}; 0.87 \\
				& FuncValue   & 0.1260           & \textbf{0.1246}           & \textbf{0.1239}           & \textbf{0.1237}           \\
				& Runtime     & \textbf{0.0788}           & \textbf{0.0712}           & \textbf{0.0687}           & \textbf{0.0682}           \\ \hline
			\end{tabular}
		}
	\end{table}
	
	\section{Conclusion}
	\label{sec:conc}
	
	In this paper, we have proposed GGDOpt, a gradient-guided diffusion framework that efficiently solves nonconvex CCP without requiring the exact distribution knowledge. By reformulating CCP as a sampling problem over the product of an unknown data distribution and a Boltzmann distribution, GGDOpt leverages both first- and second-order gradient information during reverse sampling. Theoretical convergence guarantees and practical error bounds are provided under mild assumptions. Experimental results demonstrate that GGDOpt outperforms existing methods in both solution quality and numerical stability with significant overhead reduction.
	
	\newpage
	\section*{Acknowledgments}
	The work of Boyang Zhang and Ya-Feng Liu was supported in part by the National Natural Science Foundation of China (NSFC) under Grant 12021001 and Grant 12371314. The work of Zhiguo Wang was supported in part by The National Key Research and Development Program of China under Grant 2020YFA0714003 and in part by NSFC under Grant 62203313.

	
	\bibliographystyle{plainnat}
	\bibliography{references}

	
	\newpage
	\section*{NeurIPS Paper Checklist}
	
	\begin{enumerate}
		
		\item {\bf Claims}
		\item[] Question: Do the main claims made in the abstract and introduction accurately reflect the paper's contributions and scope?
		\item[] Answer: \answerYes{} 
		\item[] Justification: The main results and contributions of this paper are all included in the abstract and introduction clearly.
		
		\item {\bf Limitations}
		\item[] Question: Does the paper discuss the limitations of the work performed by the authors?
		\item[] Answer: \answerYes{} 
		\item[] Justification: We point out all assumptions and discuss the limitations of the work thoroughly in the supplementary material.
		
		\item {\bf Theory assumptions and proofs}
		\item[] Question: For each theoretical result, does the paper provide the full set of assumptions and a complete (and correct) proof?
		\item[] Answer: \answerYes{} 
		\item[] Justification: All the assumptions used are included in the main paper, and the proofs are provided in the supplementary material.
		
		\item {\bf Experimental result reproducibility}
		\item[] Question: Does the paper fully disclose all the information needed to reproduce the main experimental results of the paper to the extent that it affects the main claims and/or conclusions of the paper (regardless of whether the code and data are provided or not)?
		\item[] Answer: \answerYes{} 
		\item[] Justification: The main configuration of experiments is claimed in the Experimental results section, and more details are provided in the supplementary material. We will release the code once the paper is published.

		\item {\bf Open access to data and code}
		\item[] Question: Does the paper provide open access to the data and code, with sufficient instructions to faithfully reproduce the main experimental results, as described in supplemental material?
		\item[] Answer: \answerNo{} 
		\item[] Justification: The data generation algorithm is provided in this paper and can be reproduced easily. The code is a straightforward implementation of the proposed framework, and will be released once the paper is published.

		\item {\bf Experimental setting/details}
		\item[] Question: Does the paper specify all the training and test details (e.g., data splits, hyperparameters, how they were chosen, type of optimizer, etc.) necessary to understand the results?
		\item[] Answer: \answerYes{} 
		\item[] Justification: The experimental settings are presented in the main paper, and full details are provided in the supplementary material.
		
		\item {\bf Experiment statistical significance}
		\item[] Question: Does the paper report error bars suitably and correctly defined or other appropriate information about the statistical significance of the experiments?
		\item[] Answer: \answerYes{} 
		\item[] Justification: In the experiments, we run multiple times for each method and the stability is shown in the main paper.
		
		\item {\bf Experiments compute resources}
		\item[] Question: For each experiment, does the paper provide sufficient information on the computer resources (type of compute workers, memory, time of execution) needed to reproduce the experiments?
		\item[] Answer: \answerYes{} 
		\item[] Justification: The information of the compute resources is provided in the supplementary material.
		
		\item {\bf Code of ethics}
		\item[] Question: Does the research conducted in the paper conform, in every respect, with the NeurIPS Code of Ethics \url{https://neurips.cc/public/EthicsGuidelines}?
		\item[] Answer: \answerYes{} 
		\item[] Justification: The research conducted in the paper conform, in every respect, with the NeurIPS Code of Ethics.

		\item {\bf Broader impacts}
		\item[] Question: Does the paper discuss both potential positive societal impacts and negative societal impacts of the work performed?
		\item[] Answer: \answerNA{} 
		\item[] Justification: This paper focus on the theoretical results of Gradient Guidance and a framework for solving chance constrained problems. There is no direct path to any negative applications of this paper.
		
		\item {\bf Safeguards}
		\item[] Question: Does the paper describe safeguards that have been put in place for responsible release of data or models that have a high risk for misuse (e.g., pretrained language models, image generators, or scraped datasets)?
		\item[] Answer: \answerNA{} 
		\item[] Justification: This paper poses no such risks.
		
		\item {\bf Licenses for existing assets}
		\item[] Question: Are the creators or original owners of assets (e.g., code, data, models), used in the paper, properly credited and are the license and terms of use explicitly mentioned and properly respected?
		\item[] Answer: \answerYes{} 
		\item[] Justification: All the original papers of used models and algorithms are properly cited in this paper.
		
		\item {\bf New assets}
		\item[] Question: Are new assets introduced in the paper well documented and is the documentation provided alongside the assets?
		\item[] Answer: \answerNA{} 
		\item[] Justification: This paper does not release new assets, and our code will be released once the paper is published.
		
		\item {\bf Crowdsourcing and research with human subjects}
		\item[] Question: For crowdsourcing experiments and research with human subjects, does the paper include the full text of instructions given to participants and screenshots, if applicable, as well as details about compensation (if any)? 
		\item[] Answer: \answerNA{} 
		\item[] Justification: This paper does not involve crowdsourcing nor research with human subjects.
		
		\item {\bf Institutional review board (IRB) approvals or equivalent for research with human subjects}
		\item[] Question: Does the paper describe potential risks incurred by study participants, whether such risks were disclosed to the subjects, and whether Institutional Review Board (IRB) approvals (or an equivalent approval/review based on the requirements of your country or institution) were obtained?
		\item[] Answer: \answerNA{} 
		\item[] Justification: This paper does not involve crowdsourcing nor research with human subjects.
		
		\item {\bf Declaration of LLM usage}
		\item[] Question: Does the paper describe the usage of LLMs if it is an important, original, or non-standard component of the core methods in this research? Note that if the LLM is used only for writing, editing, or formatting purposes and does not impact the core methodology, scientific rigorousness, or originality of the research, declaration is not required.
		\item[] Answer: \answerNA{} 
		\item[] Justification: The core method development in this research does not involve LLMs as any important, original, or non-standard components.
		
	\end{enumerate}

	\newpage
	\appendix
	
	\section*{Technical Appendices and Supplementary Material}

	\appendix
	
	
	\section{Limitations and future works}

	First, while empirical results demonstrate faster convergence with second-order guidance, theoretical guarantees of this acceleration remain to be established. Second, while the U-Net architecture serves as our baseline score estimator, it may not be optimal for all problem domains. Specialized network architectures that better capture the geometric structure of constraints may be investigated. Third, experimental results primarily focus on two specific types of problems, then further evaluation may be required to assess the effectiveness on a broader range of function types. 
	
	\section{Related works}
	
	\subsection{Chance constrained programming}
	
	CCP is a powerful modeling paradigm for optimization problems with uncertain constraints, with applications across engineering, finance, and beyond. Two common solution approaches are Convex Approximation (CA) and Sample Average Approximation (SAA). However, CA requires explicit distributional information, and SAA can be computationally expensive. Thus, designing an efficient framework for CCP under \textbf{unknown distributions} remains a pressing challenge.
	
	\subsection{Optimization via sampling}
	Traditional gradient‑based methods often converge to local minima under nonconvex settings. Sampling‑based algorithms, particularly Langevin Dynamics, have demonstrated strong performance in global optimization (\cite{ma2019sampling}). Compared to conventional optimizers, Sampling‑based algorithms can take fully advantages of data priors and solve nonconvex problems more effectively.
	
	\subsection{Learning to optimize}
	
	In order to improve the efficiency of optimization algorithms, learning-based methods are studied by \cite{chen2022learning}. Learning-based methods aim to learn a parameterized or semi-parameterized update rule of optimization without taking the form of any analytic update. Traditional learning-based methods simply learns the mapping between the input and output of the optimization algorithms, which may cause to fall into local minima. Consequently, generative sampling‑based models have attracted growing interest for optimization tasks.
	
	\subsection{Diffusion models for optimization}
	
	The rising prominence of diffusion models has spurred significant research interest in their underlying mathematical foundations and theoretical properties, as well as strategies to optimize their performance. At the same time, there are more and more researches on the application of diffusion model. How to use diffusion model to solve optimization problems is gradually attracting people's attention. In \cite{chung2022improving}, an additional correction term inspired by the manifold constraint is added into the reverse diffusion step to preserve the manifold constraint and data consistency, and used to solve the inverse problem. In \cite{krishnamoorthy2023diffusion}, a conditional diffusion model is trained via loss reweighting to map function values to corresponding points and applied for offline Black-Box Optimization. In \cite{guogradient}, a kind of Look-Ahead Guidance (LAG) is introduced to preserve the linear structure of data and then used for regularized optimization and global optimization. In \cite{li2024fast}, a diffusion-based training-to-testing (T2T) framework is used to solve new instances in combinatorial optimization while training on historical instances generated by existing algorithms.
	
	Compared with related methods, our work is the first, to the best of our knowledge, to use diffusion models to solve the general chance constrained problems. The key challenge here is the \textbf{lack of direct training data} corresponding to the product distribution of the objective and constraints. We address this through a dedicated data generation stage, followed by conditional training of the score. In contrast, \cite{guogradient} assumes access to a pre-trained unconditional diffusion model and focuses on a restricted linear-Gaussian setting. Unlike classical convex approximation approaches for CCP, our method does not require prior knowledge of the underlying distribution. Instead, we only assume access to samples from it, which makes our approach applicable to broader and more realistic settings.
	
	More specifically, our approach introduces two main innovations:
	\begin{itemize}[leftmargin=10pt]
		\item \textbf{Conditional Training and Applicability Beyond Linear-Gaussian Settings:}
		Unlike \cite{guogradient}, which applies guidance to pre-trained unconditional diffusion models and assumes a linear objective with Gaussian data, our framework involves a dedicated data generation process followed by conditional score training. This enables us to address nonlinear and structurally complex chance-constrained problems, where directly sampling from the feasible region is nontrivial.
		\item \textbf{A New Class of Guidance Derived from Product Distributions:}
		Most existing guided diffusion frameworks follow the general SDE form as follows:
		\begin{equation}
			d\boldsymbol{x}_t=[\boldsymbol{a}(\boldsymbol{x}_t,t)-b(t)^2(\boldsymbol{s}(\boldsymbol{x}_t,t)+\boldsymbol{G}_t )]dt+b(t)d\boldsymbol{\bar{B}}_t.
		\end{equation}
		In our work, we derive two types of guidance terms directly from the product distribution formulation of the target density:
		\begin{itemize}[leftmargin=10pt]
			\item a first-order guidance 
			\begin{equation}
				\boldsymbol{G}_t^{(1)} = - \beta \nabla f(\boldsymbol{x}_t),
			\end{equation}  
			\item a second-order guidance	
			\begin{equation}
				\boldsymbol{G}_t^{(2)} = -\frac{1}{\sigma_{0|t}^2}[\boldsymbol{H}^{-1}[(-\nabla_{\boldsymbol{x}_t}^2 f(\boldsymbol{x}_t) \boldsymbol{x}_t+\nabla f(\boldsymbol{x}_t)) -\frac{1}{\beta \sigma_{0|t}^2}\boldsymbol{\mu}_{0|t}]+\boldsymbol{\mu}_{0|t}],
			\end{equation}
			where the terms are computed based on a learned surrogate for the chance constraint and the posterior mean $\boldsymbol{\mu}_{0|t}$. 
		\end{itemize}
		
		In contrast, \cite{guogradient} introduces a Look-Ahead Guidance term designed for linear objectives:
		\begin{equation}
			\boldsymbol{G}_{t}^{(3)} = -\beta(t)\nabla_{\boldsymbol{x}_t}(y-\boldsymbol{g}^\top \hat{\mathbb{E}}[\boldsymbol{x}_0|\boldsymbol{x}_t])^2,
		\end{equation}
		where $\beta(t)$ and $y$ are tuning parameters, $\boldsymbol{g}$ is the gradient of the linear objective, and $\hat{\mathbb{E}}[\boldsymbol{x}_0|\boldsymbol{x}_t]$ is an approximation of the posterior mean $\boldsymbol{\mu}_{0|t}$ that can be calculated by the score network, i.e., $\hat{\mathbb{E}}[\boldsymbol{x}_0|\boldsymbol{x}_t] = \alpha^{-1}(t)(\boldsymbol{x}_t + h(t)\boldsymbol{s}_{\boldsymbol{\theta}}(\boldsymbol{x}_t,t))$. This approach is effective when the data distribution is Gaussian and the objective is linear, but may degrade under nonlinear or non-Gaussian scenarios. 
		
	\end{itemize}

	\section{Experimental details}
	
	\subsection{Experimental settings}
	
	Our neural network architecture follows the backbone of a U-Net (\cite{ronneberger2015u}) and ResNet (\cite{he2016deep}). We use group normalization (\cite{wu2018group}) to make the implementation simpler. All models use four feature map resolutions with convolutional residual blocks and self-attention blocks (\cite{vaswani2017attention}) per resolution level. Diffusion time $t$ and condition parameter $\rho$ is specified by adding the Transformer sinusoidal position embedding into each residual block.
	
	\begin{figure}[htbp]
		\centering
		\includegraphics[width=0.95\linewidth]{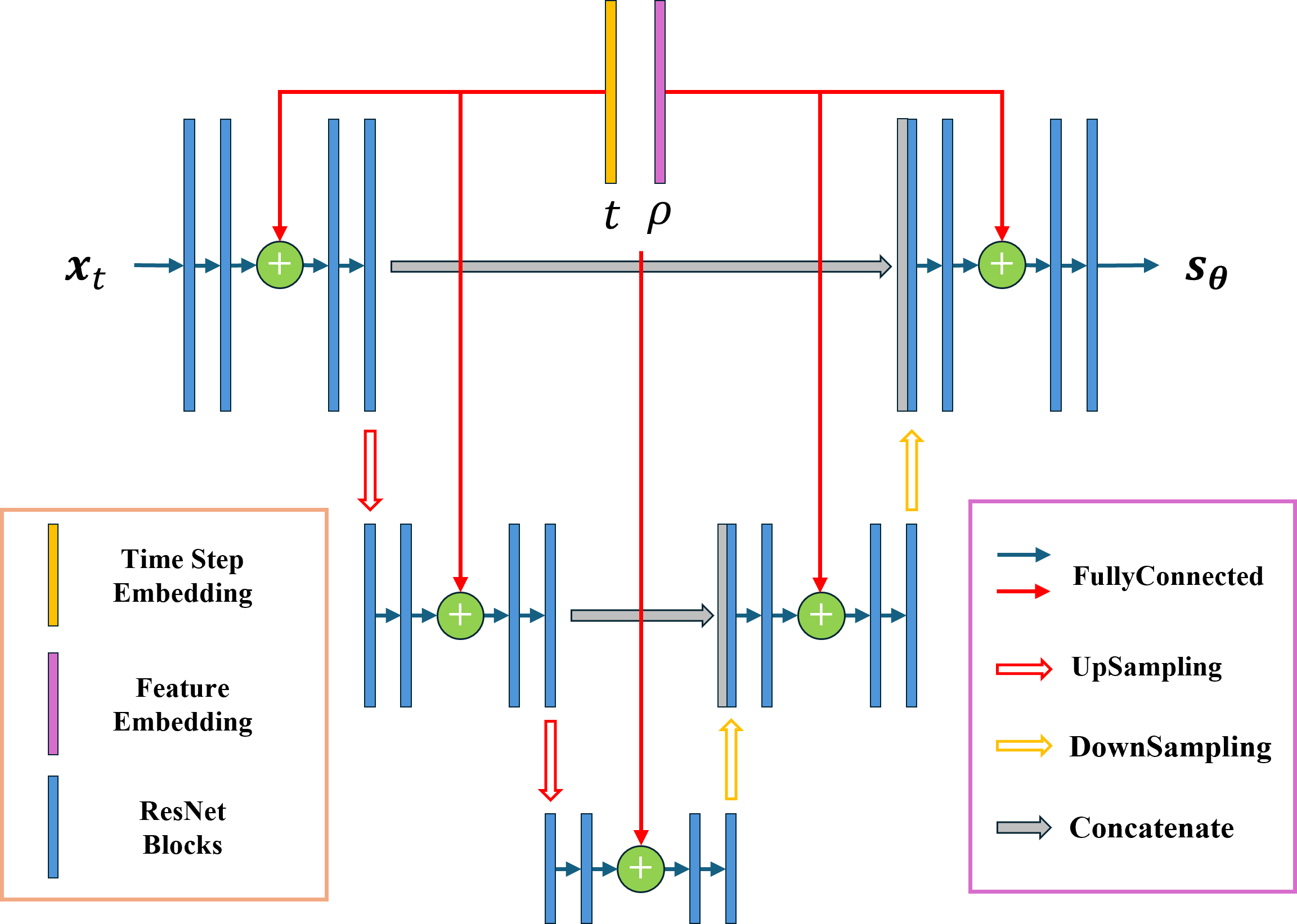}
		\caption{A sketch map for U-Net structure of GGDOpt}
		\label{afig:unet}
	\end{figure}
	
	All models are trained with 4 A800 GPUs. The training durations are approximately 0.4 hours for the linear chance constrained problem and 2 hours for the robust waveform design task. The average sampling times are listed alongside the corresponding experimental results.
	
	We set almost all our hyperparameters as default in (\cite{ho2020denoising}, \cite{guogradient}):
	\begin{itemize}[leftmargin=10pt]
		\item We test the $\eta(t)$ schedule from a set of constant, linear, quadratic and cosine schedules. We set $T=1000$ without a sweep and chose a linear schedule from $\eta(0)=10^{-4}$ to $\eta(T)=0.02$.
		\item We use Adam in our experimentation process and leave the hyperparameters to their standard values. We set the learning rate to $10^{-4}$ without any sweeping.
		\item We set the batch size to 64 for linear chance constrained problem and 128 for robust waveform design. 
	\end{itemize}
	
	To generate the dataset, we utilize CVX (\cite{grant2008cvx}) to solve the restricted problem. For the linear chance constrained problem, we generate $N=1000$ data samples, while $N=10000$ samples for the robust waveform design task. During the sampling with guidance stage, we evaluate both first- and second-order gradient guidance by implementing a DDIM-based technique (\cite{song2020denoising}) with a descaled time step $T'=100$ to accelerate the sampling process.
	
	Our code is available at https://github.com/boyangzhang2000/GGDOpt.
	
	\subsection{Effects of gradient guidance}

\begin{figure}[htbp]
		\centering
		\begin{subfigure}[b]{0.45\textwidth}
			\centering
			\includegraphics[width=\linewidth]{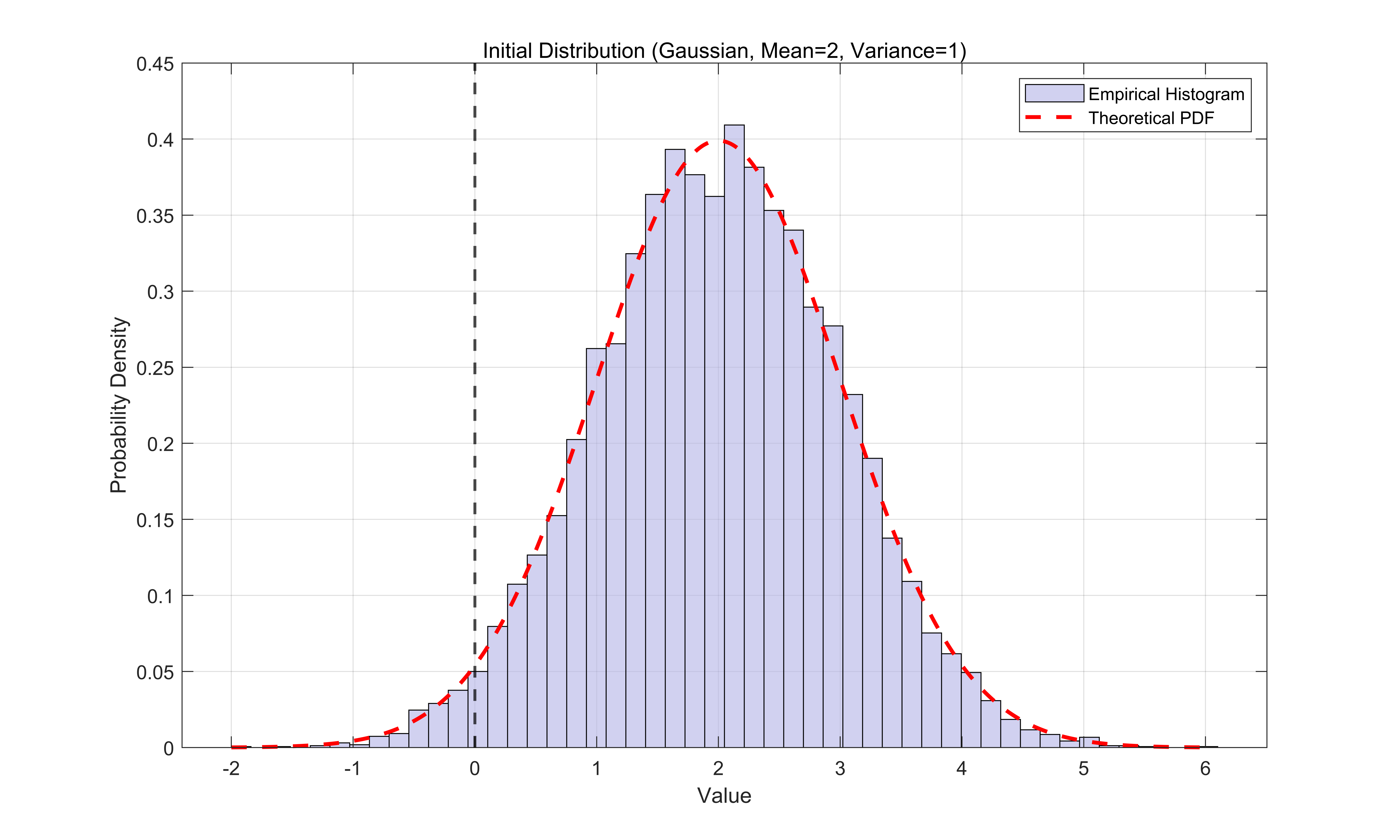}
			\caption{Initial distribution}
		\end{subfigure}
		\hfill
		\begin{subfigure}[b]{0.45\textwidth}
			\centering
			\includegraphics[width=\linewidth]{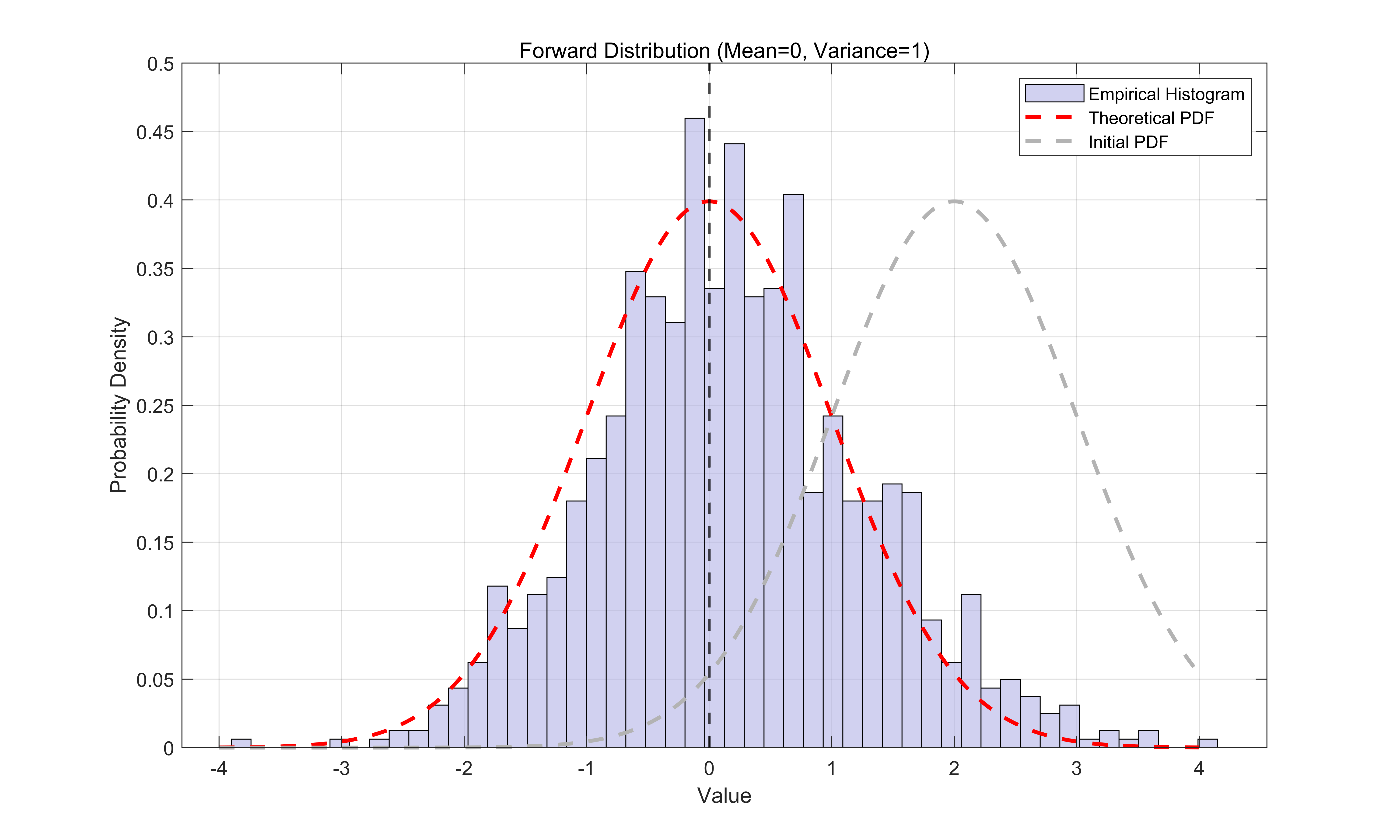}
			\caption{Forward distribution}
		\end{subfigure}	
		\caption{The forward process of GGDOpt.}
		\label{afig:initial_forward}
	\end{figure}

	First, we present an intuitive example illustrating how gradient guidance can steer the sampling trajectory toward the desired target. Specifically, we consider a one-dimensional sampling task where the initial distribution is $x_0\sim\mathcal{N}(2,1)$ and 1000 samples are drawn from it to serve as training data. We set the diffusion time step to $T =1000$ and the resulting forward process of GGDOpt is shown to closely approximate the theoretical distribution $\mathcal{N}(0,1)$ (see Figure \ref{afig:initial_forward}).
	
\begin{figure}[htbp]
		\centering
		\begin{subfigure}[b]{0.32\textwidth}
			\centering
			\includegraphics[width=\linewidth]{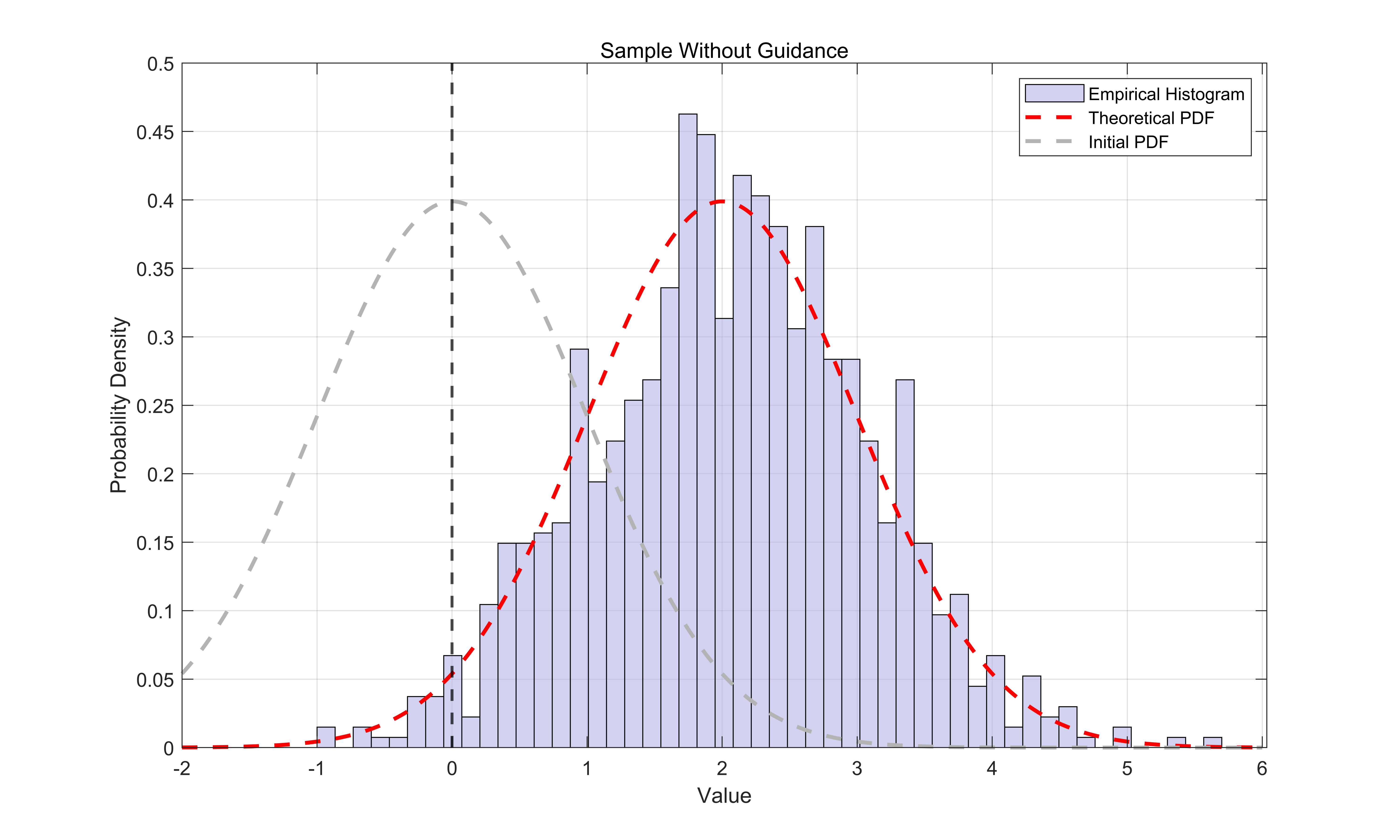}
			\caption{Sampling without Guidance}
		\end{subfigure}
		\hfill
		\begin{subfigure}[b]{0.32\textwidth}
			\centering
			\includegraphics[width=\linewidth]{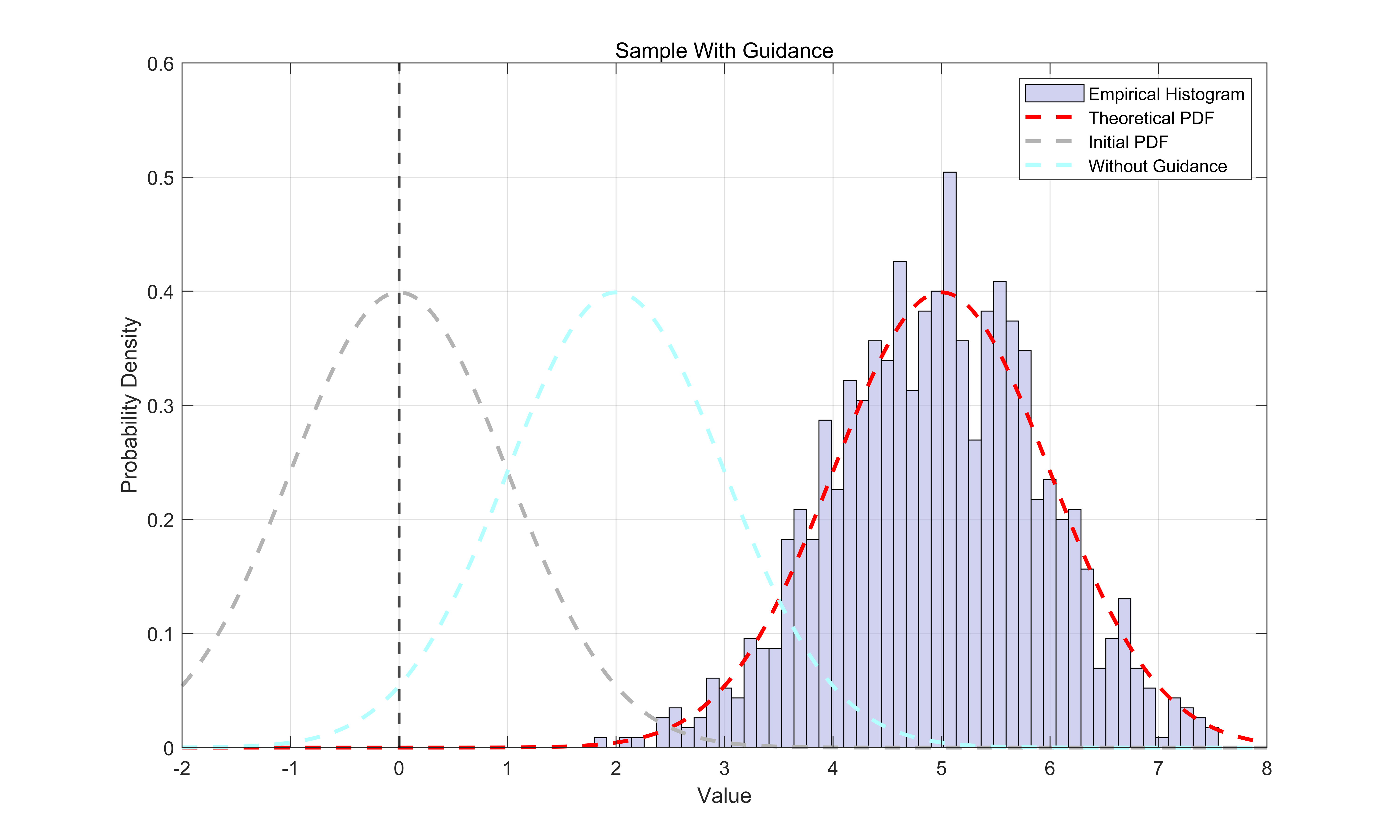}
			\caption{First-order Guidance}
		\end{subfigure}	
		\hfill
		\begin{subfigure}[b]{0.32\textwidth}
			\centering
			\includegraphics[width=\linewidth]{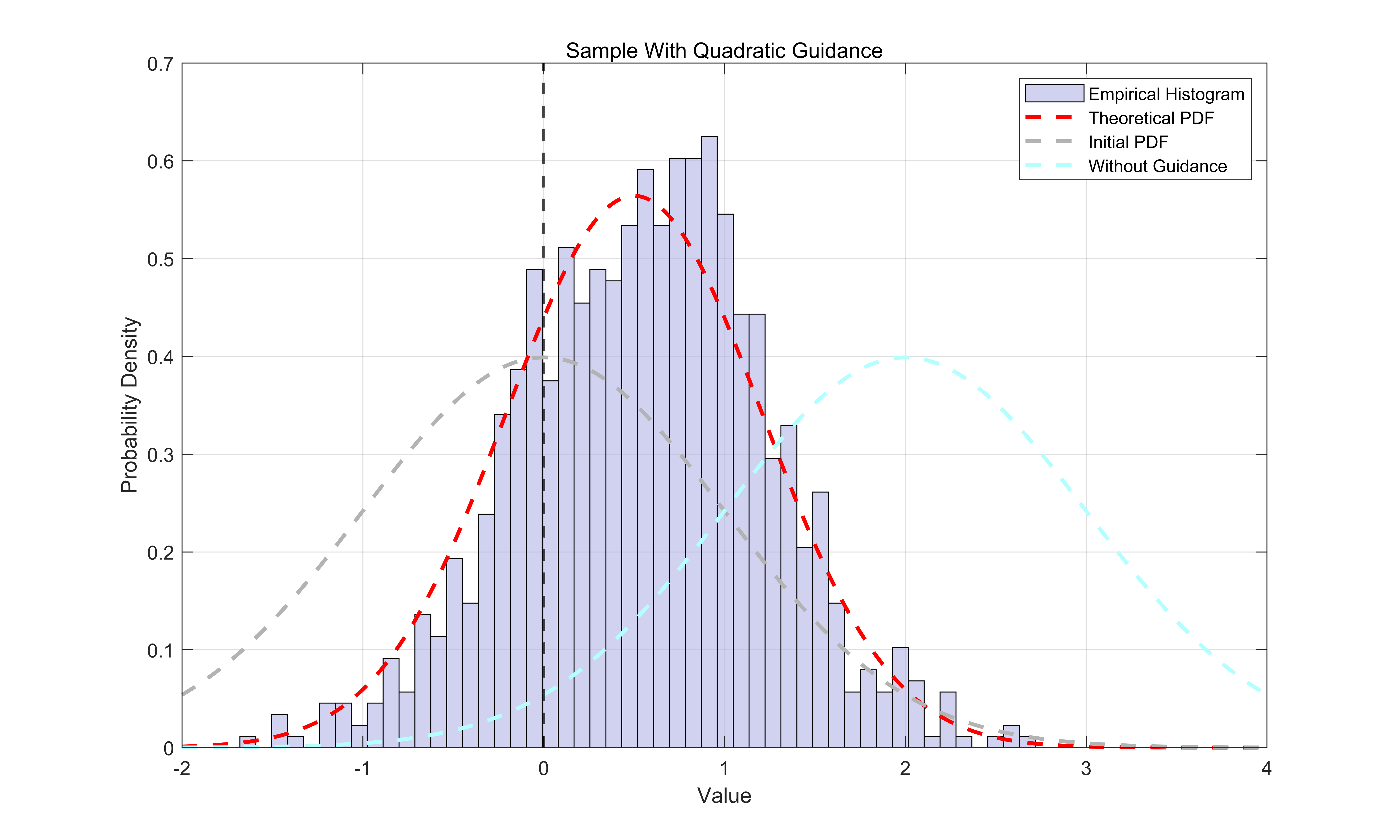}
			\caption{Second-order Guidance}
		\end{subfigure}	
		\caption{The sampling process of GGDOpt}
		\label{afig:sampling}
	\end{figure}

	Next, we compare different sampling strategies: without guidance, first-order gradient guidance, and second-order gradient guidance. Theoretically, under Gaussian assumptions, first-order gradient guidance alters only the mean of the end distribution, whereas second-order gradient guidance affects both the mean and the variance. For each method, we generate 1000 samples and the corresponding sampling results are presented in Figure \ref{afig:sampling}.

	Experimental results demonstrate that, in the absence of guidance, the sampling process shifts the distribution from the prior $\mathcal{N}(0,1)$ back to the initial distribution $\mathcal{N}(2,1)$, as expected. When applying first-order gradient guidance with $\beta = 3$, the distribution transitions from the prior $\mathcal{N}(0,1)$ to the guided distribution $\mathcal{N}(5,1)$, indicating a change in the mean while preserving the variance. In contrast, with second-order gradient guidance and $\beta = 1$, the distribution is modified to $\mathcal{N}(1/2,1/2)$, reflecting changes in both mean and variance. These results confirm that GGDOpt effectively directs the sampling process to the desired end distribution. Furthermore, setting $T=1000$ is sufficient to eliminate the limited time length error.

	\subsection{Additional experimental results}
	
	\subsubsection{Linear chance constrained problem}
	
	Consider the following linear chance constrained problem
	\begin{equation}
		\label{aeq:linear}
		\begin{aligned}
			\min_{\bm{x}\in\mathbb{R}^n}\quad&\frac{1}{2}\bm{x}^{\top}\bm{x}+\bm{b}^{\top}\bm{x}\\
			\mathrm{s.t.}\quad&\text{Prob}_{\bm{c}\sim p_{\bm{c}}}\{\bm{c}^{\top}\bm{x}+d\geq0\}\geq 1-\rho,
		\end{aligned}
	\end{equation}
	where the uncertain parameter follows a Gaussian distribution $p_{\bm{c}} = \mathcal{N}(\bm{c};\bar{\bm{c}},\bm{I})$ and the hyperparameters $(\bm{b},\bar{\bm{c}},d,\rho)$ are selected from a predefined test set. 
	
	For any $\rho<0.5$, the linear chance constraint can be expressed as
	\begin{equation}
		-\Phi^{-1}(\rho)\|\bm{x}\|_2-(\bar{\bm{c}}^{\top}\bm{x}+d)\leq0,
	\end{equation}
	where $\Phi$ denotes the standard Gaussian cumulative distribution function. Then the linear chance constrained problem (\ref{aeq:linear}) can be reformulated as the following second-order cone program:
	\begin{equation}
		\begin{aligned}
			\min_{\bm{x}}\quad&\frac{1}{2}\bm{x}^{\top}\bm{x}+\bm{b}^{\top}\bm{x}\\
			\mathrm{s.t.}\quad&-\Phi^{-1}(\rho)\|\bm{x}\|_2-(\bar{\bm{c}}^{\top}\bm{x}+d)\leq0,
		\end{aligned}
	\end{equation}
	which is solved using CVX (\cite{grant2008cvx}). In practice, we assume the distribution $p_{\bm{c}}$ is unknown and only 100 samples are available. To generate training data, we solve the restricted version of the problem for $N=1000$ values of $z$ linearly spaced in the interval $\left[0, 0.5\right]$.
	
\begin{figure}[htbp]
		\centering
		\begin{subfigure}[b]{0.32\textwidth}
			\centering
			\includegraphics[width=\linewidth]{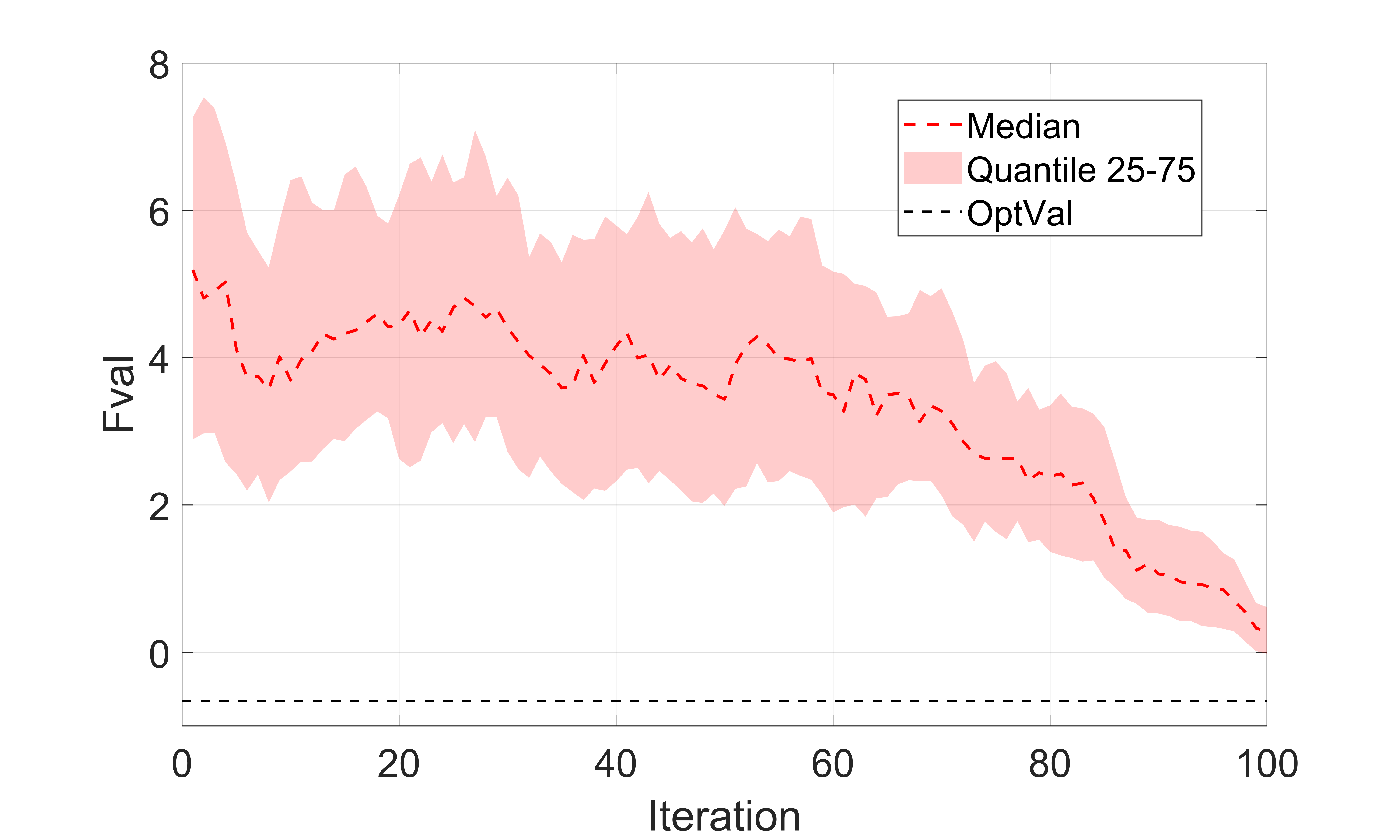}
			\caption{Without Guidance}
		\end{subfigure}
		\begin{subfigure}[b]{0.32\textwidth}
			\centering
			\includegraphics[width=\linewidth]{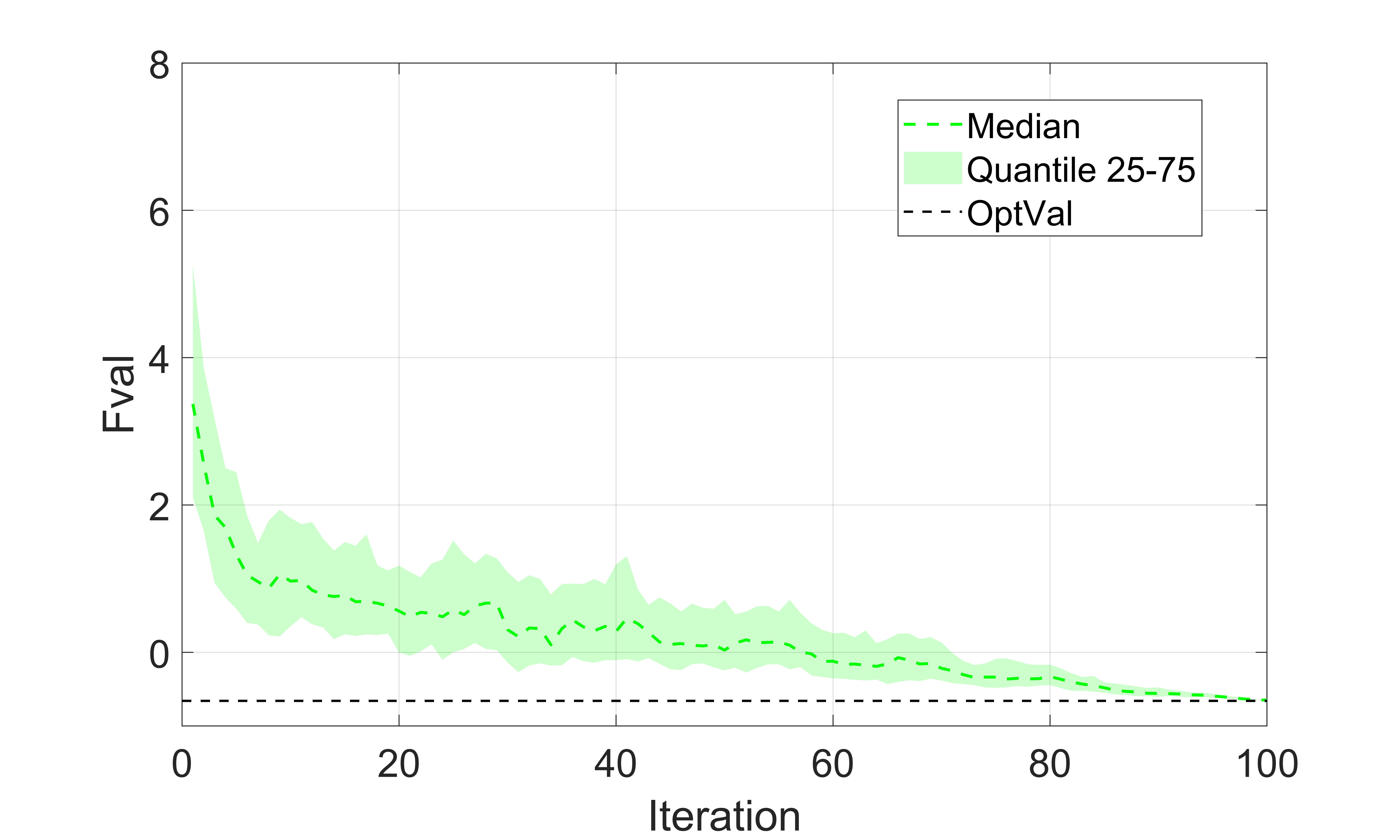}
			\caption{First-order Guidance}
		\end{subfigure}
		\begin{subfigure}[b]{0.32\textwidth}
			\centering
			\includegraphics[width=\linewidth]{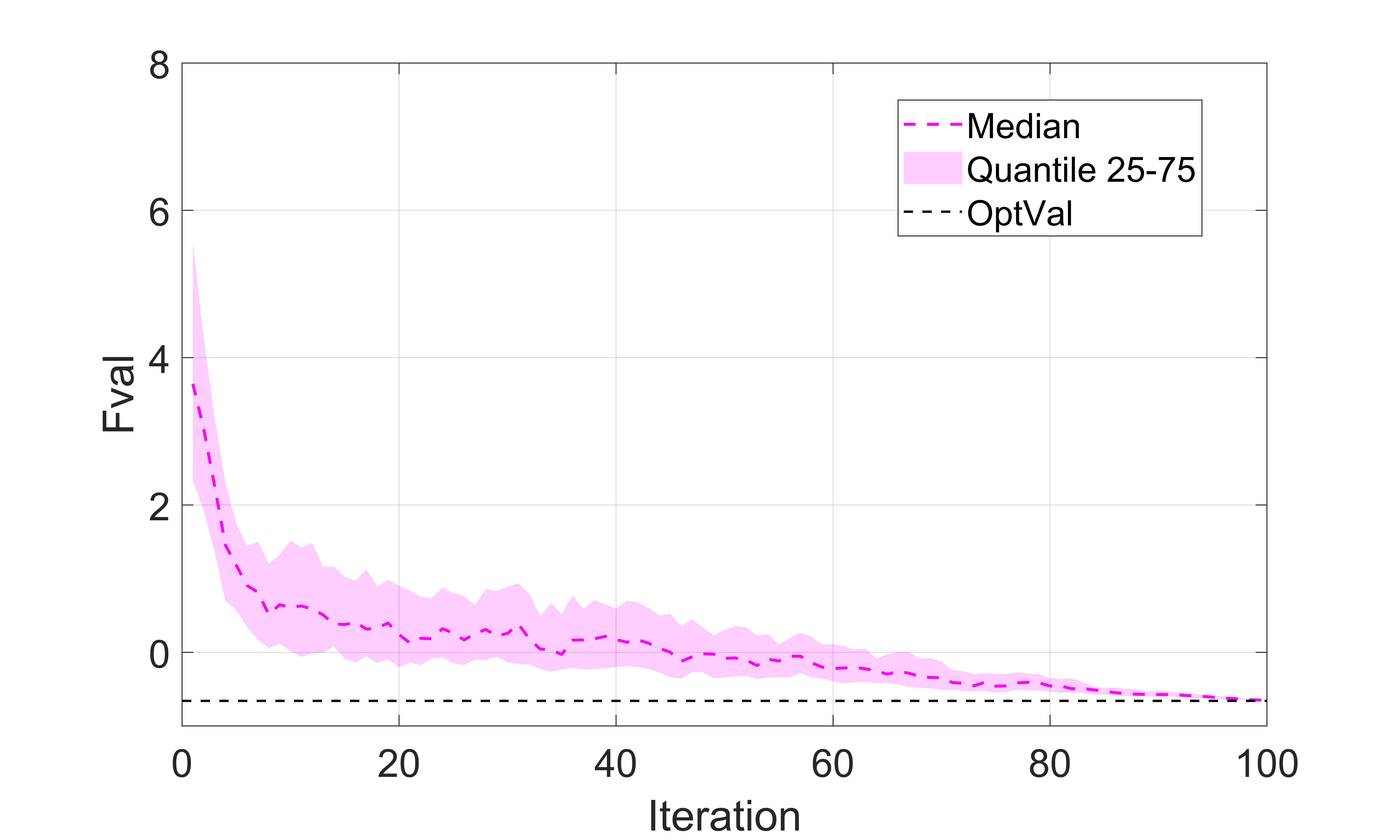}
			\caption{Second-order Guidance}
		\end{subfigure}
		\caption{Sampling process visualization of GGDOpt with median and quantiles}
		\label{afig:error}
	\end{figure}

	We evaluate the performance of the proposed GGDOpt framework by comparing it with several SAA approaches, using the CVX-based solutions as performance benchmarks. Each algorithm (excluding CVX) is run 100 times and objective values are reported after projecting the solutions onto the feasible set. Experimental results for the case with $n=8, \bm{b} = \bar{\bm{c}} = (1, 1, \ldots, 1), d=1, \rho=0.1$ are summarized in Table \ref{atab:linear}, and the sampling process characterized by median and quantiles are provided in Figure \ref{afig:error} to show the stability and fast convergence of GGDOpt.
	
	\begin{table}[htbp]
		\caption{Comparison results on the linear chance constrained problem (\ref{aeq:linear})}
		\label{atab:linear}
		\renewcommand{\arraystretch}{1.8}
		\resizebox{\textwidth}{!}{
			\begin{tabular}{ccccccc}
				\hline
				Method                   & FvalMean & FvalStd & FvalMedian & FvalQuan25 & FvalQuan75 & Runtime \\ \hline
				\begin{tabular}[c]{@{}c@{}}SOC\_CVX\\ (\cite{grant2008cvx})\end{tabular}                     & \textbf{-0.6586}  & 0       & -0.6586    & -0.6586        & -0.6586        & 0.3214  \\ \hline
				\begin{tabular}[c]{@{}c@{}}SAA\_MIP\\ (\cite{pagnoncelli2009sample})\end{tabular}                & -0.6281  & 0.0157  & -0.6318    & -0.6396        & -0.6184        & 15.4502 \\
				\begin{tabular}[c]{@{}c@{}}SAA\_CVaR\\ (\cite{nemirovski2007convex})\end{tabular}                & -0.5893  & 0.0248  & -0.5869    & -0.6021        & -0.5702        & 0.3063 \\
				\begin{tabular}[c]{@{}c@{}}SAA\_SNSCO\\ (\cite{zhou20240})\end{tabular}                & 0.8051   & 3.4014  & -0.6371    & -0.6469        & -0.6019        & 0.2793 \\
				\begin{tabular}[c]{@{}c@{}}SAA\_PDCA\\ (\cite{wang2023proximal})\end{tabular}                   & -0.6389  & 0.0314  & -0.6408    & \textbf{-0.6566}        & -0.6185        & 0.6276
				\\ \hline
				\begin{tabular}[c]{@{}c@{}}GGDOpt\\ (Without Guidance)\end{tabular}     & 0.3481   & 0.5486  & 0.2798     & -0.0181        & 0.6142         & 0.0465  \\
				\begin{tabular}[c]{@{}c@{}}GGDOpt\\ (First-order)\end{tabular}       & -0.6483  & 0.0051  & -0.6488    & -0.6525        & -0.6454        & \textbf{0.0486}  \\
				\begin{tabular}[c]{@{}c@{}}GGDOpt\\ (Second-order )\end{tabular}    & \textbf{-0.6491}  & \textbf{0.0056}  & \textbf{-0.6503}    & -0.6531        & \textbf{-0.6474}        & \textbf{0.0507}  \\ \hline
			\end{tabular}
		}
	\end{table}
	
	Based on the results presented in Table \ref{atab:linear}, we observe that SOC\_CVX is capable of exactly identifying the global minimizer of the convexified problem, given full knowledge of the underlying probability distribution. In contrast, SAA-based methods rely solely on sampled realizations and thus yield approximate solutions. Among them, SAA\_MIP requires solving a large-scale mixed-integer optimization problem, which is computationally expensive. While SAA\_SNSCO demonstrates rapid convergence to optimal solutions in most cases, its performance degrades under worst realizations of $\bm{h}$, occasionally converging to sub-optimal solutions. This leads to strong median performance but instability in statistical results.
	
	Compared with the SAA methods, our proposed GGDOpt demonstrates superior stability and yields higher-quality solutions, while also significantly reducing computational overhead.
	
\begin{figure}[htbp]
		\centering
		\includegraphics[width=0.9\linewidth]{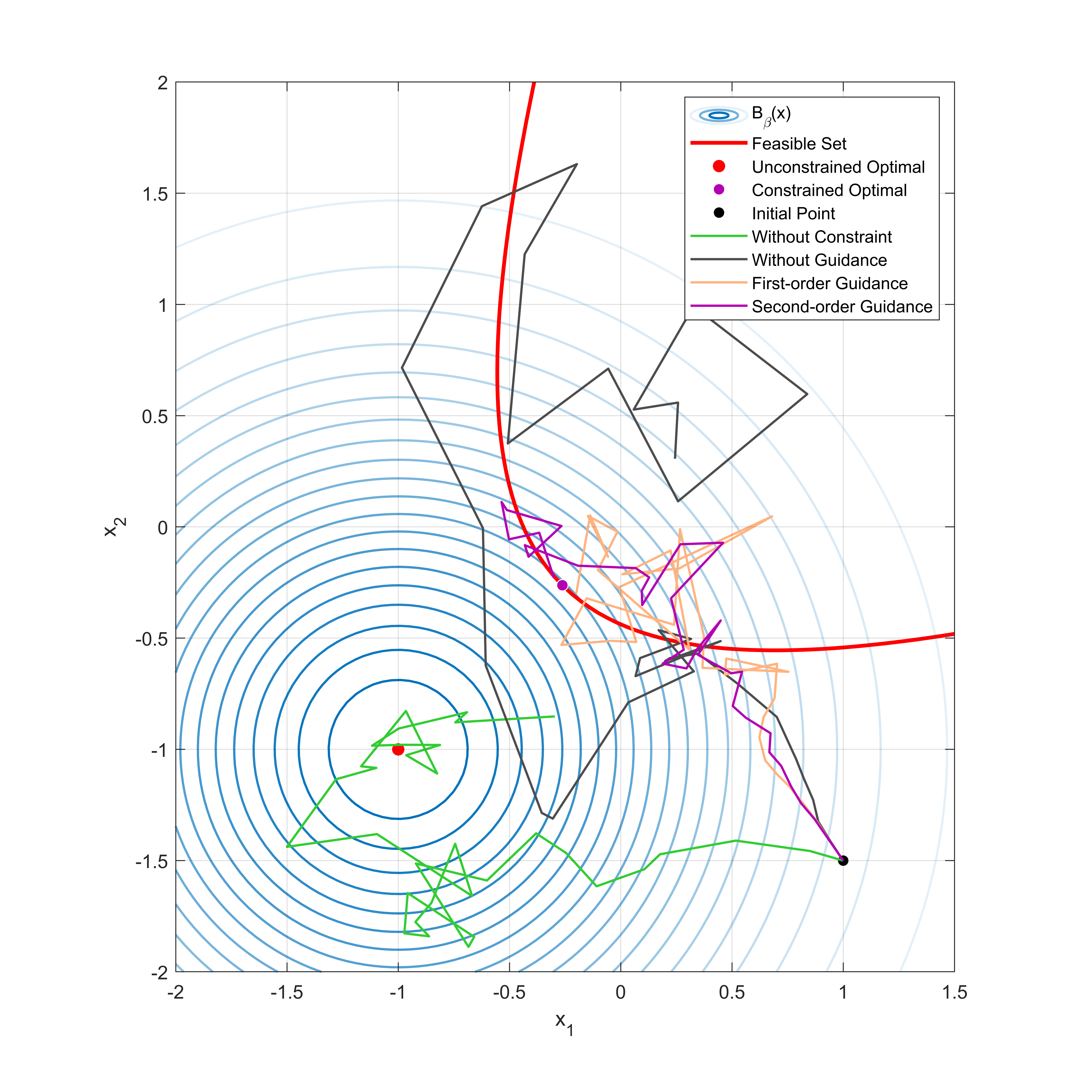}
		\caption{Sampling path of various methods}
		\label{afig:path}
	\end{figure}

\begin{figure}[htbp]
		\centering
		\begin{subfigure}[b]{0.48\textwidth}
			\centering
			\includegraphics[width=\linewidth]{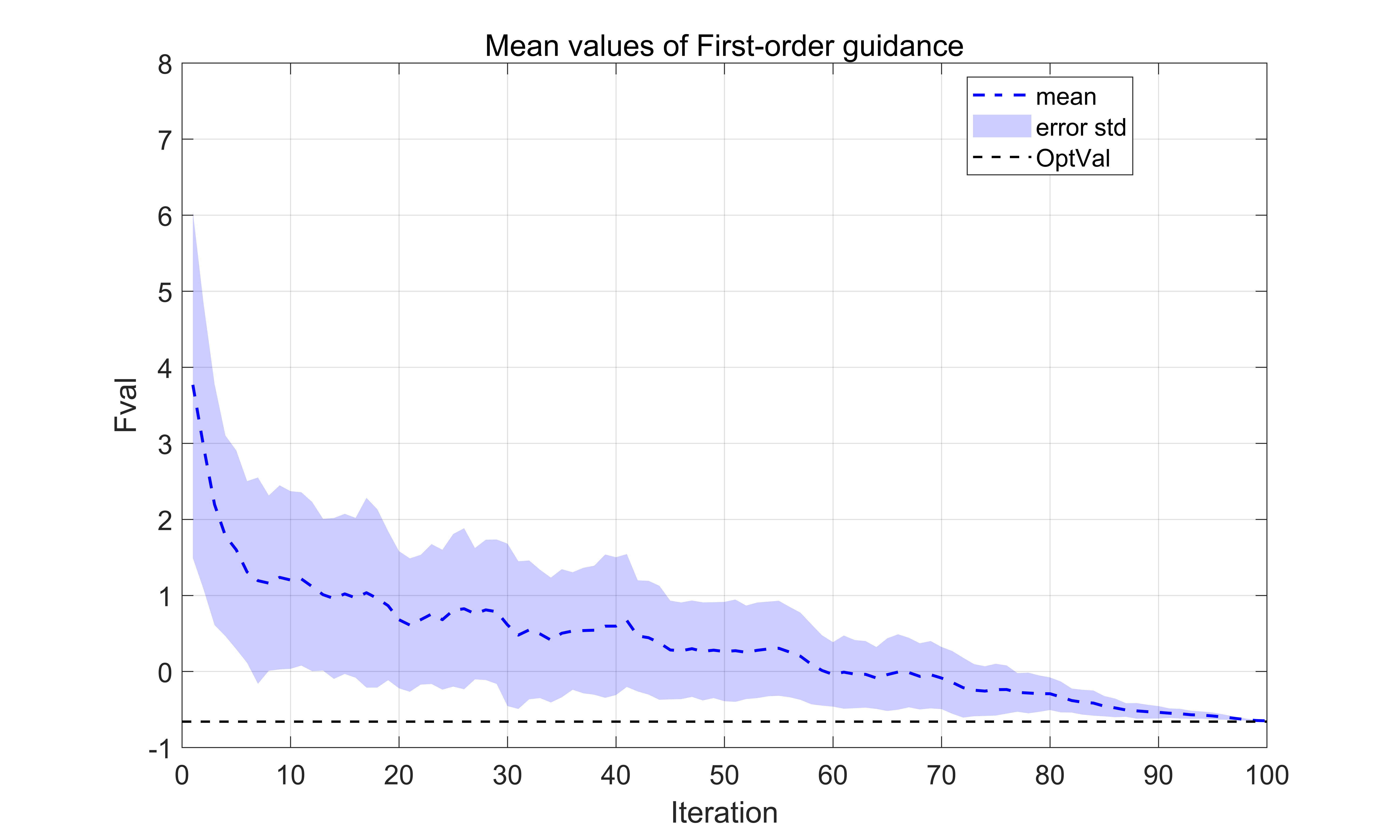}
		\end{subfigure}
		\begin{subfigure}[b]{0.48\textwidth}
			\centering
			\includegraphics[width=\linewidth]{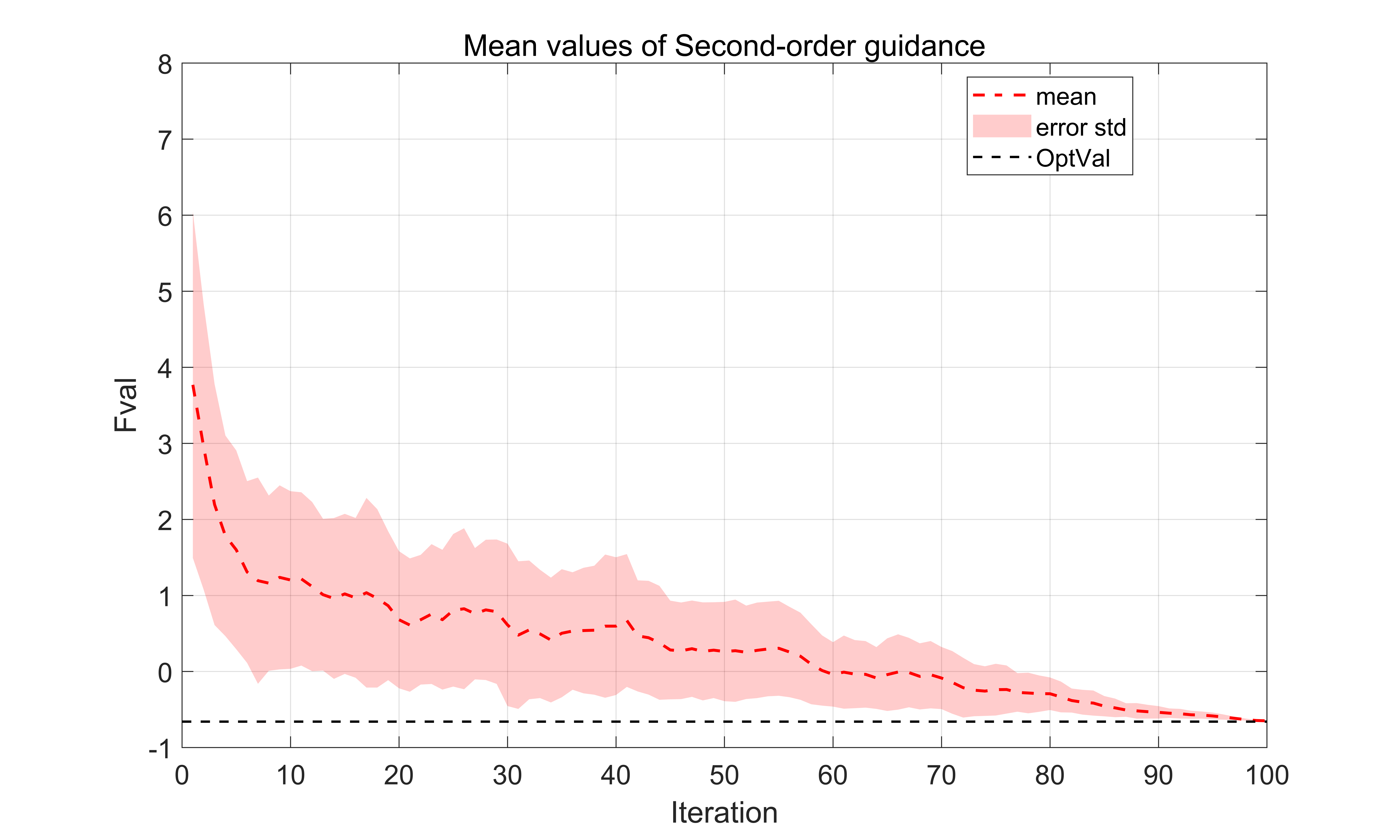}
		\end{subfigure}
		\caption{Convergence of the objective values for first- and second-order gradient guidance}
		\label{afig:conv}
	\end{figure}

	To provide an intuitive understanding of the sampling behavior in GGDOpt, we illustrate a representative sampling trajectory of different methods in Figure \ref{afig:path}. The results show that, without constraint, the sampling process will concentrate on the global minimizer of objective function. Under the influence of constraint, the samples will fall into the feasible set and gradient guidance will lead the sampling path toward the direction with lower function value. The corresponding iterations of the objective values for first-order gradient guidance and second-order gradient guidance are shown in Figure \ref{afig:conv}.
	
	Furthermore, we demonstrate that GGDOpt is capable of producing high-quality solutions across a range of values for the risk parameter $\rho$. Specifically, we vary $\rho$ from 0.05 to 0.30 while keeping all other experimental settings fixed. The corresponding results are reported in Table \ref{atab:rhos}.
	
	To further illustrate the efficiency and robustness of GGDOpt, we evaluate its performance under varying problem dimensions. In particular, we vary the number of decision variables $n$ from 2 to 1024, using the corresponding CVX solutions as performance benchmarks (normalized to 100\%). The comparative performance of GGDOpt under first-order and second-order gradient guidance is summarized in Table \ref{atab:ns}.
	
	\begin{table}[htbp]
		\caption{Comparison results on the linear chance constrained problem (\ref{aeq:linear}) with different $\rho$}
		\label{atab:rhos}
		\renewcommand{\arraystretch}{1.8}
		\centering
			\begin{tabular}{lcccccc}
				\hline
				Method               & $\rho=0.05$ & $\rho=0.10$ & $\rho=0.15$ & $\rho=0.20$ & $\rho=0.25$ & $\rho=0.30$ \\ \hline
				SOC\_CVX             & \textbf{-0.6073}  & \textbf{-0.6585}  & \textbf{-0.6983}  & \textbf{-0.7335}  & \textbf{-0.7667}  & \textbf{-0.7991}  \\ \hline
				SAA\_MIP             & -0.5729  & -0.6229  & -0.6312  & -0.6883  & -0.7053  & -0.7351  \\
				SAA\_CVaR            & -0.5787  & -0.5893  & -0.6378  & -0.6583  & -0.6769  & -0.6892  \\
				SAA\_SNSCO           & 1.0632   & 0.8051   & -0.2408  & -0.4760  & -0.5635  & -0.6617  \\
				SAA\_PDCA            & -0.5730  & -0.6283  & -0.6665  & -0.6957  & -0.7279  & -0.7639  \\ \hline
				GGDOpt\_First-order  & -0.5955  & -0.6483  & -0.6828  & -0.7032  & -0.7284  & -0.7603  \\
				GGDOpt\_Second-order & \textbf{-0.6040}  & \textbf{-0.6491}  & \textbf{-0.6944}  & \textbf{-0.7130}  & \textbf{-0.7498}  & \textbf{-0.7817}  \\ \hline
			\end{tabular}
	\end{table}
	
	\begin{table}[htbp]
		\caption{Comparison results on the linear chance constrained problem (\ref{aeq:linear}) with different $n$}
		\label{atab:ns}
		\renewcommand{\arraystretch}{1.8}
		\centering
			\begin{tabular}{lccccc}
				\hline
				\multicolumn{2}{c}{Method}     & \multirow{2}{*}{SOC\_CVX} & \multirow{2}{*}{SAA\_PDCA} & GGDOpt        & GGDOpt         \\
				\multicolumn{2}{c}{$\rho=0.1$}    &                           &                            & (First-order) & (Second-order) \\ \hline
				\multirow{2}{*}{$n=2$}    & fval & 100.00\%                  & 99.39\%                    & 99.87\%       & 100.00\%       \\
				& time & 100.00\%                  & 111.22\%                   & 7.26\%        & 7.31\%         \\ \hline
				\multirow{2}{*}{$n=4$}    & fval & 100.00\%                  & 96.86\%                    & 99.73\%       & 99.89\%        \\
				& time & 100.00\%                  & 118.26\%                   & 9.86\%        & 10.02\%        \\ \hline
				\multirow{2}{*}{$n=8$}    & fval & 100.00\%                  & 95.67\%                    & 98.44\%       & 98.56\%        \\
				& time & 100.00\%                  & 142.25\%                   & 15.12\%       & 15.77\%        \\ \hline
				\multirow{2}{*}{$n=16$}   & fval & 100.00\%                  & 90.42\%                    & 98.47\%       & 99.70\%        \\
				& time & 100.00\%                  & 160.56\%                   & 15.10\%       & 16.06\%        \\ \hline
				\multirow{2}{*}{$n=128$}  & fval & 100.00\%                  & 95.88\%                    & 98.72\%       & 99.89\%        \\
				& time & 100.00\%                  & 198.24\%                   & 23.63\%       & 25.93\%        \\ \hline
				\multirow{2}{*}{$n=1024$} & fval & 100.00\%                  & 93.19\%                    & 97.91\%       & 99.34\%        \\
				& time & 100.00\%                  & 516.16\%                   & 34.11\%       & 37.64\%        \\ \hline
			\end{tabular}
	\end{table}

	The results above demonstrate that GGDOpt effectively solves the linear chance constrained problem across varying parameter settings. Moreover, it exhibits significantly higher computational efficiency compared to alternative approaches.
	
	\subsubsection{Computational cost}
	
	Regarding the computational cost and evaluation, we test the linear chance constrained problem with $n=8$ and repeat 100 times to calculate the empirical mean of the objective value (fmean), the empirical standard deviation (fstd), and the average run time (time). The results are summarized in the following Table \ref{atab:cost}.
	
	\begin{table}[htbp]
		\caption{Computational cost of the proposed methods.}
		\label{atab:cost}
		\renewcommand{\arraystretch}{1.2}
		\centering
			\begin{tabular}{cccccc}
				\hline
				\multirow{2}{*}{Method} & \multirow{2}{*}{SOC\_CVX}  & \multirow{2}{*}{GGDOpt (First-order)} & \multicolumn{3}{c}{GGDOpt (Second-order)}                                                               \\ \cline{4-6} 
				&                           &                                                          & $\beta=0.1$ & $\beta=1$ & $\beta=10$ \\ \hline
				fmean                & -0.6586                               & -0.6483                                                & -0.6341                  & -0.6548                & -0.6585                 \\
				fstd                 & 0                                      & 0.0051                                            & 5.6726e-3                & 2.5112e-05             & 2.2329e-08              \\
				time                 & 0.3214                               & 0.0486                                             & 0.0569                   & 0.0527                 & 0.0541                  \\ \hline
			\end{tabular}
	\end{table}
	
	As observed in Table \ref{atab:cost}, the second-order method achieves lower objective values compared to the first-order method and its performance closely matches the optimal solution obtained by SOC\_CVX. Moreover, the second-order method leads to significantly lower standard deviations, particularly as $\beta$ increases.
	
	We also provide the costs of three stages for the linear chance constraint problem. For each $n$, we generate 1000 data in the training stage. During sampling, we execute 100 times of reverse process to analyze the stability of GGDOpt. The total time costed in hour is shown in Table \ref{atab:time}.
	
	\begin{table}[htbp]
		\caption{Computational time of three stages (in hours).}
		\centering
		\label{atab:time}
		\renewcommand{\arraystretch}{1.2}
			\begin{tabular}{ccccc}
				\hline
				\multicolumn{2}{c}{Stages}                       & $n=8$    & $n=16$   & $n=128$  \\ \hline
				\multicolumn{2}{c}{Data generating time}        & 0.03     & 0.06     & 0.11     \\
				\multicolumn{2}{c}{Training time}                & 0.53     & 0.96     & 11.64    \\ \hline
				\multirow{2}{*}{Total sampling time} 
				& First-order                                     & 0.0013   & 0.0017   & 0.0057   \\ 
				& Second-order                                    & 0.0014   & 0.0018   & 0.0063   \\ \hline
			\end{tabular}
	\end{table}
	
	Furthermore, our experiments indicate that increasing the quantity of training data alone does not guarantee better performance. Instead, high-quality samples closer to the true optimal solutions are the key of effective guided sampling.
	
	\subsubsection{Variance schedule}

	While Tweedie’s formula theoretically provides both the posterior mean and covariance, $\boldsymbol{\Sigma}_{0|t}=(1-\bar{\alpha}_t)(\boldsymbol{I}+(1-\bar{\alpha}_t)\nabla^2\log p(\boldsymbol{x}_t))$, computing the covariance requires evaluating the Hessian of $\log p(\boldsymbol{x})$.
	
	In our framework, the score function $\boldsymbol{s}_{\boldsymbol{\theta}}$ is parameterized by a neural network, and computing its second derivatives involves backpropagation through the network’s Jacobian, which is computationally expensive, especially in high dimensions.
	
	To strike a balance between performance and efficiency, we choose to treat the covariance as a tunable constant. This introduces an approximation, but as shown in Table \ref{atab:var}, this achieves comparable objective values to the fully Tweedie-based method, while reducing runtime by more than an order of magnitude. These results confirm that using a fixed variance can be a practical and robust alternative.
	
	\begin{table}[htbp]
		\caption{Experimental results with different variance schedules.}
		\label{atab:var}
		\renewcommand{\arraystretch}{1.4}
		\centering
			\begin{tabular}{ccccccc}
				\hline
				\multirow{2}{*}{$n=8,\rho=0.1$} & \multirow{2}{*}{Tweedie's $\boldsymbol{\Sigma}$} & \multicolumn{5}{c}{GGDOpt (Second-order)}                \\ \cline{3-7} 
				&                          & $\sigma=0.01$ & $\sigma=0.02$ & $\sigma=0.1$ & $\sigma=1$ & $\sigma=10$ \\ \hline
				fval                           & -0.6571                  & -0.6471    & -0.6457    & -0.6545   & -0.6320 & -0.6049  \\
				time                           & 1.0984                   & 0.0491     & 0.0496     & 0.0493    & 0.0492  & 0.0493   \\ \hline
			\end{tabular}
	\end{table}

	\subsubsection{Guidance term}
	
	Our experimental results in Table \ref{atab:lag} further demonstrate that for the chance constrained programming, the proposed GGDOpt consistently outperforms the Look-Ahead Guidance from \cite{guogradient} in terms of both objective value (fval) and computational efficiency (sampling time).
	
	\begin{table}[htbp]
		\caption{Comparison results with Look-Ahead Guidance \cite{guogradient}.}
		\label{atab:lag}
		\renewcommand{\arraystretch}{1.1}
		\resizebox{\textwidth}{!}{
			\begin{tabular}{ccccccccc}
				\hline
				\multirow{2}{*}{Method   ($\rho=0.1$)}                                                & \multicolumn{2}{c}{$n=2$} & \multicolumn{2}{c}{$n=4$} & \multicolumn{2}{c}{$n=8$} & \multicolumn{2}{c}{$n=16$} \\ \cline{2-9}
				& fval        & time      & fval        & time      & fval        & time      & fval        & time       \\ \hline
				SOC\_CVX                                                                           & -0.4558     & 0.2148    & -0.5630     & 0.2415    & -0.6586     & 0.3214    & -0.7394     & 0.4067     \\ \hline
				\begin{tabular}[c]{@{}c@{}}GGDOpt\\      (First-order)\end{tabular}                   & -0.4552     & 0.0156    & -0.5615     & 0.0238    & -0.6483     & 0.0486    & -0.7281     & 0.0614     \\
				\begin{tabular}[c]{@{}c@{}}GGDOpt\\      (Second-order)\end{tabular}                  & -0.4558     & 0.0157    & -0.5624     & 0.0242    & -0.6491     & 0.0507    & -0.7372     & 0.0653     \\ \hline
				\begin{tabular}[c]{@{}c@{}}LAG\\      \cite{guogradient}\end{tabular} & -0.4460     & 0.0329    & -0.5181     & 0.0738    & -0.5783     & 0.1127    & -0.6584     & 0.1436     \\ \hline
			\end{tabular}
		}
	\end{table}
	
	As shown in the table, our proposed GGDOpt consistently achieves lower objective values and the performance gap between GGDOpt and Look-Ahead Guidance increases with the problem dimension $n$. In terms of computational efficiency, GGDOpt is approximately $2\times$ faster than the Look-Ahead Guidance across all problem sizes. This performance gain stems from the computational overhead of \cite{guogradient}, where computing the guidance term $\boldsymbol{G}_{t}^{(3)}$ requires backpropagation through the score network to obtain the gradient of the posterior mean $\mathbb{E}[\boldsymbol{x}_0|\boldsymbol{x}_t]$ with respect to $\boldsymbol{x}_t$. In contrast, our first- and second-order guidance terms are derived analytically and thus do \textbf{not require any additional gradient computations through the network}, making our method more efficient and scalable.

	\subsubsection{VaR-constrained mean–variance portfolio selection problem}

	Consider a VaR-constrained mean–variance portfolio selection problem, which aims to minimize the risk while pursuing a targeted level of returns with probability at least $1-\rho$ (\cite{wang2023proximal}). Let $\boldsymbol{\mu}\in\mathbb{R}^n$ and $\boldsymbol{\Sigma}\in\mathbb{R}^{n\times n}$ denote the expectation and covariance matrix of the returns of $n$ risky assets, and $\gamma\in\mathbb{R}_+$ denote the risk aversion factor. Let $\boldsymbol{x}\in\mathbb{R}^n$ denote the allocation vector. Then this problem is formulated as follows:
	
	\begin{equation}
		\begin{aligned}
			\min_{\boldsymbol{x}\in\mathbb{R}^n}&\quad\gamma \boldsymbol{x}^{\top}\boldsymbol{\Sigma}\boldsymbol{x}-\boldsymbol{\mu}^{\top}\boldsymbol{x} \\
			\text{s.t. }&\quad\text{Prob}_{\boldsymbol{\xi}}\{\boldsymbol{\xi}^{\top}\boldsymbol{x}\geq R\}\geq 1-\rho,
		\end{aligned}
	\end{equation}
	where $R\in\mathbb{R}_+$ is a prespecified level on the return. We use 2523 daily return data of 435 stocks included in Standard \& Poor’s 500 Index between March 2006 and March 2016 and set $R=0.02\%$ and $\gamma=2$. Some results are shown in Table \ref{atab:port}:
	
	\begin{table}[htbp]
		\caption{Comparison results of the VaR-constrained mean-variance portfolio selection problem.}
		\label{atab:port}
		\renewcommand{\arraystretch}{1.2}
		\resizebox{\textwidth}{!}{
			\begin{tabular}{cccccccc}
				\hline
				($\rho, n$)                      & Metric & MIP                             & ALDM    & PDCA                          & LAG & GGDOpt (First) & GGDOpt (Second) \\ \hline
				& fval   & -0.0951 & -0.0723 & -0.0917                       & -0.0936       & -0.0904        & -0.0946         \\
				& time   & 15.58                          & 2.418   & 4.602                         & 0.9433        & 0.3768         & 0.4071          \\
				\multirow{-3}{*}{(0.05, 100)} & prob   & 0.8600                         & 0.8666  &  0.9700 & 0.8467        & 0.9200         & 0.8933          \\ \hline
				& fval   &  -0.0874 & -0.0750 & -0.0814                       & -0.0859       & -0.0827        & -0.0867         \\
				& time   & 204.2                          & 66.68   & 93.42                         & 2.7570        & 1.2732         & 1.3559          \\
				\multirow{-3}{*}{(0.05, 400)} & prob   & 0.9066                         & 0.8308  &  0.9891 & 0.8933        & 0.9533         & 0.9267          \\ \hline
				& fval   & -0.0951 & -0.0721 & -0.0856                       & -0.0927       & -0.0915        & -0.0936         \\
				& time   & 13.31                          & 2.388   & 6.258                         & 0.9365        & 0.3420         & 0.4218          \\
				\multirow{-3}{*}{(0.1, 100)}  & prob   & 0.8600                         & 0.7633  &  0.9233 & 0.8533        & 0.9067         & 0.8667          \\ \hline
				& fval   & -0.0874 & -0.0713 & -0.0826                       & -0.0864       & -0.0829        & -0.0870         \\
				& time   & 148.6                          & 67.95   & 81.95                         & 2.7323        & 1.2546         & 1.2818          \\
				\multirow{-3}{*}{(0.1, 400)}  & prob   & 0.9058                         & 0.8158  &  0.9266 & 0.8800        & 0.9267         & 0.9133          \\ \hline
			\end{tabular}
		}
	\end{table}
	
	In the above experiments, we compare our algorithm with several classical methods, including the mixed‐integer program (MIP, \cite{pagnoncelli2009sample}), the augmented Lagrangian decomposition method (ALDM, \cite{bai2021augmented}), the proximal difference‐of‐convex algorithm  (PDCA, \cite{wang2023proximal}), and the diffusion‐based Look‐Ahead Guidance (LAG, \cite{guogradient}) method. We set $\rho=0.05,0.1$ and $n=100,400$, reporting the final‐iteration objective function value (fval), total runtime (time), and the empirical probability of the chance constraint computed over randomly sampled daily returns (prob).
	
	The results show that MIP achieves the lowest objective values but incurs the highest computational cost, as it fully exploits the data by formulating CCP as mixed integer program. LAG attains competitive objectives but requires additional back‐propagation steps for guidance. In contrast, GGDOpt well balances solution quality and efficiency, significantly reducing runtime while maintaining comparable objective values and constraint satisfaction.

	\subsubsection{Robust waveform design}
	
	\begin{figure}[!b]
		\centering
		\begin{subfigure}[b]{0.95\textwidth}
			\centering
			\includegraphics[width=\linewidth]{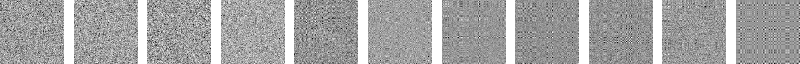}
		\end{subfigure}
		
		\vspace{1mm}
		\begin{subfigure}[b]{0.95\textwidth}
			\centering
			\includegraphics[width=\linewidth]{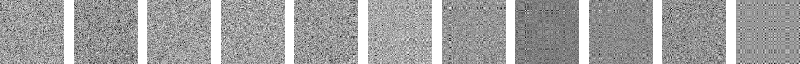}
		\end{subfigure}
		
		\vspace{1mm}
		\begin{subfigure}[b]{0.95\textwidth}
			\centering
			\includegraphics[width=\linewidth]{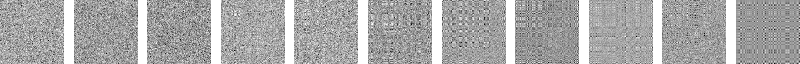}
		\end{subfigure}
		
		\vspace{1mm}
		\begin{subfigure}[b]{0.95\textwidth}
			\centering
			\includegraphics[width=\linewidth]{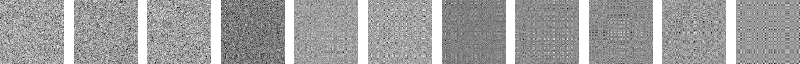}
		\end{subfigure}
		
		\vspace{1mm}
		\begin{subfigure}[b]{0.95\textwidth}
			\centering
			\includegraphics[width=\linewidth]{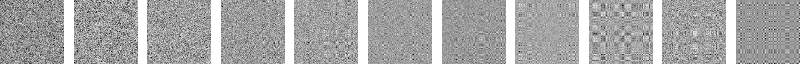}
		\end{subfigure}
		
		\vspace{1mm}
		\begin{subfigure}[b]{0.95\textwidth}
			\centering
			\includegraphics[width=\linewidth]{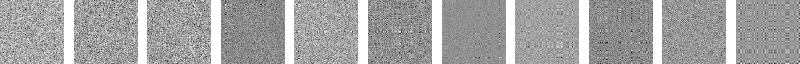}
		\end{subfigure}
		
		\vspace{1mm}
		\begin{subfigure}[b]{0.95\textwidth}
			\centering
			\includegraphics[width=\linewidth]{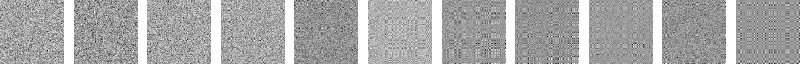}
		\end{subfigure}
		
		\vspace{1mm}
		\begin{subfigure}[b]{0.95\textwidth}
			\centering
			\includegraphics[width=\linewidth]{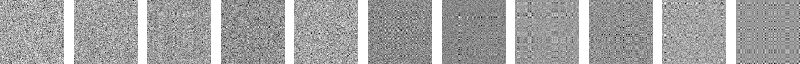}
		\end{subfigure}
		
		\vspace{1mm}
		\begin{subfigure}[b]{0.95\textwidth}
			\centering
			\includegraphics[width=\linewidth]{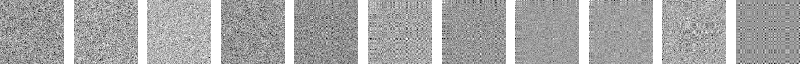}
		\end{subfigure}
		
		\vspace{1mm}
		\begin{subfigure}[b]{0.95\textwidth}
			\centering
			\includegraphics[width=\linewidth]{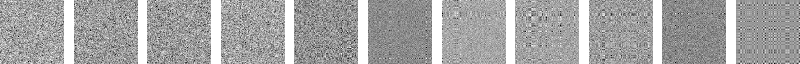}
		\end{subfigure}
		\caption{Generated 10 sampling process of GGDOpt with U-Net-2D (from left to right: $t=100, 90, 80, 70, 60, 50, 40, 30, 20, 10, 0$). }
		\label{afig:rank2d}
	\end{figure}
	
	\begin{table}[!t]
		\caption{Optimization Methods Comparison}
		\label{atab:results2}
		\renewcommand{\arraystretch}{1.6}
		\centering
			\begin{tabular}{clcccc}
				\hline
				$N_t=64,K=8$                                                                                          & Metric    & $\rho=0.05$ & $\rho=0.10$ & $\rho=0.15$ & $\rho=0.20$ \\ \hline
				\multirow{3}{*}{\begin{tabular}[c]{@{}c@{}}Empirical   Mean\end{tabular}}            & WorstProb & 0.4865   & 0.4865   & 0.4865   & 0.4865   \\
				& FuncValue & 0.0675   & 0.0675   & 0.0675   & 0.0675   \\
				& Runtime   & 4.1406   & 4.1406   & 4.1406   & 4.1406   \\ \hline
				\multirow{3}{*}{\begin{tabular}[c]{@{}c@{}}Sphere   Bounding\\\cite{ben2000robust}\end{tabular}}           & WorstProb & 0.9999   & 0.9999   & 0.9999   & 0.9999   \\
				& FuncValue & 0.0752   & 0.0750   & 0.0749   & 0.0748   \\
				& Runtime   & 1278     & 1292     & 1167     & 1131     \\ \hline
				\multirow{3}{*}{\begin{tabular}[c]{@{}c@{}}Bernstein-type   Inequality\\\cite{wang2014outage}\end{tabular}} & WorstProb & 0.9582   & 0.9335   & 0.9122   & 0.8974   \\
				& FuncValue & 0.0689   & 0.0687   & 0.0686   & 0.0685   \\
				& Runtime   & 737      & 703      & 762      & 688      \\ \hline
				\multirow{3}{*}{\begin{tabular}[c]{@{}c@{}}GGDOpt\\      (First-order)\end{tabular}}            & WorstProb & 0.9521   & 0.9097   & 0.8685   & 0.8107   \\
				& FuncValue & 0.0692   & 0.0690   & 0.0686   & 0.0685   \\
				& Runtime   & \textbf{0.6071}   & \textbf{0.5894}   & \textbf{0.5941}   & \textbf{0.5374}   \\ \hline
				\multirow{3}{*}{\begin{tabular}[c]{@{}c@{}}GGDOpt\\      (Second-order)\end{tabular}}           & WorstProb & 0.9515   & 0.9007   & 0.8573   & 0.8111   \\
				& FuncValue & \textbf{0.0688}   & \textbf{0.0685}   & \textbf{0.0684}   & \textbf{0.0684}   \\
				& Runtime   & \textbf{0.6273}   & \textbf{0.6152}   & \textbf{0.6730}   & \textbf{0.5901}   \\ \hline
			\end{tabular}
	\end{table}
	
	Consider a multiuser multiple-input single-output (MISO) downlink scenario, where a multi-antenna base station transmits independent messages to $K$ single-antenna users over a quasi-static channel. The system model adopted is standard and is briefly described as follows.
	
	Let $N_t$ denote the number of antennae at the base station and $K$ the number of users. The received signal of user $i$, $i=1,\ldots,K$, is modeled as
	\begin{equation}
		y_i(t) = \bm{h}_i^H\bm{x}(t)+\nu_i(t), 
	\end{equation}
	where $\bm{h}_i\in\mathbb{R}^{N_t}$ is the channel of user $i$; $\bm{x}(t)\in\mathbb{R}^{N_t}$ is the transmit signal from the base station; $\nu_i(t)$ is noise with distribution $\mathcal{N}(0,\sigma_i^2)$. 
	
	We assume a general vector-Gaussian linear precoding strategy, where the transmit signal is expressed as
	\begin{equation}
		\bm{x}(t)=\sum_{i=1}^{K}\bm{x}_i(t), 
	\end{equation}
	with $\bm{x}_i(t)\in\mathbb{R}^{N_t}$ representing the information-bearing signal intended for user $i$. Each $\bm{x}_i(t)$ is independently Gaussian encoded with covariance matrix $\bm{S}_i \succeq \bm{0}$, i.e., $\bm{x}_i(t)\sim\mathcal{N}(0,\bm{S}_i)$. At the receiver side, each user decodes only its own intended signal while treating the signals of other users as interference.
	
	Under this system model, the achievable rate for user $i$ can be formulated as	
	\begin{equation}
		\mathrm{R}_i=\log_2\left(1+\frac{\boldsymbol{h}_i^H\boldsymbol{S}_i\boldsymbol{h}_i}{\sum_{k\neq i}\boldsymbol{h}_i^H\boldsymbol{S}_k\boldsymbol{h}_i+\sigma_i^2}\right),i=1,\ldots,K.
	\end{equation}
	
	To formulate the rate-constrained optimization problem under imperfect channel state information (CSI), it is essential to first characterize the CSI error model. In the presence of imperfect CSI, the actual channel vector of each user can be represented as
	\begin{equation}
		\boldsymbol{h}_i=\bar{\bm{h}}_i+\boldsymbol{e}_i,i=1,\ldots,K,
	\end{equation}
	where $\bar{\bm{h}}_i\in\mathbb{R}^{N_t}$ is the presumed channel at the base station and $\bm{e}_i\in\mathbb{R}^{N_t}$ is the channel error vector. We adopt the commonly used Gaussian channel error model. Specifically, each channel error vector is assumed to have a Gaussian distribution, i.e.,
	\begin{equation}
		\bm{e}_i\sim\mathcal{N}(\bm{0},\bm{C}_i),
	\end{equation}
	for some known error covariance matrix $\bm{C}_i$. Now, consider the following probabilistically robust design formulation(\cite{wang2014outage}): 
	\begin{equation}
		\label{aeq:robust}
		\begin{aligned}
			\operatorname*{min}_{\boldsymbol{S}_{1},\ldots,\boldsymbol{S}_{K}\in\mathbb{R}^{N_{t}\times N_{t}}}\ &\sum_{i=1}^K\mathrm{Tr}(\boldsymbol{S}_i)\\
			\mathrm{s.t.}\quad \quad \ \ \ &\mathrm{Prob}_{\boldsymbol{h}_{i}\sim\mathcal{N}(\bar{\boldsymbol{h}}_{i},\boldsymbol{C}_{i})}\{\mathrm{R}_{i}\geq r_{i}\}\geq 1-\rho_{i},i=1,2,\ldots,K,\\
			&\boldsymbol{S}_{1},\ldots,\boldsymbol{S}_{K}\succeq\boldsymbol{0},i=1,2,\ldots,K.
		\end{aligned}
	\end{equation}

	To solve the aforementioned problem using GGDOpt, a naive approach is to treat each covariance matrix as a two-dimensional array and employ a 2D U-Net architecture directly. However, this approach is computationally inefficient, as it requires learning $N_t\times N_t\times K$ variables. To reduce the dimensionality of the optimization variables, we apply Cholesky factorization by expressing each covariance matrix as
	\begin{equation}
		\bm{S}_i=\bm{L}_i\bm{L}_i^T.
	\end{equation}
	This transformation reduces the number of variables per matrix from $N_t^2$ to $N_t(N_t+1)/2$, while also ensuring that $\bm{S}_i$ remains symmetric and positive semidefinite.

	Subsequently, we illustrate representative sampling trajectories of GGDOpt after training (see Figure \ref{afig:rank2d}) and observe that the generated solutions consistently approximate rank-one matrices.
	
	Remarkably, the generated samples consistently preserve the rank-one property, with the dominant eigenvalue accounting for over 99\% of the total eigenvalue. This observation suggests that solutions to the robust waveform design problem (\ref{aeq:robust}) inherently lie on a rank-one manifold with very high probability (\cite{wang2014outage}), a structure that GGDOpt can effectively captures. Consequently, rank-one decomposition can be reliably applied after generation, allowing the use of U-Net-1D as a score estimator, which substantially reduces computational costs during both training and sampling process.

	Next, we present comparative results for the case $N_t=64, K=8$ in Table \ref{atab:results2}. We compare three approximation methods with our proposed GGDOpt. The Empirical Mean approach directly utilizes the sample mean of the channel realizations $\bm{h}_i^{(\ell)}$ and solves the resulting deterministic problem. The Sphere Bounding method (\cite{ben2000robust}) and the Bernstein-type Inequality approach (\cite{wang2014outage}) construct inner convex approximations of the original nonconvex feasible region. For all users, we set $\rho_i=\rho$ for $i=1,\ldots,K$, and evaluate the worst-case outage probability using the true underlying distribution. A solution is deemed feasible if the worst-case probability exceeds $1-\rho$.
	
	The results demonstrate that across different values of $\rho$, GGDOpt consistently finds feasible solutions with lower objective values than existing convex restriction methods. Moreover, GGDOpt achieves significantly higher computational efficiency.
	
\begin{figure}[htbp]
		\centering
		\begin{subfigure}[b]{0.95\textwidth}
			\centering
			\includegraphics[width=\linewidth]{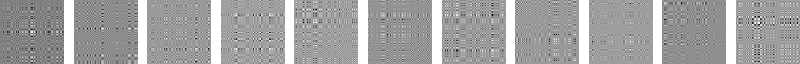}
		\end{subfigure}
		
		\vspace{1mm}
		\begin{subfigure}[b]{0.95\textwidth}
			\centering
			\includegraphics[width=\linewidth]{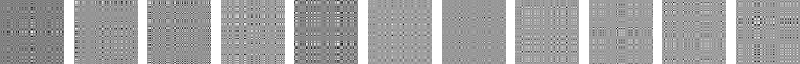}
		\end{subfigure}
		
		\vspace{1mm}
		\begin{subfigure}[b]{0.95\textwidth}
			\centering
			\includegraphics[width=\linewidth]{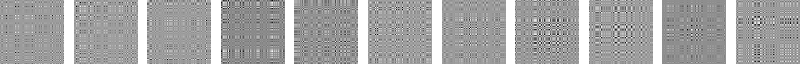}
		\end{subfigure}
		
		\vspace{1mm}
		\begin{subfigure}[b]{0.95\textwidth}
			\centering
			\includegraphics[width=\linewidth]{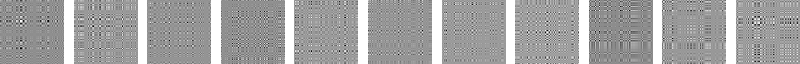}
		\end{subfigure}
		
		\vspace{1mm}
		\begin{subfigure}[b]{0.95\textwidth}
			\centering
			\includegraphics[width=\linewidth]{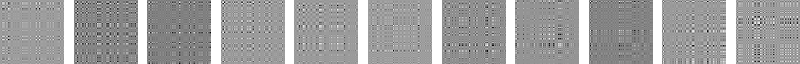}
		\end{subfigure}
		
		\vspace{1mm}
		\begin{subfigure}[b]{0.95\textwidth}
			\centering
			\includegraphics[width=\linewidth]{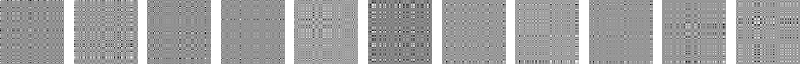}
		\end{subfigure}
		
		\vspace{1mm}
		\begin{subfigure}[b]{0.95\textwidth}
			\centering
			\includegraphics[width=\linewidth]{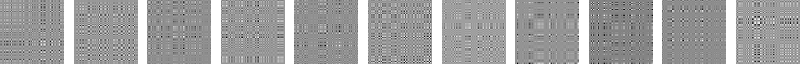}
		\end{subfigure}
		
		\vspace{1mm}
		\begin{subfigure}[b]{0.95\textwidth}
			\centering
			\includegraphics[width=\linewidth]{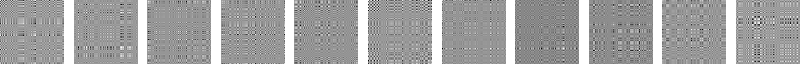}
		\end{subfigure}
		
		\vspace{1mm}
		\begin{subfigure}[b]{0.95\textwidth}
			\centering
			\includegraphics[width=\linewidth]{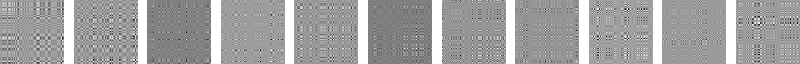}
		\end{subfigure}
		
		\vspace{1mm}
		\begin{subfigure}[b]{0.95\textwidth}
			\centering
			\includegraphics[width=\linewidth]{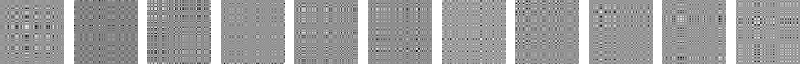}
		\end{subfigure}
		\caption{Generated 10 sampling process of GGDOpt with U-Net-1D (from left to right: $t=100, 90, 80, 70, 60, 50, 40, 30, 20, 10, 0$). }
		\label{afig:rank1d}
	\end{figure}

	By employing U-Net-1D, the sampling process is constrained to produce rank-one solutions. Representative sampling trajectories are illustrated in Figure \ref{afig:rank1d}.

	\section{Restricted problem}
	
	\subsection{Connection with CCP}
	
	In this subsection, we establish the connection between the solution of the restricted problem
	\begin{equation}
		\label{aeq:res}
		\begin{aligned}
			\min_{\bm{x}}\quad&f(\bm{x})\\
			\mathrm{s.t.}\quad& \bm{g}(\bm{x},\bar{\bm{h}})\geq\bm{z},
		\end{aligned}
	\end{equation}
	and that of the CCP
	\begin{equation}
		\label{aeq1}
		\begin{aligned}
			\min_{\bm{x}}\quad&f(\bm{x})\\
			\mathrm{s.t.}\quad&\bm{x}\in\mathcal{X}_{\rho}.
		\end{aligned}
	\end{equation}
	
	The rationale behind using the restricted problem (RP) to generate high-quality solutions is straightforward. First, solving the restricted problem (\ref{aeq:res}) is computationally more tractable than directly tackling the original CCP (\ref{aeq1}). Second, the distribution 
	$P$ of the random variable $\bm{h}$ tends to concentrate around its mean $\bm{\mu}_P$. Consequently, improving the value of $\bm{g}(\bm{x},\bm{\mu}_P)$ generally leads to an increase in the probability $\text{Prob}_{\bm{h}}\{\bm{g}(\bm{x},\bm{h})\geq\bm{0}\}$. Third, the feasible region of the RP can be viewed as an approximation of the feasible set $\mathcal{X}_{\rho}$ associated with the CCP. For a given risk level $\rho$, solving the RP yields an approximate local optimum of CCP (\ref{aeq1}). Moreover, if the global solution to (\ref{aeq1}) satisfies certain regularity conditions, this approximate local minimizer coincides with the global minimizer.
	
	In general, the quantity $\text{Prob}_{\boldsymbol{h}}\{\boldsymbol{g}(\boldsymbol{x}(\boldsymbol{z}_i),\boldsymbol{h})\geq\boldsymbol{0}\}$ is hard to compute since it requires a multidimensional integration over the distribution of $\boldsymbol{h}$. Inspired by the sample average approximation (SAA), we estimate this by an empirical average based on $L$ i.i.d. realizations of $\boldsymbol{h}$:
	\begin{equation}
		\begin{aligned}
			\text{Prob}_{\boldsymbol{h}}\{\boldsymbol{g}(\boldsymbol{x}(\boldsymbol{z}_i),\boldsymbol{h})\geq\boldsymbol{0}\}\approx 			\frac{1}{L}\sum_{l=1}^{L}\ell^{0|1}(\boldsymbol{g}(\boldsymbol{x}(\boldsymbol{z}_i),\boldsymbol{h}^{(\ell)})) {:=1-\rho^{(i)}},
		\end{aligned}
	\end{equation}
	where $\ell^{0/1}$ is the element-wise indicator function that returns 1 if all components of the argument vector are positive, and 0 otherwise. The reason why we choose this to approximate $\rho^{(i)}$ can be analyzed from the following two situations:
	
	On the one hand, if the sample size $L$ is large enough, then the empirical distribution can be regarded as a good approximation of the underlying distribution, i.e., $p(\boldsymbol{h})\approx\frac{1}{L}\sum_{l=1}^{L}\delta(\boldsymbol{h}-\boldsymbol{h}^{(\ell)})$. In this case, it is natural to replace the real value that computationally intractable with the empirical value $\rho^{(i)}=1-\frac{1}{L}\sum_{l=1}^{L}\ell^{0|1}(\boldsymbol{g}(\boldsymbol{x}(\boldsymbol{z}_i),\boldsymbol{h}^{(\ell)}))$.
	
	On the other hand, if the sample size $L$ is small, using the empirical value to estimate the real $\rho$ will cause serious distortion. In this case, a larger restriction $\boldsymbol{z}_i$ is preferred, as it will lead to $\boldsymbol{x}(\boldsymbol{z}_i)$ with greater probability of satisfying the chance constraint and better robustness to the distribution uncertainty. At this time, the empirical value $\rho^{(i)}$ is not used to approximate the real confidence, but to characterize the properties of "good" $\boldsymbol{x}(\boldsymbol{z}_i)$.
	
	To compute $\rho^{(i)}$, we proceed as follows:
	\begin{itemize}
		\item For each sampled restriction vector $\boldsymbol{z}_i \geq \boldsymbol{0}$, we solve the corresponding restricted problem, which yields a candidate solution $\boldsymbol{x}(\boldsymbol{z}_i)$.
		\item We then draw $L$ independent realizations $\boldsymbol{h}^{(\ell)}$ from the underlying distribution and evaluate the fraction of those samples for which $\boldsymbol{g}(\boldsymbol{x}(\boldsymbol{z}_i), \boldsymbol{h}^{(\ell)}) \geq 0$ holds.
	\end{itemize}
	This empirical feasible set constructed in this way provides a conservative inner approximation of the true feasible region, ensuring that the solutions obtained from the restricted problem satisfy the original chance constraint with high confidence.
	
	Next, we provide a detailed characterization of the probability $\bm{g}(\bm{x},\bm{h})\geq\bm{0}$ evaluated at the solution $\bm{x}(\bar{\bm{h}},\bm{z})$ to the restricted problem. For brevity, we denote the norm $\|\cdot\| = \|\cdot\|_{\infty}$ throughout the subsequent analysis.
	
	\textbf{Assumption 2.} Assume that 
	\begin{itemize}[leftmargin=10pt]
		\item (Lipschitz continuity) $\bm{g}(\bm{x},\cdot)$ is Lipschitz for a given $\bm{x}$, i.e.,
		\begin{equation}
			\|\bm{g}(\bm{x},\bm{h})-\bm{g}(\bm{x},\bm{h}^\prime)\|\leq L_{\bm{x}}\|\bm{h}-\bm{h}^\prime\|,\quad \forall~ \bm{h},~ \bm{h}^\prime,
		\end{equation}
		where $L_{\bm{x}}$ is the Lipschitz constant depending on $\bm{x}$.
		\item (Finite variance) The variance of the random vector $\bm{h}$ with probability $P$ is finite, i.e.,
		\begin{equation}
			\text{Var}_{P}(\bm{h})< \infty.
		\end{equation}
	\end{itemize}
	
	\textbf{Theorem 4.} Under Assumption 2, suppose that $\{\bm{h}^{(\ell)}\}_{\ell=1}^{L}$ are samples drawn from the distribution $P$ of random vector $\bm{h}$. Let $\bar{\bm{h}} = \frac{1}{L}\sum_{\ell=1}^{L} \bm{h}^{(\ell)}$ and let $z_{min}$ be the smallest element of $\bm{z}$. Suppose that $\bm{x}(\bar{\bm{h}},\bm{z})$ is the solution to the problem (\ref{aeq:res}), then we have
	\begin{equation}
		\label{eq:prob}
		\text{Prob}_{\bm{h}}\big\{\bm{g}(\bm{x}(\bar{\bm{h}},\bm{z}),\bm{h})\geq \bm{0}\big\} \geq\underbrace{ 1-\frac{\text{Var}_{P}(\bm{h})}{(z_{min}/L_{\bm{x}(\bar{\bm{h}},\bm{z})}-\|\bar{\bm{h}}-\mathbb{E}_{P}\left[\bm{h}\right]\|)^2}}_{1-\rho}.
	\end{equation}
	
	\textit{Proof. } 
	
	To characterize $\text{Prob}_{\bm{h}}\{\bm{g}(\bm{x}(\bar{\bm{h}},\bm{z}),\bm{h})\geq \bm{0}\}$, we need to consider two sources of error. The first arises from the large variance of the distribution $P$, while the second stems from the approximation of the mean of $P$ using a finite number of realizations, i.e.,
	\begin{equation}
		\begin{aligned}
			& \text{Prob}_{\bm{h}}\big\{\bm{g}(\bm{x}(\bar{\bm{h}},\bm{z}),\bm{h})\geq \bm{0}\big\} \\
			=\ & \text{Prob}_{\bm{h}}\big\{\bm{g}(\bm{x}(\bar{\bm{h}},\bm{z}),\bm{h})-\bm{g}(\bm{x}(\bar{\bm{h}},\bm{z}),\mathbb{E}_{P}\left[\bm{h}\right]) + \bm{g}(\bm{x}(\bar{\bm{h}},\bm{z}),\mathbb{E}_{P}\left[\bm{h}\right])-\bm{g}(\bm{x}(\bar{\bm{h}},\bm{z}),\bar{\bm{h}}) \\
			&  + \bm{g}(\bm{x}(\bar{\bm{h}},\bm{z}),\bar{\bm{h}}) \geq \bm{0}\big\}.
		\end{aligned}
	\end{equation}
	Since $\bm{g}(\bm{x}(\bar{\bm{h}},\bm{z}),\bar{\bm{h}})$ is the solution to the restricted problem (\ref{aeq:res}), we have $\bm{g}(\bm{x}(\bar{\bm{h}},\bm{z}),\bar{\bm{h}})\geq z_{min}\bm{1}$. Therefore, we have 
	\begin{equation}
		\begin{aligned}
			& \text{Prob}_{\bm{h}}\big\{\bm{g}(\bm{x}(\bar{\bm{h}},\bm{z}),\bm{h})\geq \bm{0}\big\} \\
			\geq\ & \text{Prob}_{\bm{h}}\big\{\|\bm{g}(\bm{x}(\bar{\bm{h}},\bm{z}),\bm{h})-\bm{g}(\bm{x}(\bar{\bm{h}},\bm{z}),\mathbb{E}_{P}\left[\bm{h}\right])\| + \|\bm{g}(\bm{x}(\bar{\bm{h}},\bm{z}),\mathbb{E}_{P}\left[\bm{h}\right])-\bm{g}(\bm{x}(\bar{\bm{h}},\bm{z}),\bar{\bm{h}})\|  \\
			& - z_{min} \leq 0\big\}.
		\end{aligned}
	\end{equation}
	According to Assumption 2, we have that 
	\begin{equation}
		\|\bm{g}(\bm{x}(\bar{\bm{h}},\bm{z}),\bm{h})-\bm{g}(\bm{x}(\bar{\bm{h}},\bm{z}),\mathbb{E}_{P}\left[\bm{h}\right])\|\leq L_{\bm{x}(\bar{\bm{h}},\bm{z})} \|\bm{h}-\mathbb{E}_{P}\left[\bm{h}\right]\|,
	\end{equation}
	and 
	\begin{equation}
		\|\bm{g}(\bm{x}(\bar{\bm{h}},\bm{z}),\mathbb{E}_{P}\left[\bm{h}\right])-\bm{g}(\bm{x}(\bar{\bm{h}},\bm{z}),\bar{\bm{h}})\|\leq L_{\bm{x}(\bar{\bm{h}},\bm{z})} \|\bar{\bm{h}}-\mathbb{E}_{P}\left[\bm{h}\right]\|.
	\end{equation}
	
	Therefore, the probability $\text{Prob}_{\bm{h}}\big\{\bm{g}(\bm{x}(\bar{\bm{h}},\bm{z}),\bm{h})\geq \bm{0}\big\}$ can be further expressed as
	
	\begin{equation}
		\begin{aligned}
			&\text{Prob}_{\bm{h}}\big\{\bm{g}(\bm{x}(\bar{\bm{h}},\bm{z}),\bm{h})\geq \bm{0}\big\} \\
			\geq \ & \text{Prob}_{\bm{h}}\big\{\|\bm{g}(\bm{x}(\bar{\bm{h}},\bm{z}),\bm{h})-\bm{g}(\bm{x}(\bar{\bm{h}},\bm{z}),\mathbb{E}_{P}\left[\bm{h}\right])\|\leq z_{min}-\|\bm{g}(\bm{x}(\bar{\bm{h}},\bm{z}),\mathbb{E}_{P}\left[\bm{h}\right])-\bm{g}(\bm{x}(\bar{\bm{h}},\bm{z}),\bar{\bm{h}})\|\big\} \\
			\geq\  & \text{Prob}_{\bm{h}}\big\{L_{\bm{x}(\bar{\bm{h}},\bm{z})} \|\bm{h}-\mathbb{E}_{P}\left[\bm{h}\right]\|\leq z_{min}-L_{\bm{x}(\bar{\bm{h}},\bm{z})} \|\bar{\bm{h}}-\mathbb{E}_{P}\left[\bm{h}\right]\| \big\} \\
			=\  & \text{Prob}_{\bm{h}}\big\{ \|\bm{h}-\mathbb{E}_{P}\left[\bm{h}\right]\|\leq z_{min}/L_{\bm{x}(\bar{\bm{h}},\bm{z})}- \|\bar{\bm{h}}-\mathbb{E}_{P}\left[\bm{h}\right]\| \big\}.
		\end{aligned}
	\end{equation}
	
	By Chebyshev’s inequality (\cite{chebyshev1867valeurs}), we obtain that 
	\begin{equation}
		\begin{aligned}
			&\text{Prob}_{\bm{h}}\big\{ \|\bm{h}-\mathbb{E}_{P}\left[\bm{h}\right]\|\leq z_{min}/L_{\bm{x}(\bar{\bm{h}},\bm{z})}- \|\bar{\bm{h}}-\mathbb{E}_{P}\left[\bm{h}\right]\| \big\} \\
			\geq\ &  1-\frac{\text{Var}_{P}(\bm{h})}{(z_{min}/L_{\bm{x}(\bar{\bm{h}},\bm{z})}-\|\bar{\bm{h}}-\mathbb{E}_{P}\left[\bm{h}\right]\|)^2}.
		\end{aligned}
	\end{equation}
	
	Hence, we have
	\begin{equation}
		\text{Prob}_{\bm{h}}\big\{\bm{g}(\bm{x}(\bar{\bm{h}},\bm{z}),\bm{h})\geq \bm{0}\big\} \geq 1-\frac{\text{Var}_{P}(\bm{h})}{(z_{min}/L_{\bm{x}(\bar{\bm{h}},\bm{z})}-\|\bar{\bm{h}}-\mathbb{E}_{P}\left[\bm{h}\right]\|)^2}.
	\end{equation}
	
	$\hfill\square$
	
	Theorem 4 demonstrates that as $z_{min}$ increases, the lower bound on the probability that the chance constraint is satisfied at the point $\bm{x}(\bar{\bm{h}},\bm{z})$ also increases. This implies that $\bm{x}(\bar{\bm{h}},\bm{z})$ is more likely to be a feasible solution to the CCP (\ref{aeq1}), while potentially achieving a lower objective value. In the following theorem, we further establish that, under certain regularity conditions, the global minimizer of the CCP (\ref{aeq1}) is contained within the set of solutions to the restricted problem (\ref{aeq:res}).
	
	\textbf{Assumption 3.} Assume that 
	\begin{itemize}[leftmargin=10pt]
		\item (Bounded bias) For any given $\rho$, denote $\bm{x}^* = \mathop{\arg\min}_{\bm{x}\in\mathcal{X}_{\rho}} f(\bm{x})$, then 
		\begin{equation}
			\|\bar{\bm{h}}-\mathbb{E}_{P}\left[\bm{h}\right]\|\leq \frac{\bm{g}(\bm{x}^*,\bar{\bm{h}})}{L_{\bm{x}(\bar{\bm{h}},\bm{z})}}-\sqrt{\frac{\text{Var}_{P}(\bm{h})}{\rho}}.
		\end{equation}
		\item (Reliable data set) For the generated data set $\mathcal{D}=\{(\bm{x}^{(i)},\rho^{(i)})\}_{i=1}^N$, $\rho^{(i)}$ is a lower bound of real probability $\text{Prob}_{\bm{h}}\{\bm{g}(\bm{x}^{(i)},\bm{h})\geq \bm{0}\}$.
	\end{itemize}

	Note that Assumption 3 can be satisfied with a sufficiently large number of realizations of $\bm{h}$ and the corresponding restriction estimator. For instance, we can choose
	\begin{equation}
		\rho^{(i)}\leq \frac{\text{Var}_{P}(\bm{h})}{(z_{min}/L_{\bm{x}(\bar{\bm{h}},\bm{z})}-\|\bar{\bm{h}}-\mathbb{E}_{P}\left[\bm{h}\right]\|)^2}.
	\end{equation}
	
	\textbf{Theorem 5.} Under Assumption 2 and Assumption 3, for any given $\rho$ and $\bar{\bm{h}}$, suppose that 
	\begin{equation}
		\mathcal{D}_{\rho}=\{\bm{x}^{(i)}\mid (\bm{x}^{(i)},\rho^{(i)})\in\mathcal{D}, \rho^{(i)}\leq \rho\},
	\end{equation}
	then we have
	\begin{equation}
		\bm{x}^* \in \mathcal{D}_{\rho}\subset\mathcal{X}_{\rho}.
	\end{equation}
	
	\textit{Proof. } 
	
	We choose $z_{min}$ as the smallest element of $\bm{g}(\bm{x}^*,\bar{\bm{h}})$, then for any $\bm{x}$ that satisfies $\bm{g}(\bm{x},\bar{\bm{h}})\geq \bm{z}$, the following inequality holds:
	\begin{equation}
		\begin{aligned}
			\text{Prob}_{\bm{h}}\{ \bm{g}(\bm{x},\bm{h})\geq \bm{0} \}\geq 1-\frac{\text{Var}_{P}(\bm{h})}{(z_{min}/L_{\bm{x}}-\|\bar{\bm{h}}-\mathbb{E}_{P}\left[\bm{h}\right]\|)^2}\geq 1-\rho.
		\end{aligned}
	\end{equation}
	This implies that 
	\begin{equation}
		\{\bm{x}\mid \bm{g}(\bm{x},\bar{\bm{h}})\geq \bm{z}\} \subset \mathcal{X}_{\rho}.
	\end{equation}
	Recall the definition of $\bm{x}(\bar{\bm{h}},\bm{z})$, which is the global minimizer of $f(\bm{x})$ over the set $\{\bm{x}\mid \bm{g}(\bm{x},\bar{\bm{h}})\geq \bm{z}\}$. Additionally, it follows naturally that $\bm{g}(\bm{x}^*,\bar{\bm{h}})\geq \bm{z}$, i.e., 
	\begin{equation}
		\bm{x}^* \in \{\bm{x}\mid \bm{g}(\bm{x},\bar{\bm{h}})\geq \bm{z}\}.
	\end{equation}
	This implies that $\bm{x}^*$ is also a global minimizer of $f(\bm{x})$ over the set $\{\bm{x}\mid \bm{g}(\bm{x},\bar{\bm{h}})\geq \bm{z}\}$. Therefore, we have 
	\begin{equation}
		\bm{x}^* \in \mathcal{D}_{\rho}\subset\mathcal{X}_{\rho}.
	\end{equation}

	$\hfill\square$

	This result plays a crucial role in the GGDOpt framework, as the sampler is inherently limited to generating solutions that are no better than the quality of the training data. Theoretical guarantees established above indicate that the data generated from the restricted problem are sufficiently informative and may contain the true global minimizer of the CCP (\ref{aeq1}). This justifies the effectiveness of using such data to train our GGDOpt.

	\subsection{Special cases}
	
	The above results provide a lower bound for the probability $\text{Prob}_{\bm{h}}\big\{\bm{g}(\bm{x}(\bar{\bm{h}},\bm{z}),\bm{h})\geq \bm{0}\big\}$. In most cases, the explicit value of this probability cannot be directly computed. However, in this subsection, we present a special case corresponding to the robust waveform design problem, where the probability can be expressed explicitly.
	
	\textbf{Theorem 6.} Suppose $\bm{x}^*(\bar{\bm{h}}_i,\bm{z})$ is the solution to the following restricted problem
	\begin{equation}
		\begin{aligned}
			\min_{\bm{x}}\quad&f(\bm{x})\\
			\mathrm{s.t.}\quad& g_i(\bm{x},\bar{\bm{h}}_i)=z_i, i=1,\ldots,K,
		\end{aligned}
	\end{equation}
	where $g_i(\bm{x},\cdot)$ is a quadratic function of $\bm{h}$ with parameters $(\bm{A}_i,\bm{b}_i, d_i)$ and the parameters $\bm{h}_i\sim\mathcal{N}(\bar{\bm{h}}_i, \bm{C}_i)$. Denote 
	\begin{equation}
		\begin{aligned}
			&\bm{Q}_i = \bm{C}_i^{1/2}\bm{A}_i\bm{C}_i^{1/2}\overset{\text{svd}}{=} \bm{U}_i\bm{\Lambda}_i\bm{U}_i^T, \\
			&\bm{r}_i = \bm{C}_i^{1/2}(\bm{A}_i\bar{\bm{h}}_i+\bm{b}_i), \\
			& s_i = \frac{1}{2}\bar{\bm{h}}_i^\top\bm{A}_i \bar{\bm{h}}_i + \bm{b}_i^\top\bar{\bm{h}}_i + d_i, \\	&\bm{c}_i = \bm{U}_i^T \bm{r}_i,
		\end{aligned}
	\end{equation}
	and let
	\begin{equation}
		\begin{aligned}
			&\bm{u}_i = \bm{U}_i^T \bm{e}_i, \bm{e}_i\sim\mathcal{N}(\bm{0}, \bm{I}), \\
			& Y_i =\frac{1}{2}\bm{u}_i^\top\bm{\Lambda}_i \bm{u}_i + \bm{c}_i^\top\bm{u}_i + s_i.
		\end{aligned}
	\end{equation}
	
	Then for $\bm{h}_i\sim\mathcal{N}(\bar{\bm{h}}_i, \bm{C}_i)$, we have
	\begin{equation}
		\text{Prob}_{\bm{h}_i}\{g_i(\bm{x}^*,\bm{h}_i)\geq0\} = 1-F_{Y_i}(0),
	\end{equation}
	where $F_{Y_i}$ is the cumulative distribution function of $Y_i$.
	
	\textit{Proof. } 
	
	For quadratic $g_i(\bm{x},\cdot)$ of $\bm{h}$ with parameters $(\bm{A}_i,\bm{b}_i,d_i)$ and given that $\bm{h}_i\sim\mathcal{N}(\bar{\bm{h}}_i, \bm{C}_i)$, the probability $\text{Prob}_{\bm{h}_i}\{g_i(\bm{x},\bm{h}_i)\geq0\}$ can be transformed into the following form:
	\begin{equation}
		\begin{aligned}
			&\text{Prob}_{\bm{h}_i\sim\mathcal{N}(\bar{\bm{h}}_i, \bm{C}_i)}\big\{g_i(\bm{x},\bm{h}_i)\geq0\big\} \\
			=\ & \text{Prob}_{\bm{h}_i\sim\mathcal{N}(\bar{\bm{h}}_i, \bm{C}_i)}\big\{ \frac{1}{2}\bm{h}_i^\top\bm{A}_i \bm{h}_i + \bm{b}_i^\top\bm{h}_i + d_i \geq0\big\} \\
			=\ & \text{Prob}_{\bm{e}_i\sim\mathcal{N}(\bm{0}, \bm{I})}\big\{ \frac{1}{2}(\bar{\bm{h}}_i+\bm{C}_i^{1/2}\bm{e}_i)^\top\bm{A}_i (\bar{\bm{h}}_i+\bm{C}_i^{1/2}\bm{e}_i) + \bm{b}_i^\top(\bar{\bm{h}}_i+\bm{C}_i^{1/2}\bm{e}_i) + d_i \geq0\big\}.
		\end{aligned}
	\end{equation}
	Denote 
	\begin{equation}
		\begin{aligned}
			&\bm{Q}_i = \bm{C}_i^{1/2}\bm{A}_i\bm{C}_i^{1/2}\overset{\text{svd}}{=} \bm{U}_i\bm{\Lambda}_i\bm{U}_i^T, \\
			&\bm{r}_i = \bm{C}_i^{1/2}(\bm{A}_i\bar{\bm{h}}_i+\bm{b}_i), \\
			&s_i = \frac{1}{2}\bar{\bm{h}}_i^\top\bm{A}_i \bar{\bm{h}}_i + \bm{b}_i^\top\bar{\bm{h}}_i + d_i,
		\end{aligned}
	\end{equation}
	then we have 
	\begin{equation}
		\begin{aligned}
			&\text{Prob}_{\bm{h}_i\sim\mathcal{N}(\bar{\bm{h}}_i, \bm{C}_i)}\big\{g_i(\bm{x},\bm{h}_i)\geq0\big\}= \text{Prob}_{\bm{e}_i\sim\mathcal{N}(\bm{0}, \bm{I})}\big\{ \frac{1}{2}\bm{e}_i^\top\bm{Q}_i \bm{e}_i + \bm{r}_i^\top\bm{e}_i + s_i \geq0\big\}. \\
		\end{aligned}
	\end{equation}
	
	Denote $\bm{Q}_i\overset{\text{svd}}{=} \bm{U}_i\bm{\Lambda}_i\bm{U}_i^T$ and let
	\begin{equation}
		\begin{aligned}
			&\bm{c}_i = \bm{U}_i^T \bm{r}_i, \\
			&\bm{u}_i = \bm{U}_i^T \bm{e}_i, \bm{e}_i\sim\mathcal{N}(\bm{0}, \bm{I}), \\
			& Y_i =\frac{1}{2}\bm{u}_i^\top\bm{\Lambda}_i \bm{u}_i + \bm{c}_i^\top\bm{u}_i + s_i.
		\end{aligned}
	\end{equation}
	
	Substituting these expressions into the above probability, we obtain that 
	\begin{equation}
		\begin{aligned}
			&\text{Prob}_{\bm{h}_i\sim\mathcal{N}(\bar{\bm{h}}_i, \bm{C}_i)}\big\{g_i(\bm{x},\bm{h}_i)\geq0\big\} \\
			=\ & \text{Prob}_{\bm{u}_i\sim\mathcal{N}(\bm{0}, \bm{I})}\big\{ \frac{1}{2}\bm{u}_i^\top\bm{\Lambda}_i \bm{u}_i + \bm{c}_i^\top\bm{u}_i + s_i \geq0\big\} \\
			=\ & \text{Prob}_{\bm{u}_i\sim\mathcal{N}(\bm{0}, \bm{I})}\{Y_i\geq 0\}. \\
		\end{aligned}
	\end{equation}
	Denote $\lambda_i^{(k)}, u_i^{(k)},$ and $c_i^{(k)}$ as the $k$-th element of $\bm{\Lambda}_i, \bm{u}_i,$ and $ \bm{c}_i$, where $k=1,\ldots,n$. Note that $Y_i$ has a quadratic form of standard Gaussian $\bm{u}_i$, which can be reformulated as a standard quadratic form: 
	\begin{equation}
		\begin{aligned}
			Y_i = \sum\limits_{\lambda_i^{(k)}\neq0}\frac{\lambda_i^{(k)}}{2}\left(u_i^{(k)}+\frac{c_i^{(k)}}{\lambda_i^{(k)}}\right)^2+\sum\limits_{\lambda_i^{(k)}=0}c_i^{(k)}u_i^{(k)} +\left(s_i-\sum\limits_{\lambda_i^{(k)}\neq0}\frac{(c_i^{(k)})^2}{2\lambda_i^{(k)}}\right),
		\end{aligned}
	\end{equation}
	where $ \left(u_i^{(k)}+\frac{c_i^{(k)}}{\lambda_i^{(k)}}\right)^2 \sim \chi^2_1((\frac{c_i^{(k)}}{\lambda_i^{(k)}})^2)$ follows noncentral chi-squared distribution and $c_i^{(k)}u_i^{(k)}\sim \mathcal{N}(0,(c_i^{(k)})^2)$ follows Gaussian distribution. 
	
	Denote $F_{Y_i}$ as the cumulative distribution function of $Y_i$, then we have
	\begin{equation}
		\text{Prob}_{\bm{h}_i\sim\mathcal{N}(\bar{\bm{h}}_i, \bm{C}_i)}\{g_i(\bm{x},\bm{h}_i)\geq0\} = 1-F_{Y_i}(0).
	\end{equation}
	
	Since $\bm{x}^*(\bar{\bm{h}}_i,\bm{z})$ is the solution to the restricted problem, by substituting $s_i = z_i$, we obtain the result of Theorem 6.
	
	$\hfill\square$
	
	Theorem 6 tells us that the probability $\text{Prob}_{\bm{h}_i}\{g_i(\bm{x}^*,\bm{h}_i)\geq0\}$ can be expressed in terms of the cumulative distribution function of $Y_i$. Note that $Y_i$ consists of $n$ independent variables. The following corollary states that, for sufficiently large $n$, $Y_i$ can be approximated as a Gaussian random variable, and the probability can be computed using the standard Gaussian cumulative distribution function $\Phi$.
	
	\textbf{Corollary 2.} For sufficiently large $n$, the probability can be approximated by
	\begin{equation}
		\text{Prob}_{\bm{h}_i}\{g_i(\bm{x},\bm{h}_i)\geq0\} \approx 1-\Phi\left(\frac{-\mu_{Y_i}}{\sigma_{Y_i}}\right),
	\end{equation}
	where $\Phi$ denotes the cumulative distribution function of the standard Gaussian distribution and 
	\begin{equation}
		\begin{aligned}
			&\mu_{Y_i} = \frac{1}{2}\text{tr}(\bm{Q}_i)+z_i, \\
			&\sigma_{Y_i}^2 = \frac{1}{2}\|\bm{Q}_i\|_F^2+\|\bm{r}_i\|^2.
		\end{aligned}
	\end{equation}
	The approximation error can be bounded by
	\begin{equation}
		\left|F_{Y_i}(0)-\Phi\left(\frac{-\mu_{Y_i}}{\sigma_{Y_i}}\right)\right| =O(n^{-1/2}).
	\end{equation}
	
	\textit{Proof. } 
	
	For sufficiently large $n$, the distribution of $Y_i$ can be approximated by Gaussian distribution $\mathcal{N}(\mu_{Y_i},\sigma_{Y_i}^2)$ with central limit theorem, where
	\begin{equation}
		\begin{aligned}
			& \mu_{Y_i} = \frac{1}{2}\text{tr}(\bm{Q}_i) + z_i, \\
			& \sigma_{Y_i}^2 = \frac{1}{2} \|\bm{Q}_i\|_F^2 + \|\bm{r}_i\|^2,
		\end{aligned}
	\end{equation}
	then the probability can be approximated by 
	\begin{equation}
		\text{Prob}_{\bm{h}_i\sim\mathcal{N}(\bar{\bm{h}}_i, \bm{C}_i)}\{g_i(\bm{x},\bm{h}_i)\geq0\} \approx 1-\Phi(\frac{-\mu_{Y_i}}{\sigma_{Y_i}}).
	\end{equation}
	The approximation error can be bounded by \cite{klartag2012variations}
	\begin{equation}
		|F_{Y_i}(0)-\Phi(\frac{-\mu_{Y_i}}{\sigma_{Y_i}})| =O(n^{-1/2}).
	\end{equation}

	$\hfill\square$

	\section{Technical appendices}

	\subsection{Proof of Theorem 1}
	
	\textbf{Theorem 1.} For any given $\beta>0$, there exists $\hat{\bm{x}}_0(\bm{x}_t)$ such that the score function of the diffused product distribution can be formulated as
	\begin{equation}
		\nabla_{\bm{x}_t}\log \tilde{p}_t(\bm{x}_t|\rho) = \nabla_{\bm{x}_t}\log p_t(\bm{x}_t|\rho) \underbrace{- \beta  \nabla_{\bm{x}_t}f\big(\hat{\bm{x}}_0(\bm{x}_t)\big)}_{\text{gradient guidance $\bm{G}_t$}},
	\end{equation}
	where $\nabla_{\bm{x}_t}\log p_t(\bm{x}_t|\rho)$ is the score function of the diffused data distribution and $\hat{\bm{x}}_0(\bm{x}_t)$ satisfies
	\begin{equation}
		f(\hat{\bm{x}}_0(\bm{x}_t)) = -\frac{1}{\beta}\log\Big(\int_{\bm{x}_0} p_{t0}(\bm{x}_0|\bm{x}_t,\rho)B_{\beta}(\bm{x}_0)d\bm{x}_0\Big).
	\end{equation}
	
	\textit{Proof. } 
	
	Given $\tilde{p}_0(\bm{x}_0|\rho)$ and the forward process $d\bm{x}=\bm{a}(\bm{x},t)dt+b(t)d\bm{B}_t$, the diffused conditional distribution of unguided distribution $p_{0}(\bm{x}_0|\rho)$ and product distribution $\tilde{p}_0(\bm{x}_0|\rho)$ satisfies
	\begin{equation}
		\begin{aligned}
			&p_t(\bm{x}_t|\rho) = \int_{\bm{x}_0} p_{0t}(\bm{x}_t|\bm{x}_0)p_{0}(\bm{x}_0|\rho)d\bm{x}_0, \\
			&\tilde{p}_t(\bm{x}_t|\rho) = \int_{\bm{x}_0} p_{0t}(\bm{x}_t|\bm{x}_0)\tilde{p}_0(\bm{x}_0|\rho)d\bm{x}_0 \propto \int_{\bm{x}_0} p_{0t}(\bm{x}_t|\bm{x}_0)p_{0}(\bm{x}_0|\rho)B_{\beta}(\bm{x}_0)d\bm{x}_0.
		\end{aligned}
	\end{equation}
	
	Consider the difference between the score function of unguided $p_t(\bm{x}_t|\rho)$ and guided $\tilde{p}_t(\bm{x}_t|\rho)$, we have that 
	\begin{equation}
		\begin{aligned}
			&\nabla_{\bm{x}_t}\log \tilde{p}_t(\bm{x}_t|\rho) - \nabla_{\bm{x}_t}\log p_t(\bm{x}_t|\rho) \\
			=\ & \nabla_{\bm{x}_t}\log \int_{\bm{x}_0} p_{0t}(\bm{x}_t|\bm{x}_0)p_{0}(\bm{x}_0|\rho)B_{\beta}(\bm{x}_0)d\bm{x}_0 - \nabla_{\bm{x}_t}\log \int_{\bm{x}_0} p_{0t}(\bm{x}_t|\bm{x}_0)p_{0}(\bm{x}_0|\rho)\\
			=\ & \nabla_{\bm{x}_t} \log \frac{\int_{\bm{x}_0} p_{0t}(\bm{x}_t|\bm{x}_0)p_{0}(\bm{x}_0|\rho)B_{\beta}(\bm{x}_0)d\bm{x}_0}{\int_{\bm{x}_0} p_{0t}(\bm{x}_t|\bm{x}_0)p_{0}(\bm{x}_0|\rho)}.
		\end{aligned}
	\end{equation}
	
	Notice that the inner fractional part can be expressed by 
	\begin{equation}
		\frac{p_{0t}(\bm{x}_t|\bm{x}_0)p_{0}(\bm{x}_0|\rho)}{\int_{\bm{x}_0} p_{0t}(\bm{x}_t|\bm{x}_0)p_{0}(\bm{x}_0|\rho)} = p(\bm{x}_0|\bm{x}_t,\rho),
	\end{equation}
	then we have
	\begin{equation}
		\begin{aligned}
			&\nabla_{\bm{x}_t}\log \tilde{p}_t(\bm{x}_t|\rho) - \nabla_{\bm{x}_t}\log p_t(\bm{x}_t|\rho)= \nabla_{\bm{x}_t}\log \int_{\bm{x}_0} p(\bm{x}_0|\bm{x}_t,\rho)B_{\beta}(\bm{x}_0)d\bm{x}_0.
		\end{aligned}
	\end{equation}
	
	One way to tackle the log integral is to use the mean value theorem. There exists $\hat{\bm{x}}_0(\bm{x}_t)$ such that
	\begin{equation}
		\int_{\bm{x}_0} p(\bm{x}_0|\bm{x}_t,\rho)B_{\beta}(\bm{x}_0)d\bm{x}_0 = B_{\beta}(\hat{\bm{x}}_0(\bm{x}_t)) \int_{\bm{x}_0} p(\bm{x}_0|\bm{x}_t,\rho)d\bm{x}_0.
	\end{equation}
	Then we have
	\begin{equation}
		\begin{aligned}
			&\nabla_{\bm{x}_t}\log \tilde{p}_t(\bm{x}_t|\rho) - \nabla_{\bm{x}_t}\log p_t(\bm{x}_t|\rho)=\nabla_{\bm{x}_t}\log B_{\beta}(\hat{\bm{x}}_0(\bm{x}_t))= -\beta \nabla_{\bm{x}_t} f(\hat{\bm{x}}_0(\bm{x}_t)), 
		\end{aligned}
	\end{equation}
	and $\hat{\bm{x}}_0(\bm{x}_t)$ satisfies
	\begin{equation}
		f(\hat{\bm{x}}_0(\bm{x}_t)) = -\frac{1}{\beta}\log\left(\frac{\int_{\bm{x}_0} p_{0t}(\bm{x}_t|\bm{x}_0)p_{0}(\bm{x}_0|\rho)B_{\beta}(\bm{x}_0)d\bm{x}_0}{\int_{\bm{x}_0} p_{0t}(\bm{x}_t|\bm{x}_0)p_{0}(\bm{x}_0|\rho)d\bm{x}_0}\right).
	\end{equation}
	
	$\hfill\square$
	
	\subsection{Proof of Corollary 1}
	
	\textbf{Corollary 1.} Assume that $p_{t0}(\bm{x}_0|\bm{x}_t,\rho) = \mathcal{N}(\bm{x}_0|\bm{\mu}_{0|t},\sigma_{0|t}^2\bm{I})$, then we have the following results.
	\begin{itemize}[leftmargin=10pt]
		\item \textbf{First-order guidance:} For $f\in\mathcal{C}^{1}(\mathbb{R}^n,\mathbb{R})$, we get
		\begin{equation}
			\bm{G}_t = - \beta \nabla_{\bm{x}_t} f(\bm{x}_t).
		\end{equation}
		\item \textbf{Second-order guidance:} For $f\in\mathcal{C}^{2}(\mathbb{R}^n,\mathbb{R})$, we get
		\begin{equation}
			\bm{G}_t = -\frac{1}{\sigma_{0|t}^2}\left[\bm{H}^{-1}\left((-\nabla_{\bm{x}_t}^2 f(\bm{x}_t) \bm{x}_t+\nabla_{\bm{x}_t} f(\bm{x}_t)) -\frac{1}{\beta \sigma_{0|t}^2}\bm{\mu}_{0|t}\right)+\bm{\mu}_{0|t}\right],
		\end{equation}
		where $\bm{H} = \nabla_{\bm{x}_t}^2 f(\bm{x}_t)+\frac{1}{\beta\sigma_{0|t}^2}\bm{I}$.
	\end{itemize}
	
	\textit{Proof. }

	Due to the implicit nature of $\hat{\bm{x}}_0(\bm{x}_t)$, directly computing $\nabla_{\bm{x}_t} f(\hat{\bm{x}}_0(\bm{x}_t))$ is intractable. Therefore, we consider an alternative approach by directly examining $\nabla_{\bm{x}_t} f(\hat{\bm{x}}_0(\bm{x}_t))$. By performing the differentiation $\nabla_{\bm{x}_t}$, we obtain
	\begin{equation}
		\begin{aligned}
			\nabla_{\bm{x}_t}\log \tilde{p}_t(\bm{x}_t|\rho) - \nabla_{\bm{x}_t}\log p_t(\bm{x}_t|\rho) &= \nabla_{\bm{x}_t}\log \int_{\bm{x}_0} p(\bm{x}_0|\bm{x}_t,\rho)B_{\beta}(\bm{x}_0)d\bm{x}_0 \\
			&= \frac{\int_{\bm{x}_0} \nabla_{\bm{x}_t}p(\bm{x}_0|\bm{x}_t,\rho)B_{\beta}(\bm{x}_0)d\bm{x}_0}{\int_{\bm{x}_0} p(\bm{x}_0|\bm{x}_t,\rho)B_{\beta}(\bm{x}_0)d\bm{x}_0}.
		\end{aligned}
	\end{equation}
	According to the assumption that $p_{t0}(\bm{x}_0|\bm{x}_t,\rho) = \mathcal{N}(\bm{x}_0|\bm{\mu}_{0|t},\sigma_{0|t}^2\bm{I})$, we have
	\begin{equation}
		\nabla_{\bm{x}_t}p(\bm{x}_0|\bm{x}_t,\rho) = \frac{\bm{x}_0-\bm{\mu }_{0|t}}{\sigma_{0|t}^2}p(\bm{x}_0|\bm{x}_t,\rho).
	\end{equation}
	Substituting into the above result, we have 
	\begin{equation}
		\begin{aligned}
			\nabla_{\bm{x}_t}\log \tilde{p}_t(\bm{x}_t|\rho) - \nabla_{\bm{x}_t}\log p_t(\bm{x}_t|\rho)=& \frac{\int_{\bm{x}_0} \frac{\bm{x}_0-\bm{\mu }_{0|t}}{\sigma_{0|t}^2}p(\bm{x}_0|\bm{x}_t,\rho)B_{\beta}(\bm{x}_0)d\bm{x}_0}{\int_{\bm{x}_0} p(\bm{x}_0|\bm{x}_t,\rho)B_{\beta}(\bm{x}_0)d\bm{x}_0} \\
			=& \frac{1}{\sigma_{0|t}^2} (\frac{\int_{\bm{x}_0} \bm{x}_0p(\bm{x}_0|\bm{x}_t,\rho)B_{\beta}(\bm{x}_0)d\bm{x}_0}{\int_{\bm{x}_0} p(\bm{x}_0|\bm{x}_t,\rho)B_{\beta}(\bm{x}_0)d\bm{x}_0}-\bm{\mu}_{0|t}) \\
			=& \frac{1}{\sigma_{0|t}^2} (\mathbb{E}\left[\tilde{\bm{x}}\right]-\bm{\mu}_{0|t}),
		\end{aligned}
	\end{equation}
	where $\tilde{\bm{x}}\sim p(\tilde{\bm{x}})\propto p(\bm{x}_0|\bm{x}_t,\rho)B_{\beta}(\bm{x}_0)$. Given an objective $f$ with the following quadratic form:
	\begin{equation}
		f(\bm{x}) = \frac{1}{2} \bm{x}^\top\bm{A}\bm{x} + \bm{b}^\top\bm{x},
	\end{equation}
	we have
	\begin{equation}
		B_{\beta}(\bm{x}_0)\propto e^{-\beta f(\bm{x})} = e^{-\beta (\frac{1}{2} \bm{x}^\top\bm{A}\bm{x} + \bm{b}^\top\bm{x})}.
	\end{equation}
	
	For $\beta\bm{A}+\frac{1}{\sigma_{0|t}^2}\bm{I}\succ \beta\bm{A}+\frac{1}{\sigma_{0}^2}\bm{I}\succ 0$, we have
	\begin{equation}
		\mathbb{E}\left[\tilde{\bm{x}}\right] = -\left(\beta\bm{A}+\frac{1}{\sigma_{0|t}^2}\bm{I}\right)^{-1}\left(\beta \bm{b}-\frac{1}{\sigma_{0|t}^2}\bm{\mu}_{0|t}\right),
	\end{equation}
	and then we have gradient guidance
	\begin{equation}
		\begin{aligned}
			\bm{G}_t &= \nabla_{\bm{x}_t}\log \tilde{p}_t(\bm{x}_t|\rho) - \nabla_{\bm{x}_t}\log p_t(\bm{x}_t|\rho) \\
			&=  -\frac{1}{\sigma_{0|t}^2}\left[(\beta\bm{A}+\frac{1}{\sigma_{0|t}^2}\bm{I})^{-1}(\beta \bm{b}-\frac{1}{\sigma_{0|t}^2}\bm{\mu}_{0|t})+\bm{\mu}_{0|t}\right].
		\end{aligned}
	\end{equation}
	
	For a general objective $f$, if we use the first-order Taylor expansion
	\begin{equation}
		f(\bm{x}) \approx f(\bm{x}_t) + \nabla_{\bm{x}_t} f(\bm{x}_t)^{\top} (\bm{x}-\bm{x}_t),
	\end{equation}
	then the Gradient Guidance can be formulated as the following form by setting $\bm{A} = \bm{0},\bm{b}=\nabla_{\bm{x}_t} f(\bm{x}_t)$:
	\begin{equation}
		\bm{G}_t = - \beta \nabla_{\bm{x}_t} f(\bm{x}_t).
	\end{equation}
	If we use the second-order Taylor expansion
	\begin{equation}
		f(\bm{x}) \approx f(\bm{x}_t) + \nabla_{\bm{x}_t} f(\bm{x}_t)^{\top} (\bm{x}-\bm{x}_t) + \frac{1}{2}(\bm{x}-\bm{x}_t)^{\top}\nabla_{\bm{x}_t}^2 f(\bm{x}_t) (\bm{x}-\bm{x}_t),
	\end{equation}
	then the Gradient Guidance can be formulated as the following form by setting $\bm{A} = \nabla_{\bm{x}_t}^2 f(\bm{x}_t),\bm{b}=-\nabla_{\bm{x}_t}^2 f(\bm{x}_t) \bm{x}_t+\nabla_{\bm{x}_t} f(\bm{x}_t))$:
	\begin{equation}
		\bm{G}_t = -\frac{1}{\sigma_{0|t}^2}\left[(\beta\nabla_{\bm{x}_t}^2 f(\bm{x}_t)+\frac{1}{\sigma_{0|t}^2}\bm{I})^{-1}\left[\beta (-\nabla_{\bm{x}_t}^2 f(\bm{x}_t) \bm{x}_t+\nabla_{\bm{x}_t} f(\bm{x}_t)) -\frac{1}{\sigma_{0|t}^2}\bm{\mu}_{0|t}\right]+\bm{\mu}_{0|t}\right].
	\end{equation}
	
	$\hfill\square$
	
	The posterior assumption in Corollary 1 can be satisfied easily. For example, with $p_{0}(\bm{x}_0|\rho)=\mathcal{N}(\bm{x}_0|\bm{\mu}_0,\sigma_0^2\bm{I})$ and forward process
	\begin{equation}
		d\bm{x}= -\theta\bm{x}dt+\sqrt{2\theta}d\bm{B}_t,
	\end{equation}
	we have 
	\begin{equation}
		p_{t}(\bm{x}_t|\rho) = \mathcal{N}\Big(\bm{x}_t | \bm{\mu}_0e^{-\theta t},(\sigma_0^2e^{-2\theta t}+1-e^{-2\theta t})\bm{I}\Big).
	\end{equation}
	Denote $\bm{\mu}_t=\bm{\mu}_0e^{-\theta t}, \sigma_t^2=\sigma_0^2e^{-2\theta t}+1-e^{-2\theta t}$, we have
	\begin{equation}
		p(\bm{x}_0|\bm{x}_t,\rho) = \mathcal{N}(\bm{x}_0|\bm{\mu}_{0|t},\sigma_{0|t}^2\bm{I}),
	\end{equation}
	where 
	\begin{equation}
		\begin{aligned}
			&\bm{\mu}_{0|t} = \bm{\mu}_{0} + \frac{\sigma_0^2}{\sigma_t^2} e^{-\theta t}(\bm{x}_t-\bm{\mu}_t), \\
			&\sigma_{0|t}^2 = \sigma_{0}^2\left(1-\frac{\sigma_{0}^2}{\sigma_{t}^2}e^{-2\theta t}\right). 
		\end{aligned}
	\end{equation}
	
	\subsection{Proof of Theorem 2}

	\textbf{Assumption 4.} For the forward process
	\begin{equation}
		d\bm{x}_t= \bm{a}(\bm{x}_t,t)dt+b(t) d\bm{B}_t,
	\end{equation}
	there is a constant $C$ such that
	\begin{itemize}
		\item [(i)] $\bm{a}(\bm{x}_t,t)$ is globally Lipschitz for any $t\in[0,T]$, i.e. $\|\bm{a}(\bm{x}_t,t)-\bm{a}(\bm{x}_t^{\prime},t)\|\leq C\|\bm{x}-\bm{x}_t^{\prime}\|$;
		
		\item [(ii)] $\bm{a}(\bm{x}_t,t)$ grows at most linearly for any $t\in[0,T]$, i.e. $\|\bm{a}(\bm{x}_t,t)\|\leq C(1+\|\bm{x}_t\|)$;
		
		\item [(iii)] $\bm{x}_t$ has a density $p_t\in\mathcal{C}^1$ for every $t>0$ and  
		\begin{equation}
			\int_{t_{0}}^{1}\int_{\|\bm{x}_t\|<R}|p_{t}(\bm{x}_t)|^{2}+\|\nabla_{\bm{x}_t}p_{t}(\bm{x}_t)\|^{2}\mathrm{d}x\mathrm{d}t<\infty,
		\end{equation}
		for any $R>0$ and $0 <t_0\leq T$;
		
		\item [(iv)] For each $S\in(0, T)$ and all $\|\bm{x}_t\|\leq N_R$ and $\|\bm{x}_t^{\prime}\|\leq N_R$, there is a constant $C_{S,N_R}$ such that $\nabla\log p_t(\bm{x}_t)$ is locally Lipschitz, i.e.,
		\begin{equation}
			\|\nabla\log p_t(\bm{x}_t)-\nabla\log p_t(\bm{x}_t^{\prime})\|\leq C_{S,N_R}\|\bm{x}_t-\bm{x}_t^{\prime}\|,
		\end{equation}
		for all $t\in(S, T)$.
	\end{itemize}

	\textbf{Remarks on Assumption 4.} Conditions (i)-(iii) are technical conditions on the forward SDE. They ensure that if we run a solution $p_t(\bm{x}_t)$ to the forward SDE, then $p_{T-t}(\bm{x}_{T_t})$ will be a solution to the reverse SDE. The last condition ensures that the solutions to the reverse SDE are unique. Assumption 4 can be expected to hold in practice, i.e., for any affine $\bm{a}(\cdot,t)$ and bounded data manifold.
	
	\textbf{Lemma 1 (Theorem 2 of \cite{pidstrigach2022score}).}
	Given a forward SDE with marginals $p_t(\bm{x}_t)$ and an approximated score $\bm{s}_{\bm{\theta}}(\bm{x}_t,t)$ to $\nabla \log p_t(\bm{x}_t)$, if the approximation error $\|\bm{s}_{\bm{\theta}}(\bm{x}_t,t)-\nabla \log p_t(\bm{x}_t)\|$ is bounded and Assumption 4 holds, then the marginal distribution of the reverse process using the approximated score starting from $p_T(\bm{x}_T)$ will have the same support as the data distribution $p_0(\bm{x}_0)$.

	\textbf{Theorem 2.} For any given $\rho\in(0,1)$, suppose that there exists a constant $\delta$ such that the error in the score estimation can be bounded as: 
	\begin{equation}
		\|\tilde{\bm{s}}_{\bm{\theta}}(\bm{x}_t,t,\rho) + \bm{G}_t - \nabla_{\bm{x}_t}\log \tilde{p}_t(\bm{x}_t|\rho)\| \leq \delta, \quad \forall~\bm{x}_t.
	\end{equation}
	For samples $\tilde{\bm{x}}_{sample}\sim p_{sample}(\bm{x}_0|\rho)$ generated by the reverse process
	\begin{equation}
		d\bm{x}_t=\left[\bm{a}(\bm{x}_t,t)-b(t)^2\big(\tilde{\bm{s}}_{\bm{\theta}}(\bm{x}_t,t,\rho)+\bm{G}_t \big)\right]dt+b(t)d\bm{\bar{B}}_t,
	\end{equation}
	with prior $p_{prior} = \mathcal{N}(\bm{0},\bm{I})$, affine drift coefficients $\bm{a}(\cdot,t)$, and
	\begin{equation}
		\tilde{\bm{s}}_{\bm{\theta}}(\bm{x}_t,t,\rho)= (1+w) \bm{s}_{\bm{\theta}}(\bm{x}_t,t,\rho) -w\bm{s}_{\bm{\theta}}(\bm{x}_t,t,\emptyset),
	\end{equation}
	as $T\rightarrow \infty$, $p_{sample}(\bm{x}_0|\rho)$ will have the same support as $\tilde{p}_{0}(\bm{x}_0|\rho)$. Further, as $\beta\rightarrow\infty$, $\tilde{\bm{x}}_{sample}$ will concentrate around $\bm{x}^* = \mathop{\arg\min}_{\bm{x}\in\mathcal{D}_{\rho}} f(\bm{x})$.
	
	\textit{Proof. } 
	
	For the forward process $d\bm{x}_t = \bm{a}(\bm{x}_t,t)dt+b(t) d\bm{B}_t, t\in [0,T]$ with affine drift coefficients $\bm{a}(\cdot,t)$, conditions (i)-(ii) in Assumption 4 are satisfied. For the given data set $\{\bm{x}^{(i)}\}_{i=1}^{N}$ contained in a ball of radius $M_R$, we have that $\log \tilde{p}_t(\bm{x}_t,t)\in \mathcal{C}^{\infty}$ in both $t$ and $\bm{x}_t$ for $t>0$ where the product distribution $\tilde{p}_0(\bm{x}_0|\rho)\propto p_{0}(\bm{x}_0|\rho)B_{\beta}(\bm{x}_0)$. Therefore we can integrate $\tilde{p}_t$ and its derivative over compact sets, implying that condition (iii) holds. Furthermore, for each $S\in(0, T)$, the Hessian w.r.t. $(\bm{x}_t, t)$ is continuous and obtains its maximum and minimum on the compact set $[S,T]\times B_{N_R}$, where $B_{N_R}$ is the ball of diameter $N_R$ around the origin. Therefore, the gradient $\nabla\log \tilde{p}_t(\bm{x}_t)$ is Lipschitz on $[S,T]\times B_{N_R}$, which proves condition (iv).
	
	The stationary distribution of the forward process is characterized by the corresponding Fokker-Planck equations, where $p_T = \mathcal{N}(\bm{0},\bm{I})$ when $T\rightarrow\infty$. Then we have that $p_{prior} = p_T$. Based on Lemma 1, if the score matching error is bounded, then the sampling distribution $p_{sample}(\bm{x}_0|\rho)$ with prior $p_{prior} = \mathcal{N}(\bm{0},\bm{I})$ will have the same support as the product distribution $\tilde{p}_0(\bm{x}_0|\rho)\propto p_{0}(\bm{x}_0|\rho)B_{\beta}(\bm{x}_0)$, where $B_{\beta}$ is the Boltzmann distribution $B_{\beta}(\bm{x}_0)\propto e^{-\beta f(\bm{x}_0)}$.
	
	Since $p_{0}(\bm{x}_0|\rho)$ has support $\mathcal{D}_{\rho}$ and the Boltzmann factor only changes the relative density within that domain, the support of $\tilde{p}_0(\bm{x}_0|\rho)$ also remains $\mathcal{D}_{\rho}$, i.e., 
	\begin{equation}
		\text{supp }p_{sample}(\bm{x}_0|\rho)=\text{supp } \tilde{p}_{0}(\bm{x}_0|\rho) = \mathcal{D}_{\rho}.
	\end{equation}
	
	As $\beta \rightarrow \infty$, sampling from the product distribution $\tilde{p}_0(\bm{x}_0|\rho)$ is equivalent to solving the optimization problem $\bm{x}^* = \mathop{\arg\min}_{\bm{x}\in\mathcal{D}_{\rho}} f(\bm{x})$. Then we have that as $T \rightarrow \infty$ and $\beta \rightarrow \infty$, the sample $\tilde{\bm{x}}_{sample}$ will concentrate around $\bm{x}^*$.
	
	$\hfill\square$

	Theorem 2 establishes that, by introducing an additional gradient guidance term into the reverse process, the sampling distribution of GGDOpt will attain the exact same support as the data distribution. Moreover, as the inverse temperature parameter $\beta$ increases, the sampling distribution becomes increasingly concentrated around points with the lowest function values within the support of the data distribution.
	
	The assumption in score estimation quantifies the approximation accuracy of the trained score network relative to the true score function. It depends on the training quality of the neural network and the expressiveness of the model class and this type of assumption is common in the theoretical analysis of diffusion models (see, e.g., \cite{pidstrigach2022score}, \cite{de2021diffusion}) and is used to establish convergence results in generative modeling and sampling.

	\subsection{Proof of Theorem 3}
	
	\textbf{Lemma 2 (\cite{bolley2005weighted}).} Let $\nu$ be a probability measure on $\mathbb{R}^d$. Assume that there exist $\bm{x}_0$ and a constant $\alpha>0$ such that
	\begin{equation}
		\int e^{\alpha\|\bm{x}-\bm{x}_0\|_2^2}d\nu(\bm{x})<\infty.
	\end{equation}
	Then for any probability measure $\mu$ on $\mathbb{R}^d$, it satisfies
	\begin{equation}
		\mathcal{W}_2(\mu,\nu)\leq C_\nu\big(\sqrt{D_{\mathrm{KL}}(\mu||\nu)}+\big(D_{\mathrm{KL}}(\mu||\nu)/2\big)^{1/4}\big),
	\end{equation}
	where $\mathcal{W}_2$ is the 2-Wasserstein distance and $C_\nu$ is defined as
	\begin{equation}
		C_{\nu}=\inf_{\bm{x}_{0}\in\mathbb{R}^{d},\alpha>0}\sqrt{\dfrac{1}{\alpha}\left(\dfrac{3}{2}+\log\int e^{\alpha\|\bm{x}-\bm{x}_{0}\|_{2}^{2}}d\nu(\bm{x})\right)}.
	\end{equation}
	
	\textbf{Lemma 3 (\cite{polyanskiy2016wasserstein}).} For any two probability density functions $\mu, \nu$ with bounded second moments, let $f:\mathbb{R}^d\to\mathbb{R}$ be a $C^1$ function such that
	\begin{equation}
		\|\nabla f(\bm{x})\|_2\leq C_1\|\bm{x}\|_2+C_2,\forall\bm{x}\in\mathbb{R}^d,
	\end{equation}
	for some constants $C_1,C_2\geq0$. Then
	\begin{equation}
		\left|\int_{\mathbb{R}^d}f(\bm{x})d\mu-\int_{\mathbb{R}^d}f(\bm{x})d\nu\right|\leq(C_1\sigma+C_2)\mathcal{W}_2(\mu,\nu),
	\end{equation}
	where $\mathcal{W}_2$ is the 2-Wasserstein distance and 
	\begin{equation}
		\sigma^2=\max\left\{\int_{\mathbb{R}^d}\|\bm{x}\|_2^2\mu(d\bm{x}),\int_{\mathbb{R}^d}\|\bm{x}\|_2^2\nu(d\bm{x})\right\}.
	\end{equation}
	
	\textbf{Lemma 4 (\cite{polyanskiy2016wasserstein}).} Let $p_{t}$ be the time $t$-marginal of a Brownian motion with initial distribution $\mu_{data}$. Denote by $c_i, i= 1, \ldots , d$ the eigenvalues of the covariance matrix $\text{Cov}(\mu_{data})$. Let $\mu_{prior}$ be the normal distribution with mean $m_T=\mathbb{E}[\mu_{data}]$ and covariance $C_T=\text{Cov}[\mu_{data}]+T\bm{I}$. Then
	\begin{equation}
		D_{KL}(p_T||\mu_{prior})\leq\frac12\log\left(\frac{\prod_{i=1}^d(c_i+T)}{T^d}\right).
	\end{equation}
	
	\textbf{Assumption 1.} We assume the following conditions hold:
	\begin{itemize}[leftmargin=10pt]
		\item The forward process is given by $d\bm{x}= b(t)d\bm{B}_t$;
		\item The reverse process starts in $p_{prior} = \mathcal{N}(\bm{m}_T, \bm{\Sigma}_T)$ where $\bm{m}_T = \mathbb{E}[\tilde{p}_0(\bm{x}_0|\rho)]$ and $\bm{\Sigma}_T = \text{Cov}(\tilde{p}_0(\bm{x}_0|\rho))+T\cdot \bm{I}$;
		\item The objective function $f(\bm{x})$ satisfies $\|\nabla_{\bm{x}}f(\bm{x}) \|_2\leq C_1\|\bm{x} \|_2+C_2$.
	\end{itemize}
	
	The first two conditions in Assumption 1 correspond to the Variance Exploding (VE) SDE in (\cite{songscore}) and are primarily used to characterize the discrepancy between the end distribution of the forward process and the prior distribution of the reverse process. Similar results can also be obtained for other forms of diffusion processes, e.g., Ornstein–Uhlenbeck processes. The third assumption imposes a growth bound on the gradient of the objective function. This type of regularity condition is common in the convergence analysis of stochastic optimization and sampling algorithms, particularly when studying stability and convergence under Langevin dynamics or diffusion-based methods (see, e.g., \cite{raginsky2017non}). In practice, this assumption holds for a broad class of functions, including smooth bounded functions and quadratic objectives, which frequently arise in real-world optimization problems.

	\textbf{Theorem 3.} Under Assumption 1, denote $\sigma^{(k)}, k=1,\ldots,n$, the eigenvalues of $\bm{\Sigma}_T$. For any given $\rho\in(0,1)$, denote $N_{\rho}=|\mathcal{D}_{\rho}|$ and $\bm{x}^* = \mathop{\arg\min}_{\bm{x}\in\mathcal{D}_{\rho}} f(\bm{x})$. Then for any given $T>0$ and $\beta>0$, the optimization error can be bounded by 
	\begin{equation}
		\begin{aligned}
			|\mathbb{E}[f(\tilde{\bm{x}}_t)]-f(\bm{x}^{*})|\leq \underbrace{C_I\big(\sqrt{C_T}+\big(C_T/2\big)^{1/4}\big)}_{I_1} + \underbrace{(N_\rho-1)\max_{\bm{x}\in\mathcal{D}_{\rho}}|f(\bm{x})-f(\bm{x}^*)|e^{-\beta \delta_{\rho}}}_{I_2},
		\end{aligned}
	\end{equation}
	where 
	\begin{equation}
		\begin{aligned}
			& C_I = \inf_{\bm{y}\in\mathbb{R}^{n},\alpha>0}\left\{\sqrt{\dfrac{1}{\alpha}\left(\dfrac{3}{2}+\log\int e^{\alpha\|\bm{x}-\bm{y}\|_{2}^{2}}\tilde{p}_0d\bm{x}\right)} (C_1\sigma_M+C_2)\right\},\\
			& \sigma_M = \max\left\{\int_{\mathbb{R}^n}\|\bm{x}\|_2^2\tilde{p}_0d\bm{x}, \int_{\mathbb{R}^n}\|\bm{x}\|_2^2p^{\pi}d\bm{x}\right\}, \\
			& C_T= \frac{1}{2}\log\left(\prod_{k=1}^{n}(\sigma^{(k)}/T)\right), \\
			& \delta_{\rho} = \min_{\bm{x}\in\mathcal{D}_{\rho},f(\bm{x})\neq f(\bm{x}^*)}|f(\bm{x})-f(\bm{x}^*)|.
		\end{aligned}
	\end{equation}
	
	\textit{Proof. } 
	
	Firstly, we give the form of $I_1$. By Lemma 4, we know that 
	\begin{equation}
		D_{KL}(\tilde{p}_0||p_{sample}) \leq D_{KL}(p_T || p_{prior}) \leq  \frac{1}{2}\log\left(\prod_{k=1}^{n}(\sigma^{(k)}/T)\right)=C_T.
	\end{equation}
	
	For $\tilde{p}_0(\bm{x}_0|\rho) \propto p_0(\bm{x}_0|\rho)B_{\beta}(\bm{x})$, there exist $\bm{y}$ and a constant $\alpha>0$ such that
	\begin{equation}
		\int e^{\alpha\|\bm{x}-\bm{y}\|_2^2}d\nu(\bm{x})<\infty.
	\end{equation}
	Then by Lemma 2, it satisfies
	\begin{equation}
		\begin{aligned}
			&\mathcal{W}_2(p_{sample},\tilde{p}_0(\bm{x}_0|\rho)) \\
			&\leq C_\nu\big(\sqrt{D_{\mathrm{KL}}(p_{sample}||\tilde{p}_0(\bm{x}_0|\rho))}+\big(D_{\mathrm{KL}}(p_{sample}||\tilde{p}_0(\bm{x}_0|\rho))/2\big)^{1/4}\big),
		\end{aligned}
	\end{equation}
	where $C_\nu$ is defined as
	\begin{equation}
		C_{\nu}=\inf_{\bm{y}\in\mathbb{R}^{d},\alpha>0}\sqrt{\dfrac{1}{\alpha}\left(\dfrac{3}{2}+\log\int\exp(\alpha\|\bm{x}-\bm{y}\|_{2}^{2})d\tilde{p}_0(\bm{x}|\rho)\right)}.
	\end{equation}
	By Lemma 3, we have that 
	\begin{equation}
		|\mathbb{E}[f(\tilde{\bm{x}}_t)]-\mathbb{E}[f(\bm{x}^{\pi})]| \leq C_I\big(\sqrt{C_T}+\big(C_T/2\big)^{1/4}\big).
	\end{equation}
	
	Next, we show the form of $I_2$. for $\bm{x}^{(i)}\in \text{supp } \tilde{p}_0(\bm{x}_0|\rho)$, the probability is given by 
	\begin{equation}
		p^{(i)} = \frac{e^{-\beta f(\bm{x}^{(i)})}}{\sum_{i=1}^{N_{\rho}}e^{-\beta f(\bm{x}^{(i)})}}.
	\end{equation}
	Denote $f^* = \min_{i=1}^{N_\rho} f(\bm{x}^{(i)})$ and $\text{Ind} = \{i\mid f(\bm{x}^{(i)}) = f^*\}$. Let $\delta^{(i)} = f(\bm{x}^{(i)}) - f^*\geq 0$, then the probability can be expressed as 
	\begin{equation}
		\begin{aligned}
			p^{(i)} = \frac{e^{-\beta f(\bm{x}^{(i)})}}{\sum_{i=1}^{N_{\rho}}e^{-\beta f(\bm{x}^{(i)})}} = \frac{e^{-\beta (f^*+\delta^{(i)})}}{\sum_{i=1}^{N_{\rho}}e^{-\beta (f^*+\delta^{(i)})}} = \frac{e^{-\beta \delta^{(i)}}}{\sum_{i=1}^{N_{\rho}}e^{-\beta \delta^{(i)}}}.
		\end{aligned}
	\end{equation}
	Then we have 
	\begin{equation}
		\mathbb{E}\left[f(\bm{x})\right] = \sum_{i=1}^{N_{\rho}} f(\bm{x}^{(i)})p^{(i)} = \sum_{i=1}^{N_{\rho}} (f^*+\delta^{(i)})\frac{e^{-\beta \delta^{(i)}}}{\sum_{i=1}^{N_{\rho}}e^{-\beta \delta^{(i)}}},
	\end{equation}
	and the limited inverse temperature error is given by
	\begin{equation}
		\begin{aligned}
			|\mathbb{E}\left[f(\bm{x})\right] - f^*| &= \left|\sum_{i=1}^{N_{\rho}} (f^*+\delta^{(i)})\frac{e^{-\beta \delta^{(i)}}}{\sum_{i=1}^{N_{\rho}}e^{-\beta \delta^{(i)}}}-\sum_{i=1}^{N_{\rho}} f^*\frac{e^{-\beta \delta^{(i)}}}{\sum_{i=1}^{N_{\rho}}e^{-\beta \delta^{(i)}}}\right| \\
			& = \sum_{i=1}^{N_{\rho}} \delta^{(i)}\frac{e^{-\beta \delta^{(i)}}}{\sum_{i=1}^{N_{\rho}}e^{-\beta \delta^{(i)}}}.
		\end{aligned}
	\end{equation} 
	Note that $\delta^{(i)}=0$ for $i\in\text{Ind}$, so we can simplify the sum as 
	\begin{equation}
		\sum_{i=1}^{N_{\rho}} \delta^{(i)}\frac{e^{-\beta \delta^{(i)}}}{\sum_{i=1}^{N_{\rho}}e^{-\beta \delta^{(i)}}} = \sum_{i=1,i\notin \text{Ind}}^{N_{\rho}} \delta^{(i)}\frac{e^{-\beta \delta^{(i)}}}{\sum_{i=1,i\notin \text{Ind}}^{N_{\rho}}e^{-\beta \delta^{(i)}}+\sum_{i=1,i\in \text{Ind}}^{N_{\rho}}e^{-\beta \delta^{(i)}}}.
	\end{equation}
	The denominator
	\begin{equation}
		\sum_{i=1,i\notin \text{Ind}}^{N_{\rho}}e^{-\beta \delta^{(i)}}+\sum_{i=1,i\in \text{Ind}}^{N_{\rho}}e^{-\beta \delta^{(i)}} = \sum_{i=1,i\notin \text{Ind}}^{N_{\rho}}e^{-\beta \delta^{(i)}}+ |\text{Ind}| \geq 1,
	\end{equation}
	so we have that
	\begin{equation}
		\begin{aligned}
			|\mathbb{E}\left[f(\bm{x})\right] - f^*| &= \sum_{i=1}^{N_{\rho}} \delta^{(i)}\frac{e^{-\beta \delta^{(i)}}}{\sum_{i=1}^{N_{\rho}}e^{-\beta \delta^{(i)}}}\\
			& \leq \sum_{i=1,i\notin \text{Ind}}^{N_{\rho}} \delta^{(i)}e^{-\beta \delta^{(i)}} \\
			& \leq (N_\rho-1)\max_{\bm{x}\in\mathcal{D}_{\rho}}|f(\bm{x})-f(\bm{x}^*)|e^{-\beta \delta_{\rho}},
		\end{aligned}
	\end{equation} 
	where 
	\begin{equation}
		\delta_{\rho} = \min_{\bm{x}\in\mathcal{D}_{\rho},f(\bm{x})\neq f(\bm{x}^*)}|f(\bm{x})-f(\bm{x}^*)|.
	\end{equation}
	
	Then the optimization error can be bounded by 
	\begin{equation}
		\begin{aligned}
			|\mathbb{E}[f(\tilde{\bm{x}}_t)]-f(\bm{x}^{*})|\leq \underbrace{C_I\big(\sqrt{C_T}+\big(C_T/2\big)^{1/4}\big)}_{I_1} + \underbrace{(N_\rho-1)\max_{\bm{x}\in\mathcal{D}_{\rho}}|f(\bm{x})-f(\bm{x}^*)|e^{-\beta \delta_{\rho}}}_{I_2},
		\end{aligned}
	\end{equation}
	where 
	\begin{equation}
		\begin{aligned}
			& C_I = \inf_{\bm{y}\in\mathbb{R}^{n},\alpha>0}\left\{\sqrt{\dfrac{1}{\alpha}\left(\dfrac{3}{2}+\log\int e^{\alpha\|\bm{x}-\bm{y}\|_{2}^{2}}\tilde{p}_0d\bm{x}\right)} (C_1\sigma_M+C_2)\right\},\\
			& \sigma_M = \max\left\{\int_{\mathbb{R}^n}\|\bm{x}\|_2^2\tilde{p}_0d\bm{x}, \int_{\mathbb{R}^n}\|\bm{x}\|_2^2p^{\pi}d\bm{x}\right\}, \\
			& C_T= \frac{1}{2}\log\left(\prod_{k=1}^{n}(\sigma^{(k)}/T)\right), \\
			& \delta_{\rho} = \min_{\bm{x}\in\mathcal{D}_{\rho},f(\bm{x})\neq f(\bm{x}^*)}|f(\bm{x})-f(\bm{x}^*)|.
		\end{aligned}
	\end{equation}
	
	$\hfill\square$

	Theorem 3 establishes that, in practical settings, the optimization error of the sampling process can be decomposed and bounded by two components: the limited time length error $I_1$ and the limited inverse temperature error $I_2$, which are given as follows:
	\begin{equation}
		\begin{aligned}
			|\mathbb{E}[f(\tilde{\bm{x}}_t)]-f(\bm{x}^{*})|&\leq|\underbrace{\mathbb{E}[f(\tilde{\bm{x}}_t)]-\mathbb{E}[f(\bm{x}^\pi)]|}_{I_1} + \underbrace{|\mathbb{E}[f(\bm{x}^\pi)]-f(\bm{x}^*)|}_{I_2}.
		\end{aligned}
	\end{equation}
	As a direct corollary, under mild assumptions, GGDOpt is shown to generate asymptotically optimal solutions to problem (\ref{aeq1}) as the time length $T$ and inverse temperature $\beta$ increase.

\end{document}